\documentclass[12pt,leqno]{book}
\usepackage{latexsym,amsfonts}
\usepackage{amsmath,amssymb,amsthm}
\usepackage[notref, notcite]{}
\usepackage[colorlinks,linkcolor=red,citecolor=blue,urlcolor=blue]{hyperref}
\usepackage[active]{srcltx} 
\usepackage{emptypage}

\newtheorem{theorem}{Theorem}[chapter]
\newtheorem{proposition}{Proposition}[chapter]
\newtheorem{lemma}{Lemma}[chapter]
\newtheorem{definition}{Definition}[chapter]
\newtheorem{example}{Example}[chapter]
\newtheorem{remark}{Remark}[chapter]
\newtheorem{corollary}{Corollary}[chapter]
\newtheorem{exercise}{Exercise}[chapter]
\newtheorem{note}{Note}
\newenvironment{Note}{\begin{note}\quad\\\rm}{\end{note}\medskip}

\def\B{{\mathbb B}}
\def\Z{{\mathbb Z}}
\def\Bb{{\cal B}}
\def\Gg{{\cal G}}
\def\Xx{{\cal X}}
\def\Aa{{\cal A}}
\def\Pp{{\cal P}}
\def\Qq {{\cal Q}}
\def\O{{\cal O}}
\def\L{{\cal L}}
\def\uX {\underline X}
\def\xv{x_1,\dots,x_n}
\def\Xv{X_1,\dots,X_n}
\def\var {{\rm Var}}

\def\hta {{\hat\theta}}
\def\Th {{\Theta}}
\def\la{\lambda}
\def\si{\sigma}
\def\tt{{\tilde \theta}}
\def\varf{\varphi}
\def\Cov{{\rm Cov}}
\def\dun {\dot\nu}
\def\dQ {\dot\Q}
\def\dP {\dot\Pp}
\def\dQ {\dot\Qq}
\def\dq {\dot q}
\def\halfe {{1\over 2}}
\def\vier {{1\over 4}}
\def\been {{\bf 1}}
\def\otee {\int_0^1}
\def\odP {\overline{\dP^0}}
\def\del{\dot\ell}
\def\tell{{\tilde \ell}}
\def\tel {\tilde \ell}
\def\stel {\ell^\star}
\def\nn {{\sum\limits_{i=1}^n}}
\def\equa{{{1\over{\sqrt n}}}}
\def\limy {\mathop{\rm\lim}\limits_{n\to\infty}}
\def\D{{\cal D}}
\def\sqn {{\sqrt n}}
\def\Nn{{\cal N}}
\def\al {\alpha}
\def\be {\beta}
\def\ga {\gamma}

\def\la {\lambda}
\def\La {\Lambda}
\def\si {\sigma}

\def\Del {\Delta}

\def\vier {{1\over 4}}

\def\zes {{1\over 6}}

\def\acht {{1\over 8}}

\def\twaalf {{1\over 12}}

\def\ene {{1\over n}}

\def\ex{{\rm E\, }}
\def\expar{{\rm E }}
\def\vars {{\varsigma}}
\def\vart {{\vartheta}}
\def\tW{\tilde W}
\def\si {\sigma}
\def\Ff{{\cal F}}
\def\half{\scriptstyle{1\over 2}}

\begin{document}

\frontmatter 

\bibliographystyle{plain}

\frontmatter

\title{\bf Semiparametric Statistics \\
and \\
the Spread Inequality \\
\quad \\
\quad \\
Lecture Notes\\
\vspace{8cm} }
\author{Chris A.J. Klaassen \\
{\normalsize Korteweg-de Vries Institute for Mathematics}\\
{\normalsize University of Amsterdam}}
\maketitle

\chapter{Preface}
Classical mathematical statistics deals with models that are parametrized
by a Euclidean, i.e. finite dimensional, parameter. Quite often such
models have been and still are chosen in practical situations for their
mathematical simplicity and tractability. However, these models are
typically inappropriate since the implied distributional assumptions
cannot be supported by hard evidence. It is natural then to relax these
assumptions. This leads to the class of semiparametric models.

An example is the classical linear regression model with normal
error distribution. If the normality assumption is replaced by the
often less questionable assumption that the errors have a density
with mean zero, we have a semiparametric linear regression model
with the regression parameters as Euclidean parameter and the
unknown error density as so called Banach parameter.

Semiparametrics has been initiated by Stein (1956), who presented
a result on Fisher information matrices and who discussed its
consequences for the asymptotic theory of estimation and of
testing hypotheses. Apparently, H\'ajek (1962) has been unaware of
this paper when he presented the first semiparametrically
efficient test. However, the main development in semiparametrics
has been in estimation theory. This development has
started with Van Eeden (1970), Beran (1974), and Stone (1975) and
got a strong impetus by Pfanzagl and Wefelmeyer (1982), Bickel
(1982), and Begun, Hall, Huang, and Wellner (1983). A
comprehensive account is given in Bickel, Klaassen, Ritov, and
Wellner (1993, 1998).

This course will treat some highlights from the theory of semiparametric
estimation as it developed during the last quarter of the past
century. More recent results will be discussed as well. Some topics are
crucial for a proper understanding of the issues of semiparametrics and
they will be treated. Others are less essential. From this last group a
few have been chosen according to our own, very personal biases.

Semiparametrics has been studied in a local asymptotic setting, in which
the Convolution Theorem yields bounds on the performance of regular
estimators. Alternatively, local asymptotics can be based on the Local
Asymptotic Minimax Theorem and on the Local Asymptotic Spread Theorem,
both valid for any sequence of estimators. This Local Asymptotic Spread
Theorem is a straightforward consequence of a Finite Sample Spread
Inequality, which has some intrinsic value for estimation theory in
general. We will discuss both the Finite Sample and Local Asymptotic
Spread Theorem, as well as the Convolution Theorem.

These notes will {\sl not} constitute a self-contained text. Often
reference will be made to the original literature for relevant technical
details. However, the main line of the argument will be understandable
without consulting the literature referred to.

Preliminary notes by Edwin van den Heuvel of this course have been transformed into an intermediate version by Bert van Es. 
I would like to thank them for their great help in putting these notes and their contents into shape. 
Furthermore, I would like to thank many people for fruitful and inspiring discussions. 
Most of all, my thesis advisor Willem van Zwet, who has put me on the
track towards the spread inequality many years ago, and Peter
Bickel, Yanki Ritov, and Jon Wellner, with whom it has been very
pleasant to collaborate over the many years that have passed since
the start of the writing of our book in 1983. \\
\medskip
\hspace{30em} {\em Chris A.J. Klaassen}

\mainmatter

\tableofcontents

\chapter{Introduction}\label{chap:1}
The random quantity $X$ takes values in the measurable space $(\Xx
, \Aa)$. It has {\sl unknown} distribution $P$. This distribution
belongs to the {\sl known} class $\Pp$. The map
$$\nu:\Pp\to{\mathbb R}^m$$
defines a Euclidean parameter. We will study estimation of the
unknown value of $\nu(P)$. Any measurable map
$$t: \Xx\to{\mathbb R}^m$$
defines an estimator $T=t(X)$ of $\nu(P)$. We will focus on
optimality and efficiency of estimators $T$ of $\nu(P)$.

To prove that an estimator is optimal, one has to show that no estimator
is better. Typically, this is done by proving the validity of a bound on
the performance of estimators and by showing subsequently that the
proposed estimator attains this bound. This is exactly what we will do
here. We will discuss such bounds in Chapter \ref{chap:2} through
\ref{chap:7}, and in Chapter \ref{chap:8} we will study construction of
estimators attaining these bounds.

We will start by a discussion of the so-called spread inequality in
Chapter \ref{chap:2}.
This inequality yields a bound on the distribution of an estimator
averaged over the class $\Pp$ of underlying distributions $P$.
It is formulated in terms of quantile functions, and it is valid without
any restriction on the estimator.

In most (but not all) of the models we will consider, the
observations $\xv$ may be viewed as realizations of i.i.d. random
variables $\Xv$. In the notation used above, this means that $X$
is an $n$-vector $(\Xv)$. Then $X$ has distribution $P^n$ for some
unknown distribution $P$ from some known class $\Pp$. Note that
this is an abuse of notation since $P^n$, and $\Pp^n$ have been
denoted by $P$ and $\Pp$ above. In fact, we will abuse notation
even to a larger extent. However, it will always be clear when we
are in the general situation with
\begin{equation}\label{generalmodel}
X, \ (\Xx,\Aa), \ P\in\Pp, \ \nu:\Pp\to{\mathbb R}^m, \ T=t(X),
t:\Xx\to{\mathbb R}^m,
\end{equation}
and when we are in the i.i.d. situation with
\begin{eqnarray}\label{iidmodel}
X,\Xv \mbox{ independent and identically distributed on }
(\Xx,\Aa), \\ X\sim P\in\Pp, \ \nu:\Pp\to{\mathbb R}^m, \
T_n=t_n(\Xv), \ t_n:\Xx^n\to{\mathbb R}^m. \nonumber
\end{eqnarray}

The spread inequality from Chapter \ref{chap:2} is a finite sample
inequality in the sense that, in the i.i.d. situation, it is valid
for $n$ finite. In the remainder of the course we will focus on
asymptotics with $n\to\infty$. First, we will study regular
parametric models. In these models $\Pp$ is parametrized by
$\theta\in\Theta$, an open subset of ${\mathbb R}^k$; so
$\Pp=\{P_\theta:\theta\in\Theta\}$. In Chapter \ref{chap:3} it
will be shown that regular parametric models are Locally
Asymptotically Normal (LAN).

Under Local Asymptotic Normality the spread inequality yields an
asymptotic bound on the performance of any sequence $\{T_n\}_{n\in
{\mathbb N}}$ of estimators. The resulting Local Asymptotic Spread
Theorem will be discussed in Chapter \ref{chap:4}. Some geometric
interpretations of it for estimation in the presence of nuisance
parameters will be presented in Chapter \ref{chap:5}.

In a Semiparametric model  the class $\Pp$ of distributions cannot be
parame\-trized as a regular  parametric model, but it contains regular
parametric submodels in a natural way.
The Local Asymptotic Spread Theorem yields lower bounds for the asymptotic
performance of an estimator sequence within each of these regular parametric
submodels.
The supremum of these lower bounds thus is an obvious lower bound to the
asymptotic performance of estimators in the semiparametric model.
This approach is the key idea in semiparametrics,  and it will be
discussed in Chapter \ref{chap:6}, together with some of its consequences.
Traditionally, the theory of efficient estimation in
semiparametric models is developed via the Convolution Theorem,
which is studied in Chapter \ref{chap:7}.

Construction of estimators attaining the semiparametric asymptotic bound
obtained in Chapters \ref{chap:6}  and \ref{chap:7}, will be studied in
Chapter \ref{chap:8}. These estimators are called (asymptotically)
efficient. In Chapter \ref{chap:9} we will extend this theory for i.i.d.
models to a large class of time series models. In Chapter
\ref{chap:10} we will discuss estimation of a Banach parameter
$\nu(P)$ where $\nu$ maps $\Pp$ into a Banach space $\Bb$. Finally, in
Chapter \ref{chap:11} we will study cross sectional sampling, a
time-saving method in survival analysis.

We close this chapter by presenting three examples of semiparametric
models for i.i.d. situations with $X, \ \Xv$ i.i.d. $X\sim P\in\Pp$,
and one example of a non-i.i.d. model.

\begin{example}[Symmetric Location]\label{exam:1.1}
{\rm Let $\Gg$ be the class of distributions on ${\mathbb R}$ with a density with respect to Lebesgue measure that is
symmetric about $0$, and let $\varepsilon$ be a random variable with
unknown distribution $G\in\Gg$. Consider
\begin{equation}\label{1.1}
X=\nu+\varepsilon
\end{equation}
with $\nu$ the unknown real location parameter of interest and $G$
the unknown nuisance parameter. Then $X$ has distribution
$P_{\nu,G}$, which  is symmetric about $\nu$, and we could
parametrize $\Pp$ as $\Pp=\{P_{\nu,G}:\nu\in{\mathbb
R},G\in\Gg\}$. Stein (1956) claimed that under this semiparametric
model $\Pp$ it should be possible to estimate the location
parameter $\nu$ asymptotically as well at any $G_0$ as under the
parametric model $\{P_{\nu,G_0}:\nu\in{\mathbb R}\}$. Beran (1974)
and Stone (1975) were the first to construct such estimators,
which were called adaptive, because they adapt to the underlying
distribution $G$. Van Eeden (1970) constructed an adaptive
estimator earlier, but under the additional assumption that $G$ is
strongly unimodal, that is, $G$ has a log-concave density.
Therefore, it is fair to say that the symmetric location model
triggered the development of semiparametrics. \hfill $\Box$}
\end{example}

\begin{example}[Linear Regression]\label{exam:1.2}
{\rm Here we take $\Gg$ to be the class of all distributions on
${\mathbb R}$ with a density and mean $0$. We extend (\ref{1.1}) to
\begin{equation}\label{1.2}
Y=\nu^T Z+\varepsilon,
\end{equation}
where $\varepsilon$ and $Z$ are independent, $\varepsilon$ has unknown
distribution $G\in\Gg$ and $Z\in{\mathbb R}^m$ has unknown
distribution such that the covariance matrix of $Z,$
\begin{equation}\label{1.3}
\Sigma_Z = \ex(Z-\ex Z)(Z -\ex Z)^T,
\end{equation}
is nonsingular. If one observes i.i.d. copies of $X=(Y,Z)$, then
it  is possible to estimate $\nu\in{\mathbb R}^m$ adaptively, as
we will see later. \hfill $\Box$ }
\end{example}

\begin{example}[Cox's Proportional Hazards Model]\label{exam:1.3}
{\rm Again, one obser\-ves i.i.d. copies of $X=(Y,Z)$ with $Z$ an
$m$-vector of covariates. Given $Z=z$ the one-dimensional survival
time $Y$ has hazard function
\begin{equation}\label{1.4}
\la(y\mid z,\nu)=\exp(\nu^T z)\la(y), \quad y>0.
\end{equation}
Here $\la(y)=g(y)/(1-G(y))$ is the unknown baseline hazard
function with $G$ a distribution function with density $g$. We
will be interested in estimation of $\nu$ in the presence of the
nuisance parameter $\la$. \hfill $\Box$ }
\end{example}

\begin{example}[ARMA]\label{exam:1.4}
{\rm
Consider the ARMA$(p,q)$ process where we  observe $X_t, \ t=1,\dots,
n$, generated by
\begin{equation}\label{1.5}
X_t=\rho_1X_{t-1}+\cdots +\rho_pX_{t-p}+\varf_1\varepsilon_{t-1}+\cdots
+\varf_q\varepsilon_{t-q}+\varepsilon_t,
\end{equation}
{\em and} the square integrable starting values
$X_{1-p},\dots,X_0, \ \varepsilon_{1-q},\dots,\varepsilon_0$. The
innovations $\varepsilon_{1-q},\dots,\varepsilon_n$ are i.i.d. with unknown
density $g$, mean $0$, and finite variance. The ${\mathbb
R}^{p+q}$-valued parameter $\nu=(\rho_1,\dots,\rho_p, \
\varf_1,\dots,\varf_q)$ is restricted in such a way that the
zeroes of $1-\rho_1B\cdots-\rho_pB^p$ and of
$1+\varf_1B\cdots+\varf_qB^q$ all lie outside the unit circle in
${\mathbb C}$. This restriction guarantees stationarity and
invertibility of the process (1.5); cf. Box, Jenkins, and Reinsel
(1994). \hfill $\Box$ }
\end{example}

\section{Exercises Chapter {\ref{chap:1}}}

\begin{exercise}[Identifiability]\label{exer:1.1}
{\rm Determine the map $\nu:\Pp\to{\mathbb R}^m$ that identifies
the Euclidean parameter $\nu \in {\mathbb R}^m$ in Example
{\ref{exam:1.1}} $(m=1)$ and Example {\ref{exam:1.2}}.
When such a map $\nu$ exists, the corresponding parameter $\nu$ is
called identifiable. \hfill $\Box$ }
\end{exercise}

\begin{exercise}[Identifiability of Cox's Model]\label{exer:1.2}
{\rm Show that the Cox Proportional Hazards Model from Example
\ref{exam:1.3} is identifiable if the covariance matrix of $Z$ is
nonsingular, by noting
\begin{equation}
\log\left(\frac{\lambda(y\,|\,z,\nu)}{\lambda(y\,|\,EZ,\nu)}\right)
= \nu^T(z - EZ).
\end{equation}
\hfill $\Box$ }
\end{exercise}

\chapter{Spread Inequality}\label{chap:2}

Consider a class $\Pp$ of probability measures dominated by a
$\si$-finite measure $\mu$ on $(\Xx,\Aa)$. The map $\nu: \
\Pp\to{\mathbb R}$ defines our one-dimensional parameter of
interest. One observes a realization $x$ of the random variable
$X$ with unknown distribution $P\in\Pp$. As in the general case of
Chapter~\ref{chap:1}, $T=t(X)$ is an estimator of $\nu(P)$.

\begin{example}[Degenerate Estimators]\label{exam:2.1}
{\rm Fix $P_0\in\Pp$. Note that the degenerate estimator
$T=\nu(P_0)$ cannot be improved in estimating $\nu(P)$ {\em at}
$P=P_0$, but it is extremely bad, of course, at all $P$ with
$\nu(P)\ne\nu(P_0)$. \hfill$\Box$ }
\end{example}

This little example illustrates that the performance of an estimator can
be judged only by considering its behavior at several (or all)  $P\in\Pp$
simultaneously.
In (most versions of) the Cram\'er-Rao
inequality this is forced by the condition of unbiasedness on the estimator.
In Bayesian statistics one chooses a prior.
The H\'ajek-Le Cam convolution theorem restricts attention to regular
estimators, as we will see in Chapter \ref{chap:7}.
The so-called Local Asymptotic Minimax theorem considers suprema over
subclasses of $\Pp$.

Here we will choose a weight function $\tW$ on $\Pp$, which is a
{\bf probability measure that describes the relative stress we want to
put on the performance of the estimator} $T$ {\bf at the respective} $P\in\Pp$.
(One may view $\tW$ as a Bayesian prior, but $\tW$ is not meant to describe prior knowledge here.)
Thus we will study the average distribution of
$T-\nu(P)$ defined by
\begin{equation}\label{2.1}
G(y)=\int\limits_\Pp P(T-\nu(P)\le y)d{\tW}(P), \quad y\in{\mathbb
R}.
\end{equation}
Roughly speaking, the spread inequality says that this distribution $G$
cannot be concentrated arbitrarily much.
In fact, it states that, whatever the estimator $T$, the quantiles of
$G$ are at least as far apart as the corresponding quantiles of a
particular distribution function $K$, so
\begin{equation}\label{2.2}
\Big|G^{-1}(v)-G^{-1}(u)\Big|\ge\Big|K^{-1}(v)-K^{-1}(u)\Big|, \quad
u,v\in (0,1).
\end{equation}
Here $K$ is a distribution defined in terms of the model and {\em not}
the estimator.

Define the random vector $(X,\vart)=(X,\nu(P))$, where $P$ is
random with distribution $\tW$. Note that our estimation problem
is completely described by the joint distribution of $(X,\vart)$,
which generates a parametric model. Furthermore, note that
(\ref{2.1}) may be rewritten as
\begin{equation}\label{2.3}
G(y)=P(T-\vart\le y), \quad y\in{\mathbb R}.
\end{equation}
We will denote the conditional density of $X$ at $x$ given
$\vart=\theta$ by $p(x\mid\theta)$ with respect to $\mu,$ the
density of $\vart$ at $\theta$ by $w(\theta)$ with respect to
Lebesgue measure, and the joint density of $(X,\vart)$ by
$f(x,\theta)=p(x\mid\theta)w(\theta)$ with respect to the product
measure of $\mu$ and Lebesgue measure.

Differentiating $G$ we see that under regularity conditions, $G$ has
density
\begin{eqnarray}
\lefteqn{g(y)=\lim_{\varepsilon\to
0}\varepsilon^{-1}(1-G(y-\varepsilon)-(1-G(y)))}\nonumber\\
&=&\lim_{\varepsilon\to 0}\ex\bigg({\been}_{[T-\vart>y]}\varepsilon^{-1}
\Big\{{f(X,\vart+\varepsilon)\over f(X,\vart)}-1\Big\}\bigg)\label{2.4}\\
&=&\ex{\been}_{[T-\vart>y]}S, \quad y\in{\mathbb R},\nonumber
\end{eqnarray}
where the {\bf Score Statistic} $S$ is defined as
\begin{equation}\label{2.5}
S={{\dot f}\over f}(X,\vart)={{\dot p}\over p} (X\mid\vart)+{{\dot
w}\over w}(\vart)
\end{equation}
with $\dot{ }$ denoting differentiation with respect to $\theta$.
Now, let $H$ be the distribution function of $S$,
\begin{equation}\label{2.6}
H(z)=P(S\le z), \quad z\in{\mathbb R}.
\end{equation}
Then (\ref{2.4}) yields (see also Exercise {\ref{exer:2.1}})
\begin{eqnarray}
g(G^{-1}(s))&=&\ex(\ex({\been}_{[T-\vart>G^{-1}(s)]}\mid S)S)\nonumber\\
&=&\otee \ex({\been}_{[T-\vart>G^{-1}(s)]}\mid
S=H^{-1}(t))H^{-1}(t)dt\label{2.7}\\
&=&\otee \varf(t)H^{-1}(t)dt.\nonumber
\end{eqnarray}
We note
\begin{equation}\label{2.8}
1-s=\otee \ex({\been}_{[T-\vart>G^{-1}(s)]}\mid
S=H^{-1}(t))dt=\otee\varf(t)dt
\end{equation}
and $0\le\varf(t)\le 1$. Maximizing the right-hand side of
(\ref{2.7}) over all critical functions $\varf$ satisfying
(\ref{2.8}) we get
\begin{equation}\label{2.9}
g(G^{-1}(s))\le\int\limits_s^1 H^{-1}(t)dt
\end{equation}
in view of the monotonicity of $H^{-1}$.
Integration yields
\begin{eqnarray}
\lefteqn{{G^{-1}(v)-G^{-1}(u)\ge\int\limits_u^v {1\over
g(G^{-1}(s))}\,ds\ge \int\limits_u^v {1\over \int\limits_s^1
H^{-1}(t)dt}ds}}\label{2.10}\\
&=&K^{-1}(v)-K^{-1}(u), \quad 0<u<v<1,\nonumber
\end{eqnarray}
with
\begin{equation}\label{2.11}
K^{-1}(u)=\int\limits_{\half}^u {1\over\int\limits_s^1
H^{-1}(t)dt}ds.
\end{equation}
All in all, we have
\begin{theorem} [One-dimensional Spread Inequality]\label{thm:2.1}
The distribution function $G$ defined in {\rm(\ref{2.1})} and
{\rm(\ref{2.3})} is at least as spread out as the distribution function
$K$ defined by its quantile function
$K^{-1}$ in {\rm(\ref{2.11})}; in formula (cf. {\rm(\ref{2.2})})
\begin{equation}\label{2.12}
G^{-1}(v)-G^{-1}(u)\ge K^{-1}(v)-K^{-1}(u), \quad 0<u<v<1.
\end{equation}
\end{theorem}
Regularity conditions, details of the proof and many implications may be
found in Klaassen (1989a,b) and Van den Heuvel and Klaassen (1997).

An alternative way to see (\ref{2.9}) is to write
\begin{eqnarray}
\lefteqn{\int_s^1H^{-1}(t)dt - g(G^{-1}(s))}\nonumber\\
&=&
\int_0^1 H^{-1}(t)\{\been_{[t>s]}-\ex(\been_{[T-\vart>G^{-1}(s)]}|S=H^{-1}(t))\}dt\nonumber\\
&=& \int_0^1
\{H^{-1}(t)-H^{-1}(s)\}\{\been_{[t>s]}-\ex(\been_{[T-\vart>G^{-1}(s)]}|S=H^{-1}(t))\}dt
\label{new2.13}\\
&=&
\int_0^1 \{H^{-1}(t)-H^{-1}(s)\}\{\been_{[H^{-1}(t)>H^{-1}(s)]}-
\ex(\been_{[T-\vart>G^{-1}(s)]}|S=H^{-1}(t))\}dt\nonumber\\
&=&
\ex \{S-H^{-1}(s)\}\{\been_{[S>H^{-1}(s)]}-\been_{[T-\vart>G^{-1}(s)]}\}\nonumber\\
&=& \ex
|S-H^{-1}(s)||\been_{[S>H^{-1}(s)]}-\been_{[T-\vart>G^{-1}(s)]}|\geq
0.\nonumber
\end{eqnarray}
Assume that $S$ has no pointmasses and note that $T-\vart$ has
neither since it has a density $g$ with respect to Lebesgue
measure. Then, equality can hold in (\ref{new2.13}) only if
\begin{equation}\label{new2.14}
\ex |\been_{[H(S)>s]}-\been_{[G(T-\vart)>s]}|=0.
\end{equation}
Consequently, the spread inequality from (\ref{2.12}) is an
equality only if
\begin{equation}\label{new2.15}
0=\ex
\int_0^1|\been_{[H(S)>s]}-\been_{[G(T-\vart)>s]}|ds=\ex|H(S)-G(T-\vart)|.
\end{equation}
We have proved the "only if" part of (See note \ref{note:2.1} in Appendix \ref{noteschap2} for the
"if" part)
\begin{theorem} [Spread Equality]\label{thm:2.2}
Let regularity conditions hold such that {\rm Theorem
\ref{thm:2.1}} is valid. $S$ has no pointmasses and equality holds
in {\rm(\ref{2.12})} iff
\begin{equation}\label{new2.16}
T-\vart=G^{-1}(H(S)),\ \mbox{a.s.}
\end{equation}
\end{theorem}

Often a model, i.e. the class $\Pp$, contains natural parametric submodels
that are of interest, or it is parametric itself. To study this
situation we consider the parametric model
\begin{equation}\label{new2.17}
\Pp=\{P_\theta:\theta\in \Theta\},\quad\Theta\subset{\mathbb
R}^k,\ \Theta\ \mbox{open},
\end{equation}
and the map $\nu:\Pp\to{\mathbb R}^m$ that defines our
$m$-dimensional parameter of interest.  Together with the
parametrization $\theta\mapsto P_\theta$ this map $\nu$ defines
the map $q:\Theta\to {\mathbb R}^m$ by $q(\theta)=\nu(P_\theta),\
\theta\in\Theta$. With $T=t(X)$ an estimator of $\nu(P)$ we are
interested in  the behavior of
\begin{equation}\label{new2.18}
T-q(\theta),\ \theta\in \Theta.
\end{equation}
The weight function $\tilde W$ on $\Pp$ can be represented now by
a density $w$ with respect to Lebesgue measure on ${\mathbb R}^k$,
such that
\begin{equation}\label{new2.19}
\int_\Theta w(\theta)d\theta =1.
\end{equation}
Let $\vart$ be a random vector with density $w$ and denote the
conditional density of $X$ given $\vart =\theta$ by
$p(x\,|\,\theta)$ and the joint density of $X$ and $\vart$ by
$f(x,\theta)=p(x\,|\, \theta)w(\theta)$. We will derive a spread
inequality for
\begin{equation}\label{new2.20}
G(y)=P(T-q(\vart)\leq y),\quad y\in {\mathbb R}^m,
\end{equation}
but we will take $m=1$ first. The crucial differences with
(\ref{2.3}) and Theorem \ref{thm:2.1} are the dimension $k$
of $\Theta$ that might be more than 1, and the map $q$ that need
not be the identity.

Again, the basic step is differentiation of $G$ as in (\ref{2.4}).
Extend $q:\Theta\to{\mathbb R}$ to $q:{\mathbb R}^k\to{\mathbb R}$.
Fix $b\in {\mathbb R}^k$ and assume that for all $\eta\in{\mathbb R}^k$ and all $\varepsilon \in {\mathbb R}$ with
$|\varepsilon|$ sufficiently small the equation
\begin{equation}\label{new2.21}
q(\eta) + \varepsilon = q(\eta+\delta b)
\end{equation}
has a solution $\delta=\delta(\varepsilon,\eta)$. In case of multiple
solutions we will take the one closest to 0. Given $q$ this
restricts the choice of $b\in {\mathbb R}^k$. Then, by the
substitution $\theta=\eta+\delta(\varepsilon,\eta)b$ we obtain (cf.
Exercise \ref{exer:2.8})
\begin{eqnarray}
\lefteqn{1-G(y-\varepsilon)}\nonumber\\
&=&
P(T-q(\vart)+\varepsilon>y)\label{new2.22}\\
&=& \int_{{\mathbb R}^k}\int_{\cal
X}\been_{[t(x)-q(\eta)>y]}f(x,\eta+\delta(\varepsilon,\eta)b)
\{1+b^T\dot\delta(\varepsilon,\eta)\}d\mu(x)d\eta,\nonumber
\end{eqnarray}
where $\dot\delta$ denotes the column $k$-vector of derivatives of
$\delta(\varepsilon,\eta)$ with respect to $\eta$. From
(\ref{new2.21}) we get $\delta(0,\eta)=0$ and by differentiation
with respect to $\eta$, $\dot\delta(0,\eta)=0$. Differentiating
(\ref{new2.21}) with respect to $\varepsilon$ at $\varepsilon=0$ we
obtain
\begin{equation}\label{new2.23}
{\partial\over\partial \varepsilon}\
\delta(\varepsilon,\eta)\Big|_{\varepsilon=0} = {1\over{b^T\dot
q(\eta)}},
\end{equation}
where again $\dot q$ denotes the column $k$-vector of derivatives.
Assuming that $q$ is even twice continuously differentiable we
denote the $k\times k$-matrix of second mixed partial derivatives
of $q$ by $\ddot q$. Differentiation of (\ref{new2.21}) with
respect to $\eta$ and $\varepsilon$ at $\varepsilon=0$ yields the column
$k$-vector
\begin{equation}\label{new2.24}
{\partial\over\partial \varepsilon}\
\dot\delta(\varepsilon,\eta)\Big|_{\varepsilon=0} =-\  {\ddot q(\eta)
b\over{(b^T\dot q(\eta))^2}}.
\end{equation}

Combining (\ref{new2.22}) through (\ref{new2.24}) we get under
regularity conditions that $G$ has density
\begin{eqnarray}
\lefteqn{g(y)=\lim_{\varepsilon\to 0}
\varepsilon^{-1}\Big(1-G(y-\varepsilon)-(1-G(y))\Big)}
\nonumber\\
&=& \lim_{\varepsilon\to 0} \ex\Big(\been_{[T-q(\vart)>y]}
\varepsilon^{-1} \Big\{{{f(X,\vart +\delta(\varepsilon,\vart)b)
\{1+b^T\dot\delta(\varepsilon,\vart)\}}\over{f(X,\vart)}}-1\Big\}\Big)
\label{new2.25}\\
&=& \ex \been_{[T-q(\vart)>y]}S\nonumber
\end{eqnarray}
with

\noindent
\parbox{1cm}{\begin{equation}\label{new2.26}\end{equation}}
\parbox{12cm}{
\begin{eqnarray*}
\lefteqn{S={\partial\over{\partial\varepsilon}} \log\Big(f(X,\vart
+\delta(\varepsilon,\vart)b)
\{1+b^T\dot\delta(\varepsilon,\vart)\}\Big)\Big|_{\varepsilon=0}}\\
&=& {b^T\over{b^T\dot q(\vart)}}\Big\{{\dot p\over p}(X|\vart)
+{{\dot w}\over w}(\vart) - {{\ddot q(\vart)b}\over{b^T\dot
q(\vart)}}\Big\}.
\end{eqnarray*}
}
\hfill

\noindent Repeating the arguments (\ref{2.6}) through
(\ref{new2.16}) with $S$ from (\ref{new2.26}) we arrive at the
spread inequality (\ref{2.12}) with $G$ the distribution function
of $T-q(\vart)$ with $q$ one dimensional. Moreover, if this $S$
has no pointmass equality holds iff
\begin{equation}\label{new2.27}
T-q(\vart)=G^{-1}(H(S)),\ \mbox{a.s.}
\end{equation}

Note  the dependence of $S$ on the quite arbitrarily chosen $b\in
{\mathbb R}^k$. Every appropriate $b$ yields a lower bound $K$ and
typically, these bounds differ.

For general $m$-dimensional functions $q$ we apply the preceding
spread inequality  to
$$
a^T(T-q(\vart)),\quad a\in {\mathbb R}^m,
$$
for every $a\in {\mathbb R}^m$. The resulting multidimensional
spread inequality is

\begin{theorem} [General Spread Inequality]\label{thm:2.3}
In the parametric model $\Pp$ from {\rm(\ref{new2.17})} with
weight density $w$ satisfying {\rm(\ref{new2.19})}, let
$f(x,\theta)=p(x|\theta)w(\theta)$ be the density of $(X,\vart)$
with respect to $\mu\times$Lebesgue measure on $({\cal
X}\times{\mathbb R}^k,{\cal A}\times{\cal B}^k),\ {\cal B}^k$
Borel. Let $\nu(P_\theta)=q(\theta)$ be the parameter of interest
with $q:{\mathbb R}^k\to{\mathbb R}^m$ twice continuously
differentiable. Fix $a\in {\mathbb R}^m,\, a \neq 0,$ and let
there exist $b\in {\mathbb R}^k$ such that
\begin{equation}\label{new2.21a}
a^Tq(\eta+\delta b)-\varepsilon = a^Tq(\eta)
\end{equation}
(cf. {\rm(\ref{new2.21})}) has a solution
$\delta=\delta(\varepsilon,\eta)$ for all $\eta\in {\mathbb R}^k$ and
$\varepsilon$ sufficiently small with $\delta(\varepsilon,\eta)$ as
close to 0 as possible, and such that $\varepsilon\mapsto f(x,
\eta+\delta(\varepsilon,\eta)b)\{1+b^T\dot \delta(\varepsilon,\eta)\}$
is absolutely continuous in a neighborhood of 0 for almost all $x$
and $\eta$. Define
\begin{equation}\label{new2.28}
S_{a,b}={b^T\over{b^T\dot q^T(\vart)a}}
\Big\{
{\dot p\over p}(X|\vart) + {\dot w\over w}(\vart) -
{{\sum_{h=1}^ma_h\ddot q_h(\vart)b}\over{b^T\dot q^T(\vart)a}}
\Big\}
\end{equation}
with $\dot q(\theta)$ the $m\times k$-matrix with
${\partial\over{\partial\theta_j}}\, q_i(\theta)$ in the $i$-th
row and $j$-th column and with $\ddot q_h(\theta)$ the $k\times
k$-matrix with ${\partial^2\over
{\partial\theta_i\partial\theta_j}} q_h(\theta)$ in the $i$-th row
and $j$-th column. Assume that $S_{a,b}$ is well-defined with
\begin{equation}\label{new2.29}
\ex |S_{a,b}|<\infty.
\end{equation}
Let $T=t(X),\ t:{\cal X}\to{\mathbb R}^m$ measurable, be any
estimator of $q(\theta)$. The distribution $G_a$ of
$a^T(T-q(\vart))$ is at least as spread out as $K_{a,b}$ in the
sense of {\rm(\ref{2.12})} with $K_{a,b}$ defined as in
{\rm(\ref{2.11})} with $H$ replaced by the distribution $H_{a,b}$
of $S_{a,b}$. Furthermore,
\begin{equation}\label{new2.30}
G_a=K_{a,b}
\end{equation}
holds and $S_{a,b}$ has no pointmass iff
\begin{equation}\label{new2.31}
a^T(T-q(\vart))=G_a^{-1}(H_{a,b}(S_{a,b}))\ \mbox{a.s.}
\end{equation}
\end{theorem}

Note that the distribution of $T-q(\vart)$ is determined, once for
each $a\in {\mathbb R}^m$ the distribution of $a^T(T-q(\vart))$ is
given.

\begin{remark}
{\rm
In the proof of Theorem \ref{thm:2.3}, the absolute continuity
and the related existence of $\dot p/p$ and $\dot w/w$ are used in
checking the validity of (\ref{new2.25}) in a way similar to the
proof of Theorem \ref{thm:2.1}. In fact, one needs

\noindent
\parbox{2cm}{\begin{equation}\label{new2.32}\end{equation}}
\parbox{12cm}{
\begin{eqnarray*}
\lefteqn{\int_{\cal X}\int_{{\mathbb R}^k}
\Big|f(x,\theta+\delta(\varepsilon,\theta)b)
\{1+b^T\dot \delta(\varepsilon,\theta)\}-f(x,\theta)  }\\
&&\quad
  -b^T\{\dot p(x\mid\theta)w(\theta) + p(x\mid\theta)\dot w(\theta)\}
\, {\partial\over {\partial \varepsilon}}\,
\delta(\varepsilon,\theta)\Big|_{\varepsilon
=0}\Big|d\mu(x)d\theta=o(\varepsilon).
\end{eqnarray*}
}
\hfill

\noindent Consequently, it also suffices to assume this
$\L_1$-differentiability itself.
}
\end{remark}

\section{Exercises Chapter {\ref{chap:2}}}

\begin{exercise}[Uniformity Generates It All]\label{exer:2.1}
{\rm Let $U$ have a uniform distribution on the unit interval and let $F$
be an arbitrary distribution function. Define the inverse distribution
function by
\begin{equation}
F^{-1}(u)= \inf \{x\, :\, F(x) \geq u\}
\end{equation}
and prove that $F^{-1}(U)$ has distribution function $F$.

\hfill $\Box$ }
\end{exercise}

\begin{exercise}[Implications Spread Ordering]\label{exer:2.2}
{\rm The class of distribution functions on the real line may be ordered
partially by the following spread order. The distribution function $G$ is
at least as spread out as $F$ if all quantiles of $G$ are at least as far
apart as the corresponding quantiles of $F$, more precisely if
\begin{equation}\label{2.12.1}
G^{-1}(v)-G^{-1}(u)\ge F^{-1}(v)-F^{-1}(u), \quad 0<u<v<1,
\end{equation}
which we will denote by
\begin{equation}\label{2.12.2}
G \ge_s F.
\end{equation}
Prove that the variance of $G$ equals at least the variance of $F$ for all
$F$ and $G$ satisfying (\ref{2.12.2}).

\hfill $\Box$ }
\end{exercise}

\begin{exercise}[Estimator Attaining the Normal Spread Lower Bound]\label{exer:2.2b}
{\rm Let $\vartheta$ have a normal distribution on ${\mathbb R}$
with mean $\mu$ and variance $\tau^2$. Given $\vartheta = \theta,$
let $\Xv$ be i.i.d. with a one-dimensional normal distribution
with mean $\theta$ and variance $\sigma^2.$ Compute the spread
lower bound for estimators $T=t(X_1, \dots, X_n)$ of $\theta,$ and
determine the estimator that attains this bound.

\hfill $\Box$ }
\end{exercise}

\begin{exercise}[Location Model]\label{exer:2.3}
{\rm Let $\Xv$ be i.i.d. with density $f(\cdot-\theta)$. Let
$T_n=t_n(\Xv)$ be a translation equivariant estimator of the
location parameter $\theta$, i.e.
\begin{equation}\label{2.13}
t_n(x_1+a,\cdots,x_n+a)=t_n(\xv)+a.
\end{equation}
With $G(\cdot)$ the distribution function from (\ref{2.3}) we
define $G_n(y)=G(y/ \sqrt n),\ y \in \mathbb R.$  Note
\begin{equation}\label{2.14}
G_n(y)=P_{f(\cdot-\theta)}({\sqrt n}(T_n-\theta)\le y), \quad {\rm
all} \ \theta\in{\mathbb R}.
\end{equation}
Choose $\vart$ normal with mean $0$ and variance $\si^2$, and let
$f$ be the standard normal density $\varf$.
Compute the lower bound $K$ from (\ref{2.11}).
This lower bound depends on $\si$, but is valid for all $\si>0$.
Take the limit for $\si\to\infty$.
Which estimator attains the resulting bound?

\hfill $\Box$ }
\end{exercise}

\begin{exercise}[More Spread Out Than Uniform]\label{exer:2.4}
{\rm
Let $\ex|S|<\infty$ and note that (\ref{2.4}) implies
\begin{equation}\label{2.15}
g(G^{-1}(s)) = \ex ({\been}_{[T-\vart>G^{-1}(s)]}S ) \le \ex|S|.
\end{equation}
Prove that $G$ is at least as spread out as the uniform distribution on
$[0,\, 1/\ex|S|]$.

\hfill $\Box$ }
\end{exercise}

\begin{exercise}[{\bf Van Zwet Inequality}]\label{exer:2.5}
{\rm
With $\ex S^2<\infty$, show that (\ref{2.4}) implies
\begin{equation}\label{2.16}
g(G^{-1}(s))\le (\ex S^2\{s\land(1-s)\})^{1/2}.
\end{equation}
Prove that $G$ is at least as spread out as the symmetric triangular
distribution with support $[-\sqrt{2/\ex S^2}, \ \sqrt{2/\ex S^2}]$.

\hfill$\Box$ }
\end{exercise}

\begin{exercise}[Trigonometric Spread Inequality]\label{exer:2.6}
{\rm With $\ex S^2<\infty$, show that (\ref{2.4}) implies
\begin{equation}\label{2.16}
g(G^{-1}(s))\le (\ex S^2\{s(1-s)\})^{1/2},
\end{equation}
which sharpens the Van Zwet inequality. Prove that $G$ is at least
as spread out as the trigonometric distribution with distribution
function $(1+\sin(\sqrt{\ex S^2}\,x))/2\,,\, |x|\leq \pi/(2
\sqrt{\ex S^2})$.

\hfill$\Box$ }
\end{exercise}

\begin{exercise}[Strong Unimodality]\label{exer:2.7}
{\rm A density $f$ with respect to Lebesgue measure on the real line is
called strongly unimodal if the convolution of any unimodal density with
$f$ yields a unimodal density. Ibragimov (1956) has shown that $f$ is
strongly unimodal iff $f$ has a log-concave version, i.e. iff the
logarithm of an appropriate version of $f$ is concave. The distribution
function of a strongly unimodal density is called strongly unimodal as
well. Prove that the lower bound $K$ of the spread inequality (\ref{2.12})
is strongly unimodal.

\hfill$\Box$ }
\end{exercise}

\begin{exercise}[Jacobian]\label{exer:2.8}
{\rm In (\ref{new2.22}) the transformation $\eta \mapsto
\eta+\delta(\varepsilon/c,\eta)b$ has been applied. Prove that the
corresponding Jacobian $|J_k+b\dot{\delta}^T(\varepsilon/c,\eta)|$
equals $1+b^T \dot{\delta}(\varepsilon/c,\eta)$, where $J_k$ denotes
the $k$x$k$ identity matrix and the $\dot{ }$ denotes
differentiation with respect to $\eta$.

\hfill$\Box$ }
\end{exercise}

\begin{exercise}[Multivariate Normal Spread Lower Bound]\label{exer:2.9}
{\rm Let $\vartheta$ have a multivariate normal distribution on
${\mathbb R}^k$ with mean vector $\mu$ and as covariance matrix
$\Pi$. Let the conditional distribution of $X$ given
$\vartheta=\theta$ be multivariate normal with mean vector
$\theta$ and covariance matrix $\Sigma$. Take $q(\theta)=\theta$
and compute the score statistic $S_{a,b}$ from Theorem
{\ref{thm:2.3}} for arbitrary $a$ and $b$ in ${\mathbb R}^k$.
Subsequently, given $a$ determine $b$ such that the variance of
the lower bound $K_{a,b}$ in the General Spread Inequality is
maximal. Call the corresponding lower bound $K_a$. Does there
exist a random vector $Z$ such that for all $a \in {\mathbb R}^k$
the random variable $a^T Z$ has distribution function $K_a$?

\hfill$\Box$ }
\end{exercise}

\begin{exercise}[Implication Spread Equality]\label{exer:2.2}
{\rm If there exists an estimator that attains equality in the
spread inequality, then the score statistic $S$ from (\ref{2.5})
cannot have atoms.

\hfill $\Box$ }
\end{exercise}

\chapter{Regular Parametric Models}\label{chap:3}

In this chapter we will consider the situation with $X, \ \Xv$
i.i.d. with unknown distribution from the class
$\Pp=\{P_\theta:\,\theta\in\Theta\}, \ \Theta\subset{\mathbb
R}^k$. This is called a parametric model, because $\Theta$ is
finite dimensional. We will assume the existence of a $\si$-finite
measure $\mu$ dominating $\Pp$, and we will represent the elements
of $\Pp$ in $\L_2(\mu)$ by the square roots
$s(\theta)=p^{1/2}(\theta)$ of their densities
$p(\theta)=dP_\theta/d\mu$ with respect to $\mu$.

\begin{definition}[Regularity]\label{defn:3.1}
{\rm The parametrization $\theta \mapsto P_\theta$ is regular if
\begin{itemize}
\item[(i)] $\Theta$ is an open subset of ${\mathbb R}^k$,
\item [(ii)] for every $\theta_0\in\Theta$ there exists a $k$-vector
${\dot\ell}(\theta_0)$ of score functions in $\L_2(P_{\theta_0})$
such that\par
\begin{equation}\label{3.1}
s(\theta)=s(\theta_0)+\halfe(\theta-\theta_0)^T{\dot\ell}(\theta_0)s(\theta_0)
+o(|\theta-\theta_0|)
\end{equation}
\item [] in $\L_2(\mu)$ as $|\theta-\theta_0|\to 0$,
\item[(iii)] for every $\theta_0\in\Theta$ the $k\times k$ Fisher
information matrix\par
\begin{equation}\label{3.2}
I(\theta_0)=\int
{\dot\ell}(\theta_0){\dot\ell}^T(\theta_0)dP_{\theta_0}
\end{equation}
is nonsingular,
\item [(iv)] the map $\theta\mapsto{\dot\ell}(\theta)s(\theta)$ is
continuous from $\Theta$ to $\L^k_2(\mu)$.
\end{itemize}
}
\end{definition}

Typically, the Fr\'echet-differentiability property of (\ref{3.1}) is
verified via pointwise differentiability.
A result to this extent is the following

\begin{proposition}\label{prop:3.1}
Let $\Theta$ be open and let for all $\theta\in\Theta$
\begin{itemize}
\item [(i)] $p(x;\theta)=p(\theta)(x)$ be continuously differentiable in
$\theta$ for $\mu$-almost all $x$ with gradient ${\dot
p}(x;\theta)$,
\item [(ii)] $|{\dot\ell}(\theta)|\in\L_2(P_\theta)$ with \par
\begin{equation}\label{3.3}
{\dot\ell}(\theta)={\dot
p}(\theta)/p(\theta){\been}_{[p(\theta)>0]},
\end{equation}
\item [(iii)] $I(\theta)$ defined by {\rm(\ref{3.2})} with ${\dot\ell}(\theta)$
as in {\rm(\ref{3.3})}, is nonsingular and continuous in $\theta$.
\end{itemize}
\medskip
Then, with ${\dot\ell}(\theta)$ as in {\rm(\ref{3.3})}, the
parametrization
$$\theta\to P_\theta$$
is regular.
\end{proposition}

\noindent{\bf Proof  of Proposition \ref{prop:3.1}.}

Fix $\theta_0\in\Theta$. In view of (i) we have for $\mu$-almost
all $x$
\begin{eqnarray}
\lefteqn{s(x;\theta)-s(x;\theta_0)}\nonumber\\
&=&\otee\halfe(\theta-\theta_0)^T
\del(x;\theta_0+\la(\theta-\theta_0))p^{1/2}(x;\theta_0+
\la(\theta-\theta_0))d\la\label{a}
\end{eqnarray}
provided $|\theta-\theta_0|$ is sufficiently small. The continuity
of ${\dot p}(\theta)$ and (\ref{a}) imply
\begin{equation}\label{b}
(s(x;\theta)-s(x;\theta_0)){\been}_{[p(x;\theta_0)>0]}
=\halfe(\theta-\theta_0)^T{\del}(x;\theta_0)s(x;\theta_0)+
o(|\theta-\theta_0|).
\end{equation}
By (\ref{a}), (ii) and (iii) it follows that (see notes \ref{note:1},
\ref{note:2} and \ref{note:3} in the appendix)
\begin{eqnarray}
\lefteqn{\int|s(x;\theta)-s(x;\theta_0)|^2d\mu(x)} \nonumber \\
&\le&\vier(\theta-\theta_0)^T\otee\int\del\del^Tp(x;\theta_0+
\la(\theta-\theta_0))d\mu(x)d\la \, (\theta-\theta_0) \nonumber \\
&=&\vier(\theta-\theta_0)^T\otee
I(\theta_0+\la(\theta-\theta_0))d\la \, (\theta-\theta_0) \label{c}\\
&=&\vier(\theta-\theta_0)^T
I(\theta_0)(\theta-\theta_0)+o(|\theta-\theta_0|^2)\nonumber\\
&=&\vier(\theta-\theta_0)^T\int\limits_{p(\theta_0)>0}
\del\del^Tp(\theta_0)d\mu \ (\theta-\theta_0)+o(|\theta-\theta_0|^2).\nonumber
\end{eqnarray}
Without loss of generality let
$(\theta-\theta_0)/|\theta-\theta_0|$ converge (see note
\ref{note:4}). Applying Vitali's theorem to (\ref{b}) and
(\ref{c}) we obtain
\begin{equation}\label{d}
\int\limits_{p(x;\theta_0)>0}
|s(x;\theta)-s(x;\theta_0)-\halfe(\theta- \theta_0)^T\del
s(x;\theta_0)|^2d\mu(x)=o(|\theta-\theta_0|^2),
\end{equation}
which combined with (\ref{c}) implies
\begin{equation}\label{e}
\int\limits_{p(x;\theta_0)=0}|s(x;\theta)-s(x;\theta_0)|^2d\mu(x)=
o(|\theta-\theta_0|^2).
\end{equation}

Together, (\ref{d}) and (\ref{e}) imply the Fr\'echet differentiability
(\ref{3.1}).

 From (i) and (iii) we obtain
\begin{equation}\label{f}
\lim_{\theta\to\theta_0}\del_i(x;\theta)s(x;\theta)=\del_i(x;\theta_0)s
(x;\theta_0)
\end{equation}
for $\mu$-almost every $x$ with $p(x;\theta_0)>0$, and
\begin{eqnarray}
\lefteqn{\limsup_{\theta\to\theta_0}\int\limits_{p(x;\theta_0)>0}\del_i^2(x
;\theta) p(x;\theta)d\mu(x)}\nonumber\\
&\leq&\limsup_{\theta\to\theta_0}\int\del_i^2(x;\theta)p(x;\theta)d\mu(x)
\label{g}\\
&=&\int\del_i^2p(x;\theta_0)d\mu(x)=\int\limits_{p(x;\theta_0)>0}
\del_i^2p(x;\theta_0)d\mu(x).\nonumber
\end{eqnarray}
The continuity from (iv) of Definition 3.1 follows from (\ref{f})
and (\ref{g}) by another application of Vitali's theorem. \hfill
$\Box$

\begin{example}[Normal distribution]\label{exam:normal}{\rm
Let us show that the family of normal distributions is a regular parametric
family.
Here the densities are given by
\begin{equation}
p_{(\mu,\sigma)}(x) = {1\over{\sqrt{2\pi}\sigma}}e^{-{1\over 2}{(x-\mu)^2\over
\sigma^2}},\quad -\infty<x<\infty.
\end{equation}
Note that the parameter $\theta$ is equal to the vector
$(\theta_1,\theta_2)= (\mu,\sigma)$ and that the parameter space
is given by $\Theta= {\mathbb R}\times (0,\infty)$. We will check
the conditions of Proposition \ref{prop:3.1}.

First we see that the density is continuously differentiable in
$\theta$ for all $x$. The gradient is equal to
\begin{equation}
\dot p(x;(\mu,\sigma)) =
\left(\begin{array}{l}
{1\over{\sqrt{2\pi}\sigma^3}}(x-\mu)e^{-{1\over 2}{(x-\mu)^2\over \sigma^2}}\\
{1\over{\sqrt{2\pi}\sigma^2}}\Big({(x-\mu)^2\over \sigma^2}-1\Big)
e^{-{1\over 2}{(x-\mu)^2\over \sigma^2}}
\end{array}
\right).
\end{equation}
As score function we get
\begin{equation}\label{3.13}
\dot \ell(x;(\mu,\sigma)) =
{{\dot p(x;(\mu,\sigma))}\over { p(x;(\mu,\sigma))}}=
\left(\begin{array}{l}
{(x-\mu)\over \sigma^2}\\
{1\over\sigma}\Big({(x-\mu)^2\over \sigma^2}-1\Big)
\end{array}
\right).
\end{equation}
Computing the Fisher information (\ref{3.2}) we get
\begin{equation}
I((\mu,\sigma)) =
\left(\begin{array}{cc}
1\over\sigma^2&0\\
0&{2\over\sigma^2}
\end{array}
\right).
\end{equation}
Since all the conditions of Proposition \ref{prop:3.1} hold this family is
a regular parametric model.
}
\end{example}

The proposition cannot be applied to the family of uniform
distributions on $(0,\theta)$ and to the  Laplace$(\mu)$
distributions. In the first case this is essential. A Local
Asymptotic Spread Theorem as we will prove it, is valid for
regular parametric models, but it does not hold for uniform
distributions. For the Laplace distributions regularity can be
shown by adapting the proof of the proposition.

The notation
\begin{equation}\label{3.4}
\ell(\theta)=\log p(\theta)=2\log s(\theta)
\end{equation}
explains the notation ${\dot\ell}(\theta)$ in (\ref{3.1}). We
define the log-likelihood and the score function of $(\Xv)$ by
\begin{equation}\label{3.5}
L_n(\theta)=\nn\ell (X_i;\theta)\quad {\rm and} \quad
S_n(\theta)=\equa\nn\del(X_i;\theta)
\end{equation}
respectively. A fundamental consequence of regularity is {\sl Local
Asymptotic Normality}, formulated in

\begin{theorem}[LAN]\label{thm:3.1}
Write
\begin{equation}\label{3.6}
L_n(\theta+{t\over{\sqrt n}})-L_n(\theta)=t^TS_n(\theta)- \halfe
t^T I(\theta)t+R_n(\theta,t).
\end{equation}
If $\Pp=\{P_\theta:\, \theta\in\Theta\}$ is a regular parametric
model, then for any compact $K\subset\Theta$, any
$M\in(0,\infty)$, and any $\varepsilon>0$
\begin{equation}\label{3.7}
\limy \,\sup_{|t|\le M}\,\sup_{\theta\in K}
P_\theta(|R_n(\theta,t)|>\varepsilon)=0
\end{equation}
holds.
Moreover
\begin{equation}\label{3.8}
S_n(\theta)\ {\buildrel\D\over\rightarrow_{\theta+{t\over\sqn}}} \
\Nn(I(\theta)t,\, I(\theta))
\end{equation}
uniformly in $\theta\in K$ and in $|t|\le M$, for every compact
$K\subset\Theta$ and for every finite $M\ge 0$. Finally,
\begin{equation}\label{3.9}
\limy\,\sup_{|t|\le M}\,\sup_{\theta\in K}
P_\theta(|S_n(\theta+{t\over{\sqrt
n}})-S_n(\theta)+I(\theta)t|>\varepsilon)=0,
\end{equation}
for all $K, \ M$, and $\varepsilon$ as above.
\end{theorem}

By the uniform convergence of (\ref{3.8}) we mean that for any
compact $K\subset\Theta$, any $M\in(0,\infty)$ and for every
bounded, continuous function $g$ on ${\mathbb R}$ we have
\begin{equation}\label{uncon}
\limy\,\sup_{|t|\le M}\,\sup_{\theta\in K} |\expar_{\theta +
{t\over \sqn}}\ g(S_n(\theta))-\expar_{\Nn(I(\theta)t,I(\theta))}\
g(S)|=0.
\end{equation}
For $\theta$ and $t$ fixed (\ref{uncon}) is equivalent to
(\ref{3.8}) itself, see for instance Grimmett and Stirzaker (1992)
Theorem 7.2.19.

\bigskip

\noindent{\bf Proof of Theorem \ref{thm:3.1}}

We will only sketch the main ideas.
A precise proof may be found on pages 509--513 of Appendix A.9 of Bickel,
Klaassen, Ritov and Wellner (1993) (henceforth BKRW).
Define
\begin{equation}\label{3.1.a}
T_{ni}=2\{s(X_i;\, \theta+{t\over {\sqrt n}})/s(X_i;\theta)-1\},
\quad i=1,\dots,n.
\end{equation}
First, one shows that $P_\theta(\max_{1\le i\le
n}|T_{ni}|\ge\varepsilon)\to 0$ as $n\to\infty$, uniformly in
$\theta\in K$ and $|t|\le M$ (see note \ref{note:6}). Then, on the
event $\{\max_{1\le i\le n}|T_{ni}|<\varepsilon\}$, one uses
\begin{equation}\label{3.1.b}
\log(1+x)=\int_0^x(1-y+{y^2\over 1+y})dy
=x-\halfe x^2+{2\over 3}\al|x|^3
\end{equation}
with $|\al|\le 1$ for $|x|\le\halfe$, to write
\begin{eqnarray}
\lefteqn{L_n(\theta+{t\over{\sqrt
n}})-L_n(\theta)=2\nn\log(1+\halfe
T_{ni})}\nonumber\\
&=&2\nn\Big\{\halfe T_{ni}-\acht
T^2_{ni}+\twaalf\al_{ni}\Big|T_{ni}\Big|^3\Big\}\nonumber\\
&=&\nn T_{ni}-\vier\nn T^2_{ni}+\zes\nn\al_{ni}\Big|T_{ni}\Big|^3\cr
&=&t^TS_n(\theta)-\halfe t^TI(\theta)t\label{3.1.c}\\
&&\quad +\Big\{\nn T_{ni}-\Big(t^TS_n(\theta)-\vier
t^TI(\theta)t\Big)\Big\}\nonumber\\
&&\quad -\vier\Big\{\nn
T^2_{ni}-t^TI(\theta)t\Big\}+\zes\nn\al_{ni}\Big|T_{ni}\Big|^3\nonumber
\end{eqnarray}
with $|\al_{ni}|\le 1$ a.s. Subsequently, one proves that the last
three terms in (\ref{3.1.c}) all converge to $0$ in probability
under $\theta$ uniformly in $\theta\in K, \ |t|\le M$(on the event
$\{\max_{1\le i\le n}|T_{ni}|<\varepsilon\}$) (see note \ref{note:6}).

One needs a uniform central limit theorem to prove (\ref{3.8}) for
$M=0$ (see BKRW, p.513). To prove (\ref{3.8}) in its full
generality we need (\ref{3.9}). To prove ({\ref{3.9}) we first
have to study contiguity. \hfill $\Box$

\begin{definition}[Contiguity]\label{cont}
{\rm The sequence of probability measures $\{Q_n\}$ is contiguous with
respect to
$\{P_n\}$ (both $P_n$ and $Q_n$ on $(\Xx_n,\Aa_n))$, if for all
$\{A_n\}$ with $P_n(A_n)\to 0$ also $Q_n(A_n)\to 0$ holds.
}
\end{definition}

\begin{corollary}\label{cor:3.2}
The Local Asymptotic Normality as formulated in {\rm(\ref{3.7})}
and {\rm(\ref{3.8})} with $M=0$ and hence $t=0$, implies that
$\{P^n_{\theta_n}\}$ and $\{P^n_{\theta_n+t_n/{\sqrt n}}\}$ are
mutually contiguous for all convergent sequences $\{\theta_n\}$
and $\{t_n\}$.
\end{corollary}

\noindent{\bf Proof}

This follows from Le Cam's first and third lemma; see p.17 of
BKRW. \hfill$\Box$

\bigskip

We will use Corollary \ref{cor:3.2} to complete the proof of Theorem
\ref{thm:3.1}.

\bigskip

\newpage

\noindent{\bf Proof of Theorem \ref{thm:3.1}.  Continued}

Contiguity as in Corollary \ref{cor:3.2}, continuity and symmetry
of $I(\theta)$ and (\ref{3.7}) show
\begin{eqnarray}
\lefteqn{0=L_n(\theta+{t\over{\sqrt
n}})-L_n(\theta)+L_n\Big(\theta+{t+h\over{\sqrt
n}}\Big)-L_n\Big(\theta+{t\over{\sqrt n}}\Big)}\nonumber\\
&&\quad -\Big\{L_n\Big(\theta+{t+h\over{\sqrt
n}}\Big)-L_n(\theta)\Big\}
\nonumber\\
&=&t^TS_n(\theta)-\halfe
t^TI(\theta)t+h^TS_n\Big(\theta+{t\over{\sqrt
n}}\Big)-\halfe h^TI(\theta+{t\over{\sqrt n}}\Big)h\label{3.1.d}\\
&&\quad -(t+h)^TS_n(\theta)+\halfe(t+h)^TI(\theta)(t+h)+o_P(1)\nonumber\\
&=&h^T\Big\{S_n(\theta+{t\over{\sqrt
n}}\Big)-S_n(\theta)+I(\theta)t\Big\}+o_P(1).\nonumber
\end{eqnarray}
Since (\ref{3.1.d}) holds for any $h\in{\mathbb R}^k$, this yields
(\ref{3.9}). This enables us to complete the proof of (\ref{3.8})
to the case $M>0$. Consider $\theta_n\to\theta, \ t_n\to t$. It
suffices to prove
\begin{equation}\label{3.1.e}
S_n(\theta_n){\buildrel\D\over\rightarrow_{\theta_n+{t_n\over\sqn}}}\Nn(I(
\theta)t, I(\theta))
\end{equation}
in view of the continuity of $\theta\to I(\theta)$.

Now, (\ref{3.8}) with $M=0$ implies
\begin{equation}\label{3.1.f}
S_n(\theta_n+{t_n\over\sqn}){\buildrel\D\over\rightarrow_{\theta_n+{t_n
\over\sqn}}} \Nn(0, I(\theta)),
\end{equation}
whereas (\ref{3.9}) yields
\begin{equation}\label{3.1.g}
S_n\Big(\theta_n+{t_n\over\sqn}\Big)=S_n(\theta_n)-I(\theta_n)t_n+o_P(1)
\end{equation}
both under $\theta_n$ and, by contiguity, under
$\theta_n+t_n/\sqn$. Combining (\ref{3.1.f}) and (\ref{3.1.g}) we
obtain (\ref{3.1.e}) and complete the proof of the LAN theorem.
\hfill $\Box$

\begin{example}[Normal Shift]\label{exam:3.1}
{\rm Let $I$ be nonsingular and let $X$ be $\Nn(It,I)$, $t
\in{\mathbb R}^k$. Note that in view of (\ref{3.8}) the
log-likelihood ratio
\begin{eqnarray}
\lefteqn{\ell(X;t)-\ell(X;0)}\nonumber\\
&=&-\halfe(X-It)^TI^{-1}(X-It)+
\halfe X^TI^{-1}X\label{3.10}\\
&=&t^TX-\halfe t^TIt\nonumber
\end{eqnarray}
has the same structure as the right-hand side of (\ref{3.6}).
This explains the name LAN.

Clearly, the best (shift equivariant) estimator of $t$ based on
$X$ is $I^{-1}X$. Comparing (\ref{3.6}) and (\ref{3.10}) we see
that this indicates that good estimators of $t$ should be equal to
$I^{-1}(\theta)S_n(\theta)$ approximately and hence good
estimators of $\theta$ to
\begin{equation}\label{3.11}
\theta+{1\over{\sqrt n}}I^{-1}(\theta)S_n(\theta)=\theta+\ene\nn
I^{-1}(\theta)\del(X_i;\theta).
\end{equation}
This is exactly what we will prove in Chapter \ref{chap:6}. \hfill
$\Box$ }
\end{example}

\section{Exercises Chapter {\ref{chap:3}}}

\begin{exercise}[Location Model]\label{exer:3.1}
{\rm Let $g$ be a density on ${\mathbb R}$, which is absolutely
continuous with respect to Lebesgue measure with Radon-Nikodym
derivative $g'$, such that the Fisher information for location is
finite, that is
\begin{equation}\label{3.12}
I(g)=\int(g'/g)^2g<\infty.
\end{equation}

An example of such a density is the Laplace or double exponential
density $g(x)=\halfe\exp(-|x|)$. Note that this density does {\sl
not} yield a location model satisfying the conditions of
Proposition \ref{prop:3.1}. Nevertheless, it can be shown along
the lines of the proof of Proposition \ref{prop:3.1}, that the
location model is regular for every $g$ satisfying (\ref{3.12}).
See Example 2.1.2, p.15 of BKRW. \hfill $\Box$ }
\end{exercise}

\begin{exercise}[Linear Regression]\label{exer:3.2}
{\rm Compute score function $\del$ and Fisher information matrix
for the linear regression model of Example \ref{exam:1.2} with $G$
fixed, and formulate conditions needed for  regularity.

\hfill $\Box$ }
\end{exercise}

\begin{exercise}[Cox's Proportional Hazards Model]\label{exer:3.3}
{\rm Do the same as in Exercise \ref{exer:3.2}, for the Cox model
of Example \ref{exam:1.3} with $\lambda$ fixed.

{\em Hint} Introduce the cumulative hazard function
$\Lambda\,,\,\Lambda(y)=\int_0^y\lambda\,,$ to simplify
computations and to make them more transparent. Verify that the
conditional distribution of $\exp(\nu^TZ)\Lambda(Y)$ given $Z=z$
is exponential with shape parameter 1. \hfill $\Box$ }
\end{exercise}

\begin{exercise}[Exponential Shift Model]\label{exer:3.4}
{\rm Let $X_1, X_2, \dots,X_n$ be i.i.d. exponential random
variables with scale parameter 1 and location parameter $\theta$,
i.e. with density $p(\theta)$ with respect to Lebesgue measure on
${\mathbb R}$ given by
$$p(x;\theta)= e^{-(x-\theta)}{\bf 1}_{[x>\theta]}.$$
Study the asymptotic behavior of the log-likelihood ratio
$$L_n(\theta+{t\over n})-L_n(\theta)$$ under $\theta$ for $t\in{\mathbb R}$ with
$L_n(\theta)$ defined as in (\ref{3.5}) and note that the limit
experiment may be described via one exponential random variable
with scale parameter 1 and location parameter $t$.
 \hfill $\Box$ }
\end{exercise}

\chapter{Local Asymptotic Spread Theorem}\label{chap:4}

The Local Asymptotic Normality derived in the previous chapter for
regular parametric models suggests a local asymptotic approach as
follows. Fix $\theta_0\in \Theta$ and consider parameter values at
distance $O(1/\sqrt{n})$ from $\theta_0$. In the framework of the
spread inequality this can be implemented by choosing
\begin{equation}\label{new4.1}
\vart = \theta_0+{\sigma\over \sqrt{n}}\, \vars,
\end{equation}
where the random vector $\vars$ has density $w_0$ on $({\mathbb
R}^k,{\cal B}^k)$ and $\sigma >0$ will be taken to tend to
infinity after the limit for $n\to\infty$ has been taken. As we
will prove, the General Spread Inequality from Theorem
\ref{thm:2.3} implies

\begin{theorem}[Local Asymptotic Spread Theorem]\label{thm:new4.1}
Let $\Pp=\{P_\theta :\theta\in\Theta\}$ be a regular parametric
model and let $\{T_n\}_{n\in{\mathbb N}}$ be a sequence of
estimators of $\nu(P_\theta)=q(\theta)$ with $q:{\mathbb
R}^k\to{\mathbb R}^m$ twice continuously differentiable. Fix
$\theta_0\in\Theta$ and let $\vart$ be defined as in
{\rm(\ref{new4.1})}, where $w_0$ is a differentiable density with
bounded support and with derivative vector $\dot w_0$ satisfying
\begin{equation}\label{new4.2}
\int |\dot w_0|^2/w_0 < \infty.
\end{equation}
Assume that $\dot q(\theta_0)I^{-1}(\theta_0)\dot q^T(\theta_0)$
is nonsingular. For $a\in {\mathbb R}^m$, denote the distribution
function of $\sqrt{n}\, a^T (T_n - q(\vart))$ by $G_{n,\sigma,a}$.
Asymptotically, as $n\to \infty$ and subsequently
$\sigma\to\infty$, $G_{n,\sigma,a}$ is at least as spread out as
${\cal N}(0,a^T\dot q(\theta_0)I^{-1}(\theta_0)\dot
q^T(\theta_0)a)$, more precisely
\begin{eqnarray}
\lefteqn{\limsup_{\sigma\to\infty}\, \limsup_{n\to\infty}\
G^{-1}_{n,\sigma,a}(v)-G^{-1}_{n,\sigma,a}(u)}
\nonumber\\
&\geq&
\liminf_{\sigma\to\infty\, }\liminf_{n\to\infty}\
G^{-1}_{n,\sigma,a}(v)-G^{-1}_{n,\sigma,a}(u)\label{new4.3}\\
&\geq& (a^T\dot q(\theta_0)I^{-1}(\theta_0)\dot
q^T(\theta_0)a)^{1/2}(\Phi^{-1}(v) - \Phi^{-1}(u)),\quad 0<u<v<1.
\nonumber
\end{eqnarray}
Moreover, equalities hold in {\rm(\ref{new4.3})} for all $a\in
{\mathbb R}^m$ iff under $\theta_0$
\begin{equation}\label{new4.4}
\sqrt{n}\, (T_n-q(\theta_0) - {1\over n}\sum_{i=1}^n \dot
q(\theta_0)I^{-1}(\theta_0)\dot\ell(X_i;\theta_0))
\stackrel{P}{\to} 0.
\end{equation}
\end{theorem}

\noindent{\bf Proof}. As shown in note \ref{note:6.5}
in Appendix \ref{noteschap4} it follows
from the regularity of $\Pp$ and the properties of $w_0$ and $q$
that the General Spread Inequality, Theorem \ref{thm:2.3}, may be
applied here, if given $\sigma>0$, the sample size $n$ is large
enough. This yields as a lower bound to $G_{n,\sigma,a}$ the
distribution function $K_{n,\sigma,a}$, the score statistic of
which is (take $a\not= 0$)
\begin{eqnarray}
\lefteqn{S_{n,\sigma,a} = {{a^T\dot
q(\theta_0)I^{-1}(\theta_0)}\over{a^T\dot q(\theta_0)
I^{-1}(\theta_0)\dot q^T(\vart) a}}}\label{new4.5}\\
&\times&\Big\{ {1\over \sqrt{n}}\sum_{i=1}^n \dot\ell(X_i;\vart)
+{1\over\sigma}\, {\dot w_0\over w_0} (\vars) - {1\over
\sqrt{n}}\, {{\sum_{h=1}^m a_h\ddot q_h(\vart)I^{-1}(\theta_0)\dot
q^T(\theta_0) a} \over{a^T\dot q(\theta_0) I^{-1}(\theta_0)\dot
q^T(\vart) a}} \Big\}\nonumber
\end{eqnarray}
and has distribution function $H_{n,\sigma,a}$. We have chosen
here $b=I^{-1}(\theta_0)\dot q^T(\theta_0) a$ and we have replaced
$a$ from Theorem \ref{thm:2.3} by $\sqrt{n}a$. Since the support
of $w_0$ is bounded, $\dot q$ and $I$ are continuous, and $\dot
q(\theta_0)I^{-1}(\theta_0)\dot q^T(\theta_0)$ is nonsingular, any
choice of $b\not= 0$ is allowed in the sense of (\ref{new2.21}),
provided $n$ is large enough.

By the LAN Theorem \ref{thm:3.1}, in particular (\ref{3.8}) with
$t=0$, by the finiteness of $\ex |\dot w_0(\vars)/w_0(\vars)|^2$,
by the boundedness of the support of $w_0$, and by the continuity
of $\dot q$ and $\ddot q$ we see that $S_{n,\sigma,a}$ converges
in distribution to ${\cal N}(0,(a^T\dot
q(\theta_0)I^{-1}(\theta_0)\dot q^T(\theta_0) a)^{-1})$ as
$n\to\infty$ and subsequently $\sigma\to\infty$. Since (recall
$a\not= 0$ and the nonsingularity of $ \dot
q(\theta_0)I^{-1}(\theta_0)\dot q^T(\theta_0)$)
\begin{equation}\label{new4.6}
\lim_{\sigma\to\infty}\lim_{n\to\infty} \ex S^2_{n,\sigma,a} =
(a^T\dot q(\theta_0)I^{-1}(\theta_0)\dot q^T(\theta_0)
a)^{-1}<\infty
\end{equation}
holds, $|S_{n,\sigma,a}|$ is uniformly integrable and this
suffices to show that $\int_s^1 H^{-1}_{n,\sigma,a}(t)dt,
K^{-1}_{n,\sigma,a}$, and $K_{n,\sigma,a}$ converge (cf. Lemma 2.1
and its proof of Klaassen (1989a)). A straightforward computation
shows that the limit of $K_{n,\sigma,a}$ is ${\cal N}(0,a^T\dot
q(\theta_0)I^{-1}(\theta_0)\dot q^T(\theta_0) a)$, thus completing
the proof of (\ref{new4.3}).

By an asymptotic version of the argument leading to the Spread
Equality Theorem \ref{thm:2.2}, we see that the equalities hold in
(\ref{new4.3}) for all $a\in {\mathbb R}^m$ iff
\begin{equation}\label{new4.7}
\sqrt{n}\,a^T(T_n-q(\vart))-G_{n,\sigma,a}^{-1}(H_{n,\sigma,a}(S_{n,\sigma,a}))
\stackrel{P}{\to}0, \ \mbox{as}\ n\to\infty,\ \sigma\to\infty.
\end{equation}
However, in view of the smoothness (\ref{3.9}) and the contiguity
in Corollary \ref{cor:3.2}, both following from the LAN Theorem
\ref{thm:3.1}, and in view of the asymptotic normality of both
$G_{n,\sigma,a}$ and $H_{n,\sigma,a}$ the convergence in
(\ref{new4.7}) is equivalent to the one in (\ref{new4.4}).
\hfill$\Box$

\begin{remark}
{\rm Let the random $m$-vector $Z$ be ${\cal N}(0,\dot
q(\theta_0)I^{-1}(\theta_0)\dot q^T(\theta_0))$ distributed. The
random variable $a^TZ$ has a normal distribution with variance
$a^T\dot q(\theta_0)I^{-1}(\theta_0)\dot q^T(\theta_0)a$.
Therefore, we might say that asymptotically as $n\to\infty$ and
subsequently $\sigma\to\infty$, the random vector
\begin{equation}\label{new4.8}
\sqrt{n}\,(T_n-q(\vart))
\end{equation}
is at least as spread out as a normal distribution with covariance
matrix $\dot q(\theta_0)I^{-1}(\theta_0)\dot q^T(\theta_0)$; see
Exercise \ref{exer:2.9}. It is an open question if in the General
Spread Inequality there exist choices $b\in {\mathbb R}^k$ for
each $a\in{\mathbb R}^m$ such that the bounds $K_{a,b}$ can be
viewed as stemming from an $m$-dimensional random vector $Z$ via
the linear combinations $a^TZ$ for all $a\in {\mathbb R}^m$. }
\end{remark}

\begin{remark}
{\rm If we would not have chosen $b=I^{-1}(\theta_0)\dot
q^T(\theta_0)a$ in (\ref{new4.5}), then we would have gotten
\begin{equation}\label{new4.9}
\Big({{(b^T\dot
q^T(\theta_0)a)^2}\over{b^TI(\theta_0)b}}\Big)^{1/2}(\Phi^{-1}(v)-\Phi^{-1}(u))
\end{equation}
as a lower bound in (\ref{new4.3}). If we differentiate,
supervised by Lagrange, the expression $b^TI(\theta_0)b+\lambda
(b^T\dot q^T(\theta_0)a-1)$ with respect to $b$ we obtain the
right hand side of (\ref{new4.3}) as the maximum of
(\ref{new4.9}).

}
\end{remark}

The Local Asymptotic Spread Theorem renders meaning to the
following definitions.

\begin{definition}[Linearity]\label{newdefn:4.1}
{\rm The sequence of estimators $\{T_n\}$ is {\em  asymptotically
linear} at $P_0$ in the influence function $\psi:\,
\Xx\times\Pp\to{\mathbb R}^m$ with $\expar_P|\psi(X,P)|^2<\infty,\
\expar_P\psi(X,P)=0$, if
\begin{equation}\label{new4.10}
\limy\, P_0\Big({\sqrt n}\Big|
T_n-\Big\{\nu(P_0)+\ene\nn\psi(X_i,P_0)\Big\}\Big|>\varepsilon\Big)=0
\end{equation}
for every positive $\varepsilon$.
}
\end{definition}

\begin{definition}[ Efficient Influence Function,
 Information  Bound]\label{newdefn:4.2}
{\rm \quad\\
The {\em efficient influence function} at $P_\theta$ in estimating
$\nu$ within the model $\Pp$ is defined by
\begin{equation}\label{new4.11}
\tell(x;P_\theta\mid \nu,\Pp)={\dot
q}(\theta)I^{-1}(\theta)\del(x;\theta), \quad x\in\Xx.
\end{equation}
The {\em information bound} at $P_\theta$ for estimating $\nu$
within the model $\Pp$ is defined as
\begin{equation}\label{new4.12}
I^{-1}(P_\theta\mid\nu,\Pp)={\dot q}(\theta)I^{-1}(\theta){\dot
q}^T(\theta).
\end{equation}
Note that
\begin{equation}\label{new4.13}
\expar_\theta\tell\tell^T(X;P_\theta\mid\nu,\Pp)=I^{-1}(P_\theta\mid\nu,\Pp).
\end{equation}
}
\end{definition}

\begin{definition}[Efficiency]\label{newdefn:4.3}
{\rm The sequence of estimators $\{T_n\}_{n\in {\mathbb N}}$ is
(locally asymptotically) {\em efficient} at $P_0$ in estimating
$q(\theta)$ if it is asymptotically linear at $P_0$ in the
efficient influence function, i.e. if it satisfies (\ref{new4.4}),
or equivalently, if for each $\{\theta_n\}$ with
$\sqrt{n}(\theta_n-\theta_0)$ bounded
\begin{equation}\label{new4.14}
\sqrt{n}(T_n-q(\theta_n)-{1\over n}\sum_{i=1}^n \dot
q(\theta_n)I^{-1}(\theta_n)\dot\ell(X_i;\theta_n))\stackrel{P}{\to}
0
\end{equation}
under $\theta_n$ (cf. the last sentence of the proof of Theorem
\ref{thm:new4.1}). }
\end{definition}

\section{Exercises Chapter {\ref{chap:4}}}

\begin{exercise}[Efficient Influence Functions]\label{exer:4.1}
{\rm Compute the efficient influence functions for estimation of
the Euclidean parameter $\nu$ in the models of Exercises
\ref{exer:3.1}, \ref{exer:3.2}, and \ref{exer:3.3} with $G$ fixed.
\hfill $\Box$ }
\end{exercise}

\begin{exercise}[Weak Convergence of Quantile Functions]\label{exer:4.2}
{\rm Let $F_n, n=1,2,\dots,$ and $F$ be distribution functions
with $F_n$ weakly converging to $F$. Prove that for all $u$ at
which $F^{-1}(\cdot)$ is continuous, $F_n^{-1}(u)$ converges to
$F^{-1}(u)$. We will say that $F_n^{-1}$ converges weakly to
$F^{-1}.$ Also prove that the weak convergence of $F_n^{-1}$ to
$F^{-1}$ implies that $F_n$ converges weakly to $F.$

\hfill $\Box$ }
\end{exercise}

\begin{exercise}[Weak Convergence of Spread Bounds]\label{exer:4.3}
{\rm Let $K_n^{-1}, n=1,2,\dots,$ and $K^{-1}$ be quantile
functions as defined by (\ref{2.11}) with $H$ replaced by $H_n$
for $K_n^{-1}, n=1,2,\dots.$ Here $H$ and $H_n$ are distribution
functions with mean 0. If $H_n$ converges weakly to $H$, then
$K_n$ converges weakly to $K.$ Prove this.

 \hfill $\Box$ }
\end{exercise}

\chapter{Geometric Interpretation}\label{chap:5}

The main issue in semiparametrics is estimation  in the presence of infinite
dimensional  nuisance parameters.
In the remainder of this section we will prepare for this by studying some
consequences of the Local Asymptotic Spread Theorem for estimation under finite dimensional
nuisance parameters.

Consider the regular parametric model  $\Pp$  from Definition \ref{3.1} with parameter
 $\theta\in{\mathbb R}^k$.
This parameter  $\theta={\nu\choose\eta}$  is split up into the
parameter of interest  $\nu\in{\mathbb R}^m$  and the nuisance
parameter  $\eta\in{\mathbb R}^{k-m}$. Accordingly, we write the
score function  $\del$  as
 $\del(\theta)=\del={\del_1 \choose \del_2}$  and the Fisher information matrix as
\begin{equation*}
I(\theta) = I =
\begin{pmatrix}
   I_{11} & I_{12} \\
   I_{21} & I_{22}
\end{pmatrix}
\ {\rm with} \ I_{ij}=E_\theta\del_i\del_j(X;\theta).
\end{equation*}
In the situation with
$\eta$ known we have a submodel of  $\Pp$  which is again regular
parametric with score function  $\del_1$,  and Fisher information
matrix  $I_{11}$. In view of the Local Asymptotic Spread theorem
an estimator sequence $\{T_n\}$  of  $\nu$ is uniformly
asymptotically efficient within this submodel if it is uniformly
asymptotically linear in the efficient influence function
$I_{11}^{-1}\del_1$, that is, if for all  $\nu$  with
$\theta=(\nu^T,\eta^T)^T\in\Theta$  and every sequence  $\nu_n$
with $\nu_n\to\nu$  and
 $\theta_n=\Big(\nu^T_n,\eta^T\Big)^T\in\Th$,
$${\sqrt n}\Big(T_n-\nu_n-\ene\nn I_{11}^{-1}(\theta_n)\del_1(X_i;\theta_n)\Big)
{\buildrel P\over\rightarrow_{\theta_n}}\,0.$$ Furthermore, the
information bound is  $I_{11}^{-1}$.

Within the full model  $\Pp$  with  $\eta$  unknown the Local
Asymptotic Spread theorem shows that  $\{T_n\}$  is uniformly
asymptotically efficient in estimating $\nu\,,$ if it is uniformly
asymptotically linear in the efficient influence function
\begin{equation}\label{4.12}
\tell_1(\theta)=(J \ 0)I^{-1}(\theta)\del(\theta),
\end{equation}
where  $J$  is the  $m\times m$  identity matrix and $0$ the
$m\times(k-m)$ null matrix. By computing $I^{-1}I$ it  is
straightforward to check that we may write
\begin{equation}\label{4.13}
I^{-1}(\theta)=I^{-1}=
\begin{pmatrix}
   I^{-1}_{11.2}                & -I^{-1}_{11.2}I_{12}I^{-1}_{22} \\
-I^{-1}_{22.1}I_{21}I^{-1}_{11} & I^{-1}_{22.1}
\end{pmatrix}
\end{equation}
with
\begin{equation}\label{4.14}
I_{11.2}=I_{11}-I_{12}I^{-1}_{22}I_{21}, \
I_{22.1}=I_{22}-I_{21}I^{-1}_{11}I_{12}.
\end{equation}
Furthermore, we will write
\begin{equation}\label{4.15}
\stel_1=\del_1-I_{12}I^{-1}_{22}\del_2.
\end{equation}
In this notation (\ref{4.12}) becomes
\begin{eqnarray*}
\tell_1(\theta)&=&\tell_1=I^{-1}_{11.2}(J\,,\,-I_{12}I^{-1}_{22})\del
=I^{-1}_{11.2}(\del_1-I_{12}I^{-1}_{22}\del_2)=I^{-1}_{11.2}\stel_1\\
&=&(E\stel_1{\stel}_1^T)^{-1}\stel_1,
\end{eqnarray*}
since  $\ex\stel_1{\stel}^T_1=\ex(\del_1-I_{12}I^{-1}_{22}\del_2)
(\del_1-I_{12}I^{-1}_{22}\del_2)^T=I_{11}-I_{12}I^{-1}_{22}I_{21}=I_{11.2}$
holds.
Summarizing, we have:
In the restricted model the efficient influence function equals
\begin{equation}\label{4.16}
(\ex\del_1\del_1^T)^{-1}\del_1 \ {\rm with} \
(\ex\del_1\del^T_1)^{-1}=I_{11}^{-1}
\end{equation}
as information bound.
In the full model the  efficient influence function equals
\begin{equation}\label{4.17}
(\ex\stel_1{\stel}^T_1)^{-1}\stel_1 \ {\rm with} \
(\ex\stel_1{\stel}^T_1)^{-1}= (I^*)^{-1}
\end{equation}
as information bound.
In view of the similarities between (\ref{4.16}) and (\ref{4.17}) we will call
 $\stel_1$  the {\sl efficient score function}  for estimation of  $\nu$  within
 $\Pp$ and $I^*$ the {\sl efficient Fisher information matrix}.

In the space of mean  $0$  random variables of functions of  $X$
with  $X\sim P_\theta$  we define an inner product by
\begin{equation}\label{4.18}
\langle f,g\rangle=\ex_\theta f(X)g(X).
\end{equation}
In the resulting Hilbert space each component of  $\stel_1$  may be viewed
as the projection on  $[\del_2]^{\perp}$  of the corresponding components of
 $\del_1$.
Here,  $[\del_2]$  is the (closed) linear span of all components of  $\del_2$
and  $[\del_2]^{\perp}$  is its orthocomplement.
Indeed,
$$
 \langle\stel_1,\, \del^T_2\rangle=
\langle\del_1 - I_{12} I^{-1}_{22}\del_2,\,\del^T_2\rangle
=\ex(\del_1-I_{12}I^{-1}_{22}\del_2)\del^T_2=0
$$
and hence  $\stel_1\in [\del_2]^{\perp}$.
Furthermore,
$$\del_1-\stel_1=I_{12}I^{-1}_{22}\del_2\perp[\del_2]^{\perp},$$
since each component of  $I_{12}I^{-1}_{22}\del_2$  belongs to
 $[\del_2]^{\perp\perp}=[\del_2]$.
In formula,
\begin{equation}\label{4.19}
\stel_1=\Pi\Big(\del_1\mid [\del_2]^{\perp}\Big)=\del_1-\Pi\Big(\del_1\mid
[\del_2]\Big).
\end{equation}
Furthermore, we have
\begin{equation}\label{4.20}
I^{-1}_{11}\del_1=\Pi\Big(\tell_1\mid [\del_1]\Big),
\end{equation}
as can be seen as follows
\begin{eqnarray}
\lefteqn{\langle\tell_1-I^{-1}_{11}\del_1,\,\del^T_1\rangle}\nonumber\\
&=&\ex I^{-1}_{11.2}\Big(\del_1-I_{12}I^{-1}_{22}\del_2\Big)
\del_1^T-I^{-1}_{11}\ex\del_1\del^T_1\label{4.21}\\
&=&I^{-1}_{11.2}I_{11.2}-I^{-1}_{11}I_{11}=0.\nonumber
\end{eqnarray}

\setlength{\unitlength}{0.8cm}
\begin{figure}[h]
\begin{picture}(10,6)(-8,0)
\thicklines \put(0,0){\line(1,1){6}}
\multiput(-2,-2)(0.8,0.8){4}{\line(1,1){0.5}}
\put(0,0){\vector(1,1){5}} \put(0,0){\vector(0,1){5}}
\put(0,0){\vector(3,1){6}} \put(5,5){\vector(0,-1){3.32}}
\put(6,2){\vector(-1,1){2}} \put(0,0){\vector(3,1){5}}
\put(-3,-1){\line(1,2){2}} \put(7,-1){\line(1,2){2}}
\put(-1,3){\line(1,0){10}} \put(-3,-1){\line(1,0){10}}
\put(-0.7,4){$\dot\ell_2$} \put(4.3,5){$\dot\ell_1$}
\put(2.6,4){$I^{-1}_{11}\dot\ell_1$} \put(6,1.3){$\tilde\ell_1$}
\put(5,1){$\ell_1^*$} \put(6.2,-0.5){$[\dot\ell_2]^\perp$}
\put(6.2,6.2){$[\dot\ell_1]$}
\end{picture}
\vspace{1.5cm}
\caption{\label{fig:5.1}Efficient influence functions and efficient score functions}
\end{figure}
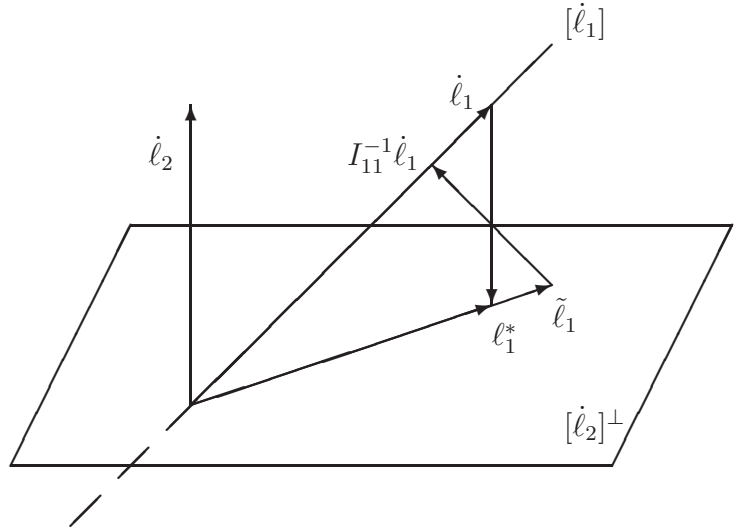

\noindent Formula (\ref{4.19}) shows that, knowing the score function within the smaller
model, we can obtain the score function of the same parameter within the bigger
model by projection.
Similarly, (\ref{4.20}) shows that the efficient influence function for the
smaller model may be obtained by projection of the efficient influence
function of
the bigger model; see Figure \ref{fig:5.1}.
Both devices are useful in semiparametrics as we will see in the next
Chapter \ref{chap:6}

\begin{example}[Normal distribution]\label{exam:4.1}
{\rm \quad\\
Consider the normal family discussed in Example \ref{exam:normal}.
The parameter $\theta$ equals $(\mu,\sigma)$. For the inverse of
the Fisher information matrix for one observation we get
\begin{equation}
I(\theta)^{-1} = \left(\begin{array}{cc}
\sigma^2&0\\
0&\sigma^2\over2
\end{array}
\right).
\end{equation}
The efficient influence function for estimating $\mu$ with
$\sigma$ known, as given by (\ref{4.16}), equals (cf.
(\ref{3.13}))
\begin{equation}
\sigma^2{(x-\mu)\over \sigma^2} = x-\mu.
\end{equation}
The efficient influence function for estimating $\mu$ with
$\sigma$ unknown (\ref{4.17}) is also equal to $ x-\mu$ since the
off diagonal elements of the Fisher information matrix vanish. The
information lower bound for estimating $\mu$ is equal to
$\sigma^2$ in both cases. Here $\mu$ can be estimated just as well
whether we know $\sigma$ or not. The bounds are attained by the
sample mean. \hfill $\Box$ }
\end{example}

\begin{example}[Cox's Proportional Hazards Model (parametric)]
\label{exam:4.2}
{\rm \quad\\
The semiparametric version of the model is described in Example
\ref{exam:1.3}. One obser\-ves i.i.d. copies of $X=(Y,Z)$ with $Z$ an
$m$-vector of covariates.
Given $Z=z$ the survival time $Y$ has hazard function
\begin{equation}\label{ex:3.1}
\la(y\mid z,\nu)=\exp(z^T\nu)\la(y), \quad y>0.
\end{equation}
Here $\la(y)=g(y)/(1-G(y))$ is the unknown baseline hazard function
with $G$ a distribution function with density $g$.
We get a parametric version of this model by assuming that $g$ is a
density from a parametric family. Suppose for instance that $g$ is equal
to an Exp($\lambda$) density, for some constant $\lambda >0$. The corresponding
hazard function is constant $\la$. We will be interested in estimation of
$\nu$ in the presence of the nuisance parameter $\la$.
Given $Z=z$ the hazard function of $Y$ is
\begin{equation}\label{ex:3.2}
\la(y\mid z,\nu)=\exp(z^T\nu)\la, \quad y>0.
\end{equation}
Since this hazard is independent of $y$ the conditional distribution
of the survival time $Y$ given $Z=z$ is Exp($\exp(z^T\nu)\la$).
Hence the unconditional density of the pair $(Y,Z)$ is
\begin{equation}\label{ex:3.3}
p_\theta(y,z)=\exp(z^T\nu)\la\exp\Big(-\exp(z^T\nu)\la
y\Big)p_Z(z),
\end{equation}
where $p_Z$ is the known density of the covariate vector $Z$ and
$\theta= (\nu^T,\lambda)^T$.

For convenience we assume that $Z$ and $\nu$ are one dimensional. Let us
consider the case of one observation first.
We get
$$
\ell(y,z)=z\nu + \log(\la) - \exp(z\nu)\la y + \log(p_Z(z))
$$
and
\begin{eqnarray*}
\del_1(y,z)&=& {\partial\over\partial \nu}\ell(y,z) = z(1-\exp(z\nu)\la y),\\
\del_2(y,z)&=& {\partial\over\partial \la}\ell(y,z) ={1\over\la}-\exp(z\nu)y.
\end{eqnarray*}
Note that indeed the expectations of $\del_1(Y,Z)$ and $\del_2(Y,Z)$ vanish.

Computing the Fisher information matrix for one observation we get
\begin{equation}
I(\theta) = I((\nu,\la)) = \left(\begin{array}{cc}
\ex Z^2&{1\over \la}\ex Z\\
{1\over \la}\ex Z&{1\over \la^2}
\end{array}
\right).
\end{equation}
Its inverse is given by
\begin{equation}
I(\theta)^{-1}  = {1\over \var Z} \left(\begin{array}{cc}
1&-\la\ex Z\\
-\la\ex Z&\la^2 \ex Z^2
\end{array}
\right).
\end{equation}
The efficient influence function for estimating $\nu$ in the
restricted model where $\lambda$ is known, is given by
(\ref{4.16}) and equals
\begin{equation}
{1\over \ex Z^2} z(1-\exp(z\nu)\la y).
\end{equation}
In the larger model with both $\nu$ and $\la$ unknown the  score
function $\ell_1^*$ for estimating $\nu$ equals (cf. (\ref{4.17}))
\begin{eqnarray*}
\lefteqn{z(1-\exp(z\nu)\la y) -\la \ex Z({1\over\la}-\exp(z\nu)y)}\\
&=& (z- \ex Z) (1-\exp(z\nu)\la y).
\end{eqnarray*}
The information lower bounds for the restricted and unrestricted model
 are equal to
\begin{equation}
{1\over \ex Z^2}\quad\mbox{and}\quad {1\over \var Z}\ .
\end{equation}
Since $\var Z \leq \ex Z^2$ the lower bound in the unrestricted
model is indeed larger than the lower bound for the restricted
model. Equality occurs if and only if $\ex Z=0$. It can be shown
that maximum likelihood estimators attain the lower bounds.
 \hfill $\Box$
}
\end{example}

\section{Exercises Chapter {\ref{chap:5}}}

\begin{exercise}[Bivariate Normality]\label{exer:4.1}
{\rm Consider the regular parametric mo\-del  $\Pp$  with
$\theta=(\nu,\Sigma)$ and
 $P_\theta$  the normal distribution with mean vector  $\nu={\nu_1\choose
\nu_2}$  and covariance matrix  $\Sigma=({\eta^2_1\atop\rho\eta_1\eta_2}
{\rho\eta_1\eta_2\atop\eta^2_2})$.
Show that knowledge of  $\Sigma$  does not help asymptotically in estimating
 $\nu$,   as compared to the situation with  $\Sigma$  unknown.
Prove that the information lower bound in estimating the
correlation coefficient $\rho$ is
 $(1-\rho^2)^2$;  see Example 2.4.6 of BKRW (1993).\hfill $\Box$
}
\end{exercise}

\begin{exercise}[Plug-in Estimators]\label{exer:4.2}
{\rm Consider the regular parametric mo\-del  $\Pp$  from
Definition \ref{defn:3.1} with the parameter
$\theta={\nu\choose\eta}$  split up into the parameter of interest
 $\nu\in{\mathbb R}^m$  and the nuisance parameter  $\eta\in{\mathbb R}^{k-m}$.
Let  ${\hat\eta}_n$  be an efficient estimator of  $\eta$  based
on  $n$  i.i.d. observations  $\Xv$  from  $P_\theta$. For  $\eta$
known let  ${\hat\nu}_n(\eta)$  be an efficient estimator of
 $\nu$
within the restricted model  $\Pp_1(\eta)=\{P_\theta:\,\nu\in
N\},\eta$  fixed. Use (\ref{3.9}) to show heuristically that
${\hat\nu}_n({\hat\eta}_n)$  is an efficient estimator of  $\nu$
within the full model  $\Pp$.

A sample splitting technique as will be discussed in Chapter \ref{chap:8} enables one to
construct a version of  ${\hat\nu}_n({\hat\eta}_n)$  that can be
proved rigorously to be efficient. \hfill $\Box$ }
\end{exercise}

\chapter{Semiparametric Models}\label{chap:6}

In Chapter \ref{chap:4} we have discussed the Local Asymptotic Spread Theorem as a bound
 on
the asymptotic performance of estimators of a parameter
$\nu(P_\theta)=q(\theta)$ in regular parametric models. The
regularity condition we needed on  $\nu:\,\Pp\to{\mathbb R}^m$,
was that
 $q:\,\Th\to{\mathbb R}^m$  be twice (continuously) differentiable.
For the more general nonparametric and semiparametric models $\Pp$
that we will study in this chapter, we will need a generalization
of the differentiability of  $q(\cdot)$,  to wit pathwise
differentiability of $\nu(\cdot)$. To this end we need the concept
of tangent space.

\begin{definition}[Paths and Tangent Space]\label{defn:5.0}
{\rm Let $\Pp$ be a model and fix $P_0\in\Pp$. A one-dimensional
regular parametric submodel $\{P_\eta:|\eta|<\varepsilon\}, \
\varepsilon>0$, containing $P_0$ is called a {\em curve} or {\em path}
through $P_0$. Its score function $h$ at $P_0$ is called a tangent
and the collection of all such possible tangents $h$ is called the
{\sl tangent set}  $\dP^0$  of  $\Pp$  at  $P_0$. The {\em tangent
space}  $\dot\Pp$  is defined as the closed linear span
 $[\dP^0]$  of  $\dP^0$  in  $\L_2(P_0)$.
}
\end{definition}

For $m=1$ our definition of pathwise differentiability reads as
follows.

\begin{definition}[Pathwise Differentiability]\label{defn:5.1}
{\rm The parameter  $\nu:\,\Pp\to{\mathbb R}$  is {\em pathwise
differentiable on  $\Pp$  at} $P_0\in\Pp$, if there exists a
bounded linear functional $\dun(P_0):\,{\dot\Pp}\to{\mathbb R}$
such that for every path $\{P_\eta:|\eta|<\varepsilon\}, \ \varepsilon>0$,
with tangent $h$ at $P_0$
\begin{equation}\label{5.1}
\nu(P_\eta)=\nu(P_0)+\eta\dun(P_0)(h)+o(\eta), \ {\rm as} \ \eta\to0.
\end{equation}
}
\end{definition}

We will denote the inner product and norm of the Hilbert space
$\L_2(P_0)$  by
 $\langle\cdot,\cdot\rangle_0$  and  $\|\cdot\|_0$  respectively.
Here  $\langle f,g\rangle_0=\expar_{P_0}f(X)g(X)$.
In fact, we are dealing with  $\L^0_2(P_0)=\{f\in\L_2(P_0):\expar_{P_0}f(X)=0\}$
usually.
Within this subspace  $\langle f,g\rangle_0$  may be interpreted as
 $\Cov_{P_0}(f(X),g(X))$  and  $\|f\|_0$  becomes the
standard deviation of  $f(X)$  under  $P_0$.

By the Riesz representation theorem, see for example Rudin (1966),
 there exists a unique
 $\dun\in{\dot\Pp}\subset\L^0_2(P_0)$  such that
\begin{equation}\label{5.2}
\dun(P_0)(h)=\langle\dun,h\rangle_0,\quad\mbox{for  all}\
h\in\dP.
\end{equation}
We will call this  $\dun\in{\dot\Pp}$  the {\em derivative} of
$\nu$.

In general,  $\nu=(\nu_1,\dots,\nu_m)^T:\,\Pp\to{\mathbb R}^m$
will be called {\em pathwise differentiable on  $\Pp$  at}  $P_0$
with derivative $\dun=(\dun_1,\dots,\dun_m)^T$  if each  $\nu_i:\,\Pp\to{\mathbb R}$  is
pathwise differentiable with derivative
$\dun_i\in{\dot\Pp}\subset\L^0_2(P_0)\,, i=1,\dots,m$.

\begin{example}[Regular Parametric Models]\label{exam:5.1}
{\rm Let  $\Pp=\{P_\theta:\,\theta\in\Theta\}, $
$\Theta\subset{\mathbb R}^k$,  be a regular parametric model with
score function  $\del$ at  $P_0=P_{\theta_0}, \ \theta_0$  fixed
in $\Theta$. The tangent space $\dP=[\del]$  at  $P_0$  is the
($k$-dimensional) closed linear span of the components of  $\del$.
Let $\nu(P_\theta)=q(\theta)$, where $q(\cdot)$  is differentiable
at $\theta_0$  with $(m{\rm x}k)$-matrix
$${\dot q}(\theta_0)=\left({{\partial q_i(\theta)} \over {\partial \theta_j}}\right)_{j=1,\dots,k}^{i=1,\dots,m}$$
of derivatives. Choose  $a\in{\mathbb R}^k, \ a\ne 0$. Then,
$P_\eta=P_{\theta_0+\eta a}$  defines a one-dimensional regular
parametric submodel of  $\Pp$  with score function
$a^T\del=\del^Ta$ at  $P_0$. Note
\begin{eqnarray}
\nu(P_\eta)&=&q(\theta_0+\eta a)=q(\theta_0)+\eta{\dot
q}(\theta_0)a+o(\eta)\label{5.3}\\
&=&\nu(P_0)+\eta\langle{\dot q}(\theta_0)I^{-1}(\theta_0)\del, \
\del^Ta\rangle_0+o(\eta).\nonumber
\end{eqnarray}
This shows that  $\nu:\, \Pp\to{\mathbb R}^m$  is pathwise
differentiable on  $\Pp$  at
 $P_0$  with derivative
\begin{equation}\label{5.4}
\dun={\dot q}(\theta_0)I^{-1}(\theta_0)\del=\tel(\cdot,P_0\mid\nu,
\Pp)
\end{equation}
the efficient influence function! \hfill $\Box$ }
\end{example}

In fact,  $\dun$  may be called the efficient influence function at  $P_0$  for
estimating  $\nu$  within  $\Pp$,  whether  $\Pp$  is regular parametric or
 not.
As an illustration we prove

\begin{theorem}[Regular Parametric Submodels]\label{thm:5.1}
 Let  $\nu:\,\Pp\to{\mathbb R}^m$  be pathwise differentiable at  $P_0$  with
derivative  $\dun$.
Let  $\Qq$  be any regular parametric submodel of  $\Pp$  with
 $I^{-1}(P_0\mid\nu,\Qq)$  well-defined.
If  $P_0\in\Qq$,  then
\begin{equation}\label{5.5}
\tel(\cdot,P_0\mid\nu,\Qq)=\Pi_0(\dun\mid\dQ),
\end{equation}
where  $\Pi_0$  denotes projection under  $\langle\cdot,\cdot\rangle_0$  within
 $\L^0_2(P_0)$.
Furthermore,
\begin{equation}\label{5.6}
I^{-1}(P_0\mid\nu,\Qq)\le\langle\dun,\dun^T\rangle_0
\end{equation}
with equality iff  $\dun\in\dQ^m$.
\end{theorem}

\noindent{\bf Proof.}

Let  $\Qq$  be {\em any} submodel of  $\Pp$  and let
$\nu_\Qq:\,\Qq\to{\mathbb R}^m$ denote the restriction of  $\nu$
to  $\Qq$. Then  $\nu_\Qq$  is pathwise differentiable with
derivative
 $\dun_\Qq\in\dQ\subset\dP$.
Note that (\ref{5.1}) implies that $\dun_\Qq(P_0):\,\dQ\to{\mathbb
R}^m$  is unique and equals the restriction of
$\dun(P_0):\,\dP\to{\mathbb R}^m$  to  $\dQ$. Consequently, we
have
\begin{eqnarray}
\lefteqn{\langle\dun_\Qq-\Pi_0(\dun\mid\dQ),h\rangle_0=\langle\dun_\Qq-
\dun,h\rangle_0}\label{5.a}\\
&=&\dun_\Qq(P_0)(h)-\dun(P_0)(h)=0, \ {\rm for\ all} \
h\in\dQ.\nonumber
\end{eqnarray}
Together with  $\dun_\Qq-\Pi_0(\dun\mid\dQ)\in\dQ$  this yields
\begin{equation}\label{5.b}
\dun_\Qq=\Pi_0(\dun\mid\dQ).
\end{equation}
If  $\Qq$  is a regular parametric submodel, then combination of (\ref{5.b})
and
(\ref{5.4}) of Example \ref{exam:5.1} proves (\ref{5.5}).

Furthermore, for any  $a\in{\mathbb R}^m$  this implies
\begin{eqnarray}
a^TI^{-1}(P_0\mid\nu,\Qq)a&=&a^T \expar_{P_0}\left(\tilde
l_Q\tilde l_Q^T\right)a= \expar_{P_0}\left( a^T \tilde l_Q \tilde
l_Q^T a \right)= \expar_{P_0}\left( a^T \tilde l_Q (\tilde l_Q
a)^T\right)
\nonumber\\
&=& \expar_{P_0} \left((a^T\tilde l_Q)^2\right) =
\|a^T\Pi_0(\dun\mid\dQ)\|_0^2=\|\Pi_0
(a^T\dun\mid\dQ)\|_0^2\label{5.c}\\
&\le&\|a^T\dun\|_0^2=a^T\langle\dun,\dun^T\rangle_0a. \nonumber
\end{eqnarray}
The resulting inequality with $a$ arbitrary is the meaning and
content of inequality (\ref{5.6}). \hfill$\Box$

\bigskip

Theorem \ref{thm:5.1} suggests the following definitions.

\begin{definition}[Efficient Influence Function and Information Bound]
\label{defn:5.2} {\rm If  $\nu:\,\Pp\to{\mathbb R}^m$  is pathwise
differentiable with derivative  $\dun$  at
 $P_0$,  then  $\dun$  is called the {\em efficient influence function} at
 $P_0$  for estimating  $\nu$  within  $\Pp$;  notation
\begin{equation}\label{5.7}
\tel=\tel(\cdot,P_0\mid\nu,\Pp)=\dun.
\end{equation}
The {\em information bound} at  $P_0$  for estimating  $\nu$  within  $\Pp$  is
defined as
\begin{equation}\label{5.8}
I^{-1}(P_0\mid\nu,\Pp)=\langle\dun,\dun^T\rangle_0.
\end{equation}
}
\end{definition}

A justification  for these definitions is given by the following
generalization of Theorem \ref{thm:new4.1}.

\begin{theorem}[LAS Theorem]\label{thm:new6.2}
 Let  $\nu:\,\Pp\to{\mathbb R}^m$  be pathwise differentiable on  $\Pp$  at  $P_0$
and let  $T_n$  be any estimator of  $\nu$. Let $G_{{\cal Q},
n,\sigma,a}$ for $a\in {\mathbb R}^m$ denote the distribution
function of $\sqrt{n}\, a^T (T_n - q(\vart))$ for any regular
parametric submodel ${\cal Q}=\{P_\theta\,:\, \theta \in\Theta\}$
with $q(\theta)=\nu(P_\theta)$ and with $\vartheta$ and $\sigma$
as in (\ref{new4.1}). If
\begin{equation}\label{new6.12}
[\dun]\subset\odP,
\end{equation}
holds with $\odP$ the closure of the tangent set  $\dP^0$,  then
asymptotically, as $n\to \infty$ and subsequently
$\sigma\to\infty$, $G_{{\cal Q},n,\sigma,a}$ is at least as spread
out as ${\cal N}(0,I^{-1}(P_0|\nu,\Pp))$, i.e.
\begin{eqnarray}
\lefteqn{\sup_{\cal Q}\limsup_{\sigma\to\infty}\,
\limsup_{n\to\infty}\ G^{-1}_{{\cal
Q},n,\sigma,a}(v)-G^{-1}_{{\cal Q},n,\sigma,a}(u)}
\nonumber\\
&\geq& \sup_{\cal Q}\liminf_{\sigma\to\infty\,
}\liminf_{n\to\infty}\
G^{-1}_{{\cal Q},n,\sigma,a}(v)-G^{-1}_{{\cal Q},n,\sigma,a}(u)\label{lasinequality}\\
&\geq& (a^TI^{-1}(P_0|\nu,\Pp)a)^{1/2}(\Phi^{-1}(v) -
\Phi^{-1}(u)),\quad 0<u<v<1. \nonumber
\end{eqnarray} Moreover,  equalities hold
iff $T_n$ is asymptotically linear with influence function
$\tel(\cdot,P_0\mid\nu,\Pp)$ under  $P_0$,  and  $T_n$ is called
efficient then at $P_0$.
\end{theorem}

The proof of this LAS Theorem is similar to the proof of
the LAS Theorem \ref{thm:new4.1} and involves an approximation procedure
based on (\ref{new6.12}).

\bigskip

In the remainder of this chapter we will study consequences of the
above results for semiparametric models  $\Pp=\{P_{\nu,G}:\,\nu\in
N, \ G\in\Gg\}$, where  $N\subset{\mathbb R}^m$  is open and
$\Gg$  is general. Fix  $\nu_0\in N$  and  $G_0\in\Gg$  and denote
$P_0=P_{\nu_0,G_0}$. We will assume that the parametric submodel
\begin{equation}\label{5.14}
\Pp_1=\Pp_1(G_0)=\{P_{\nu,G_0}:\,\nu\in N\}
\end{equation}
of  $\Pp$  is regular with score function $\del_1$. Furthermore,
we will denote the submodel with  $\nu=\nu_0$  fixed by
\begin{equation}\label{5.15}
\Pp_2=\Pp_2(\nu_0)=\{P_{\nu_o,G}:\, G\in\Gg\}
\end{equation}
and the tangent space of  $\Pp_2$  at  $P_0$  by  $\dP_2$.
Analogously to (\ref{4.19}) we introduce the {\em efficient score function} for
 $\nu$  by
\begin{equation}\label{5.16}
\stel_1(\cdot,P_0\mid\nu,\Pp)=\stel_1=\del_1-\Pi(\del_1\mid\dP_2).\
\end{equation}
This terminology is motivated by the following result.

\begin{theorem}[Efficient Score Function]\label{thm:5.3}
Let  $\Pp_1$  be regular and  $\stel_1$  be defined by {\rm(\ref{5.16})}.
If
$\Qq=\{P_{\nu,G_\ga}:\, \nu\in N, \ \ga\in\Gamma\}$  is a regular
parametric submodel of  $\Pp$  containing  $P_0$  and hence  $\Pp_1$,
then
\begin{equation}\label{5.17}
I(P_0\mid\nu,\Qq)\ge \ex\stel_1{\stel_1}^T
\end{equation}
with equality iff  $[\stel_1]\subset\dQ$. Furthermore, if
$\dP=\dP_1+\dP_2=[\del_1]+\dP_2$  and {\rm(\ref{5.16})} are fulfilled for
$\nu(P_{\nu,G})=\nu$, then
\begin{equation}\label{5.18}
\dun=\tel(\cdot,P_0\mid\nu,\Pp)=\tel_1=(\ex\stel_1{\stel_1}^T)^{-1}\stel_1
\end{equation}
and the information bound equals
\begin{equation}\label{5.19}
I^{-1}(P_0\mid\nu,\Pp)=(\ex\stel_1{\stel}_1^T)^{-1}.
\end{equation}
\end{theorem}

\noindent{\bf Proof.}

In view of (\ref{4.17}), (\ref{4.15}) and (\ref{4.14}),
the Pythagorean  Theorem yields for all
 $a\in{\mathbb R}^m$
\begin{eqnarray}
\lefteqn{a^TI(P_0\mid\nu,\Qq)a=\|a^T(\del_1-\Pi(\del_1\mid\dQ_2))\|^2_0}
\nonumber\\
&=&\|a^T\stel_1+\Pi(a^T\del_1\mid\dP_2)-\Pi(a^T\del_1\mid\dQ_2)\|^2_0
\label{5.3.a}\\
&=&\|a^T\stel_1\|_0^2+\|\Pi(a^T\del_1\mid\dP_2)-\Pi(a^T\del_1\mid\dQ_2)\|_0^2
\nonumber\\
&\ge& a^T(\ex\stel_1{\stel}^T_1)a,\nonumber
\end{eqnarray}
and hence (\ref{5.17}).

The pathwise differentiability of  $\nu$  with derivative  $\dun$  at  $P_0$
yields
\begin{equation}\label{5.3.b}
\nu-\nu_0=\nu(P_{\nu,G_0})-\nu(P_0)=\langle\dun,{\del}^T_1\rangle_0(\nu-\nu_0)
+o(|\nu-\nu_0|)
\end{equation}
and hence
\begin{equation}\label{5.3.c}
\langle\dun,{\del}^T_1\rangle_0 \ {\rm equals\ the} \ m\times m \ {\rm
identity\ matrix}.
\end{equation}
The differentiability also yields for any path  $\{P_{\nu_0,G_\ga}:\,
|\ga|<\varepsilon\}\subset\Pp_2$  with tangent  $h$,
\begin{equation}\label{5.3.d}
0=\nu(P_{\nu_o,G_\ga})-\nu(P_0)=\ga\langle\dun,h\rangle_0+o(|\ga|)
\end{equation}
and hence
\begin{equation}\label{5.3.e}
[\dun]\perp\dP_2.
\end{equation}

Since the efficient influence function  $\tel=\dun$  belongs to
 $\dP=[\del_1]+\dP_2$  there exists an  $m\times m$  matrix A and
 $h\in\dP_2^m$,  such that
\begin{equation}\label{5.3.f}
\dun=A\stel_1+h.
\end{equation}
By (\ref{5.3.e}) this implies
\begin{equation}\label{5.3.g}
0=\langle\dun,h^T\rangle_0=A\langle\stel_1,h^T\rangle_0+\langle
h,h^T\rangle_0=\langle h,h^T\rangle_0,
\end{equation}
that is  $h=0$  or in other words
\begin{equation}\label{5.3.h}
\dun=A\stel_1.
\end{equation}
But, by (\ref{5.3.c}) this shows that the identity matrix equals
\begin{equation}\label{5.3.i}
\langle\dun,{\del}^T_1\rangle_0=A\langle\stel_1,{\del}^T_1\rangle_0=
A\langle\stel_1,{\stel_1}^T\rangle_0
\end{equation}
and hence  $A=(\langle\stel_1,{\stel_1}^T\rangle_0)^{-1}$  or
\begin{equation}\label{5.3.j}
\dun=(\langle\stel_1,{\stel_1}^T\rangle_0)^{-1}\stel_1.
\end{equation}
\hfill $\Box$

\bigskip

\begin{remark}
Like we noted about (\ref{4.19}), the projection in (\ref{5.16}) and
Theorem \ref{thm:5.3} show that,
knowing the score function within the smaller, parametric submodel, we can
obtain the score function of the same parameter within the bigger,
semiparametric
model by projection.
Similarly, as we noted about (\ref{4.20}), Theorem \ref{thm:5.1}
 and the LAS Theorem \ref{thm:new6.2}
show that the efficient influence function for the smaller
(parametric) submodel may be obtained  by projection of the
efficient influence function of the bigger, nonparametric
model.
\end{remark}

\bigskip

\begin{example}[Construction of one dimensional families]\label{exam:5.2a}{\rm
Let $\cal G$ be the family of {\em all probability measures
dominated by $\mu$}. Let $G_0$ be a fixed element of
${\cal G}$ with density $g_0$ and let $h$ be an element
of $\L_2(G_0)$ such that $\int h(t)g_0(t)dt=0$. We construct a one
dimensional regular parametric subfamily with tangent $h$ as follows.

Suppose that the function $\Psi : {\mathbb R} \to (0,\infty)$
is bounded and continuously
differentiable with a bounded derivative $\Psi'$. Furthermore, let
$\Psi(0)=\Psi'(0)=1$ and let $\Psi'/\Psi$ be bounded. An example of such a
function is
$$\Psi(x) = 2(1+e^{-2x})^{-1}.$$
Define
\begin{equation}
g_\eta(x)={{g_0(x)\Psi(\eta h(x))}\over{\int_{-\infty}^\infty
g_0(t)\Psi(\eta h(t))dt}}, -\infty<\eta<\infty.
\end{equation}
We have
\begin{eqnarray}
{d\over{d \eta}} g_0(x)\Psi(\eta h(x))
&=&g_0(x)h(x)\Psi'(\eta h(x))\label{part:1}\\
{d\over{d \eta}}\int_{-\infty}^\infty
g_0(t)\Psi(\eta h(t))dt
&=&\int_{-\infty}^\infty g_0(t)h(t)\Psi'(\eta h(t))dt.
\label{part:2}
\end{eqnarray}
Note that for $\eta=0$ the derivatives (\ref{part:1}) and (\ref{part:2}) equal
$g_0(x)h(x)$ and zero respectively.
The derivative at zero of $g_\eta(x)$ equals
\begin{equation}
{d\over{d \eta}} g_\eta(x)\Big|_{\eta=0} = {{1.g_0(x)h(x) -
g_0(x).0}\over{1^2}}= g_0(x)h(x).
\end{equation}
Computing the score function for $\eta=0$ we get
\begin{equation}
{\dot l}_0(x) = {d\over d\eta} \log(g_\eta (x))\Big|_{\eta=0}=
{{{d\over d\eta}g_\eta (x)\Big|_{\eta=0}}\over{g_0 (x)}}=h(x).
\end{equation}
By checking the conditions of Proposition \ref{prop:3.1} it follows that
this one dimensional model is regular, with score function $h$ at $\eta=0$.
i.e.
at $G_0$. By Lemma \ref{expscore} all score functions at $\eta=0$ satisfy
$\int h(t)g_0(t)dt=0$. Hence the tangent set at $G_0$ is equal to
\begin{equation}
\dot {\cal P}^0 =
\{h\in {\cal L}_2(G_0): \int h(t)g_0(t)dt=0\} ={\cal L}_2^0
(G_0) .
\end{equation}
Since  the tangent set is already a closed linear space here we also have
\begin{equation}
\dot {\cal P} ={\cal L}_2^0(G_0) .
\end{equation}

Next consider $\cal G$ to be the family of {\em all probability measures
dominated by $\mu$ that are symmetric around $\nu_0$}.
For each $h$ in ${\cal L}_2(G_0)$ that is symmetric around $\nu_0$,
 the construction above yields a one dimensional regular submodel of
 $\cal G$. On the other hand it is easy to see that all score functions
for this model have to be symmeric around $\nu_0$. Hence here we have
\begin{eqnarray}
\lefteqn{\dot {\cal P} =\dot {\cal P}^0
= \{h\in {\cal L}_2(G_0): \int h(t)g_0(t)dt=0,}\\
&&\quad\quad h(\nu_0-x)=h(\nu_0 +x), a.e.\}.\nonumber
\end{eqnarray}
}
\end{example}

\begin{example}[Symmetric Location]\label{exam:5.2}
{\rm Let  $\Pp=\{P_{\nu,G}:\, \nu\in{\mathbb R}, G\in\Gg\}$  with
$P_{\nu,G}$  the distribution with density $g(\cdot-\nu)$;  here
$g(-x)=g(x)$  and  $g$  has finite Fisher information for location
(cf. (\ref{3.12})).

For this model the submodels ${\cal P}_1$ and ${\cal P}_2$ are defined by
\begin{eqnarray}
{\cal P}_1 &=& \{P_{\nu,G_0} : \nu\in {\mathbb R}\},\\
{\cal P}_2 &=& \{P_{\nu_0,G} : G\in {\cal G},\}
\end{eqnarray}
where ${\cal G}$ denotes the family of distribution functions of
distributions which are symmetric around zero and have finite Fisher
information for location.
We can compute the tangent space $\dot{\cal P}_2$ by the same construction
as in the previous example for the symmetric distributions, see Exercise
\ref{exer:5.1}.
It then turns out that the efficient influence function for estimating $\nu$
equals $\stel_1=\dot l_1$ and that the information lower bound equals
\begin{equation}
(\ex {\stel_1}^2)^{-1} = I(g_0)^{-1},
\end{equation}
where $I(g)=\int(g'/g)^2g$ is the Fisher information for location.
The lower bounds for estimating $\nu$ with $G$ known or unknown coincide.
}
\end{example}

\section{Exercises Chapter {\ref{chap:6}}}

\begin{exercise}[Symmetric Location]\label{exer:5.1}
{\rm
Prove the following statements for the symmetric location model of
Example \ref{exam:5.2} (see Example 3.2.4 of BKRW (1993))
\begin{eqnarray}
\dP&=&[\del_1]+\dP_2,\nonumber\\
\del_1(x)&=&-\,{g'_0\over g_0}(x-\nu_0),\label{5.20}\\
\dP_2&=&\{h:h(\nu_0-x)=h(\nu_0+x), \ \int h g_0 = 0, \ \int h^2g_0<\infty\},\nonumber\\
 \del_1&\perp&\dP_2.\nonumber
\end{eqnarray}
Determine  $\stel_1$ and the information lower bound for estimating $\nu$.}
 \hfill $\Box$
\end{exercise}

\begin{exercise}[Linear Regression]\label{exer:5.2}
{\rm Compute the efficient influence function for estimating
$\nu$  within the semiparametric linear regression model of
Example \ref{exam:1.2}; see also Exercise \ref{exer:3.2}} \hfill
$\Box$
\end{exercise}

\chapter{Convolution Theorem}\label{chap:7}

In general, if we  don't impose restrictions on the class of
estimators considered, the lower bound for estimating a parameter
will be equal to zero. For instance an estimator with constant
value $\theta_0$ has zero loss under the parameter value
$\theta_0$. Under the other parameter values this estimator has a
nonvanishing bias, independent of the sample size $n$. To avoid
this problem we can either restrict the class of estimators or
consider an average loss instead of the loss at a fixed parameter
value. Indeed, we have considered the average distribution of an
arbitrary estimator in the spread inequality. The approach of
restriction is used in the Convolution Theorem \ref{thm:4.1} below
where the estimators are required to be uniformly regular.

Suppose that, instead of fixed deterministic, the value $\theta$
is random, that the resulting random variable is denoted by
$\vart$ as before, and that it has a probability density $w$.
Further assume that, given $\vartheta=\theta$ the observation $X$
has probability density $p(\cdot\mid\theta)$. The squared error
loss of an estimator $T=t(X)$ is then equal to
\begin{eqnarray}
\lefteqn{\ex (T-\vart)^2=\int\ex ((T-\theta)^2|\vartheta=\theta)
w(\theta)d\theta}
\nonumber\\
&=&\int\Big(\int (t(x)-\theta)^2p(x\mid\theta)d\mu(x)\Big)
w(\theta)d\theta.
\end{eqnarray}
So averaging over the parameter value can be viewed as considering the
parameter random just as we have done in the preceding chapters. Before we
state the Convolution Theorem we will adopt this approach in some
heuristics.

 Consider the normal shift experiment of Example \ref{exam:3.1}.
We will choose a weight function for $t$, that is, we will view
$t$ as random. To that end we rename $t$ into $U$ with $U$
normally distributed with mean vector $0$ and nonsingular
covariance matrix $\Sigma$. Since  the conditional  distribution
of $X$ given $U=u$ is $\Nn(Iu,I)$, we obtain writing $X=IU+V,
V\sim {\cal N}(0,I)$ with $U$ and $V$ independent,
\begin{equation}\label{4.1}
{X \choose U} \sim \Nn \bigg( 0,
\begin{pmatrix}
   I \Sigma I + I & I \Sigma \\
   \Sigma I       & \Sigma
\end{pmatrix}
\bigg)
\end{equation}
and hence
\begin{equation}\label{4.2}
{X \choose CX-U} \sim \Nn \bigg( 0,
\begin{pmatrix}
  I \Sigma I + I & 0 \\
  0              & C
\end{pmatrix}
\bigg), \quad C = \Sigma I (I \Sigma I + I)^{-1}.
\end{equation}
The proofs of (\ref{4.1}) and (\ref{4.2}) are given in note \ref{note:7}.
 Note that (\ref{4.2}) implies that $CX-U$ and any function of $X$ are
independent.
Let  $T=t(X)$ be an estimator of $u$.
Writing
\begin{equation}\label{4.3}
T-U=(t(X)-CX)+(CX-U)
\end{equation}
we see that $T-U$  is the sum of two independent terms, namely  $CX-U$  and
 $t(X)-CX$.
Consequently, the distribution of  $T-U$  is the convolution of  $\Nn(0,C)$  and
another distribution.
Moreover, $T=t(X)=CX$  is the best estimator of  $u$  here.
For the particular choice  $\Sigma=\si^2I^{-1}$  we obtain
 $C=\si^2/(\si^2+1)I^{-1}$.
Consequently, as  $\si^2\to\infty$  the distribution of  $T-U$  is the
convolution of  $\Nn(0,I^{-1})$  and another distribution, and the best
estimator
is  $T=I^{-1}X$.

This argument is based on (\ref{4.2}) and (\ref{4.3}), and it
shows the origin of the convolution structure in a normal shift
limit experiment. In view of (\ref{3.10}) and (\ref{3.11}) of
Example \ref{exam:3.1} it is clear that LAN should imply a
convolution structure. A precise formulation  of this has been
given by H\'ajek (1970). This result is called the H\'ajek-Le Cam
convolution theorem. A proof along the lines of the above argument
may be found in Chapter \ref{chap:3} of Van den Heuvel (1996).
However, we will formulate and prove it along classical lines for
regular parametric models and with  $\nu:\, \Pp\to{\mathbb R}^m$
as parameter of interest. We will metrize  $\Pp$  with the total
variation distance, see Section \ref{tv}.

\begin{definition}[Regularity]\label{defn:4.1}
{\rm
The sequence of estimators $\{T_n\}$ with $T_n=t_n(\Xv)$ is {\em uniformly
regular} in estimating $\nu(P)$ if $\Big\{\L_P\Big({\sqrt
n}(T_n-\nu(P))\Big)\Big\} = \Big\{\L_P(Z_n)\Big\}$ converges uniformly on compact subsets
$K$ of $\Pp$. In other words, there exists a family $\{\L_P(Z)\}$
of distributions on ${\mathbb R}^m$ such that for every bounded,
continuous function $g$ on ${\mathbb R}^m$ and for every compact
$K\subset\Pp$,
\begin{equation}\label{4.4}
\limy\,\sup_{P\in K} |\expar_Pg(Z_n)-\expar_Pg(Z)|=0.
\end{equation}
}
\end{definition}

\begin{example}
{\rm Returning to the example of an estimator $T_n$ of $\theta$
that has a constant value $\theta_0$, we see that such an
estimator is not regular. Even for fixed $P$ the random variable
$Z_n=\sqrt{n}(T_n-\theta) = \sqrt{n}(\theta_0-\theta)$ does not
converge in distribution for all $\theta$ values. By imposing
uniform regularity on the estimators these constant estimators
will therefore be excluded. }
\end{example}

\begin{definition}[Linearity]\label{defn:4.2}
{\rm The sequence of estimators $\{T_n\}$ is {\em uniformly
asymptotically linear} in the influence function $\psi:\,
\Xx\times\Pp\to{\mathbb R}^m$ with $\expar_P|\psi(X,P)|^2<\infty,\
\expar_P\psi(X,P)=0$, if
\begin{equation}\label{4.5}
\limy\,\sup_{P\in K} P\Big({\sqrt n}\Big|
T_n-\{\nu(P)+\ene\nn\psi(X_i,P)\Big\}\Big|>\varepsilon\Big)=0
\end{equation}
for every positive $\varepsilon$ and all compacts $K\subset\Pp$.
}
\end{definition}

Note that this definition extends the pointwise asymptotic
linearity introduced in Definition \ref{newdefn:4.1}. The map
$\nu:\,\Pp\to{\mathbb R}^m$ and the parametrization $\theta\to
P_\theta$ together constitute a map $q:\, \Theta\to{\mathbb R}^m$.
As before, we will assume that $q$ is differentiable with ${\dot
q}(\theta)$ the $m\times k$-matrix of partial derivatives
$\partial q_i(\theta)/\partial\theta_j, \ i=1,\dots,m, \
j=1,\dots,k$.

\bigskip

These definitions enable us to formulate the convolution theorem and the
resulting concept of efficiency as follows.

\begin{theorem}[Convolution Theorem]\label{thm:4.1}
 Let $\{T_n\}$ be a uniformly regular sequence of estimators of
$q(\theta)=\nu(P_\theta)$ in the regular parametric model $\Pp$.
Let $A(\theta)$ be an $\ell\times k$ matrix which is continuous in
$\theta$ and let $S_n(\theta)$ be the score function defined in
{\rm(\ref{3.5})}. If ${\dot q}(\theta)$ is continuous in $\theta$,
then uniformly on compact subsets $K$ of $\Theta$
\begin{equation}\label{4.9}
{{\sqrt
n}(T_n-\nu(P_\theta)-\ene\nn\tell(X_i;P_\theta\mid\nu,\Pp))\choose
A(\theta)S_n(\theta)}{\buildrel\D\over\rightarrow_\theta} \
{\Del_\theta\choose S_\theta},
\end{equation}
where the random $m$-vectors $\Del_\theta$ and $\ell$-vectors
$S_\theta$ are independent. Consequently, the limit distribution
of ${\sqrt n}(T_n-\nu(P_\theta))$ under $P_\theta$ is the
convolution of the normal $\Nn(0,I^{-1}(P_\theta\mid\nu,\Pp))$
distribution and the distribution of $\Del_\theta$.
\end{theorem}

\begin{remark}[Efficiency]\label{rem:4.1}
{\rm Adding up the two components of the vectors in (\ref{4.9})
with $A(\theta)={\dot q}(\theta)I^{-1}(\theta)$, we see that
\begin{equation}\label{4.10}
{\sqrt
n}\Big(T_n-\nu(P_\theta)\Big){\buildrel\D\over\rightarrow_\theta}\ S_\theta+\Del_\theta
\end{equation}
with $S_\theta$ normal
$\Nn\Big(0,I^{-1}(P_\theta\mid\nu,\Pp)\Big)$ and with $S_\theta$
and $\Del_\theta$ independent. In terms of covariance matrices
this independence yields
\begin{equation}\label{4.11}
\var(S_\theta+\Del_\theta)=\var S_\theta+\var \Del_\theta\ge\var
S_\theta
\end{equation}
in the partial ordering of matrices defined by $A\ge B$ iff $A-B$
is positive semidefinite. This shows that a sequence $\{T_n\}$ of
estimators is asymptotically optimal if and only if $\Del_\theta$
is degenerate at $0$. In view of (\ref{4.9}) this is equivalent to
uniform asymptotic linearity of $\{T_n\}$ in the efficient
influence function $\tell(\cdot;P_\theta\mid\nu,\Pp)$. Such
estimators are called {\em uniformly asymptotically efficient} and
this also explains the terminology of Definition \ref{newdefn:4.2}
once more. }
\end{remark}

\noindent{\bf Proof of Theorem \ref{thm:4.1}}

The concept of tightness and Prohorov's theorem which are needed in this proof
are given in Section \ref{proh}.

Consider a sequence $\{\theta_n\}$ converging to a fixed $\theta$.
Denote
\begin{equation}\label{4.a}
(U_n,V_n)=\Big({\sqrt n}\Big(T_n-q(\theta_n)\Big), \
S_n(\theta_n)\Big).
\end{equation}
By (\ref{4.4}) and (\ref{3.8}) $\{U_n\}$  and  $\{V_n\}$  are marginally
convergent in
distribution, hence marginally tight, and consequently jointly tight.
In view of Prohorov's theorem any subsequence  $\{n'\}$   of  $\{n\}$  has  a
further subsequence  $\{n''\}$  such that
\begin{equation}\label{4.b}
(U_{n''},V_{n''})\to (U,V).
\end{equation}
By an abuse of notation we will write  $n$  instead of  $n''$.
Let
$$W_n=L_n\Big(\theta_n+{t\over{\sqrt n}}\Big)-L_n(\theta_n),$$
and note that by (\ref{3.7}) of the LAN-Theorem \ref{thm:3.1}
\begin{equation}\label{4.c}
(U_n,W_n)\to \Big(U,t^TV-\halfe t^TI(\theta)t\Big)=(U,W).
\end{equation}
By the continuous differentiability of  $q(\cdot)$  and the
uniform regularity of  $\{T_n\}$  we have, for all  $a\in{\mathbb
R}^m$,
\begin{eqnarray}
\lefteqn{\limy \expar_{\theta_n+t/{\sqrt n}}\exp\{ia^TU_n\}}\nonumber\\
&=&\limy \expar_{\theta_n+t{\sqrt n}}\exp\Big\{i a^T{\sqrt
n}\Big(T_n-q(\theta_n+t/{\sqrt
n})\Big)+i a^T{\dot q}(\theta)t\Big\}\label{4.d}\\
&=&\expar\exp\Big\{i a^TU+i a^T{\dot q}(\theta)t\Big\}.\nonumber
\end{eqnarray}
On the other hand we have
\begin{equation}\label{4.e}
\expar_{\theta_n+t/{\sqrt n}}\exp\{i a^TU_n\}=
\expar_{\theta_n}\exp\{i a^TU_n+W_n\land\la\}+R_{n\la}+R_n
\end{equation}
with
\begin{eqnarray}
\Big|R_{n\la}\Big|&=&
\Big|\expar_{\theta_n}\exp\{ia^TU_n\}\Big(e^{W_n}
-e^{W_n\land\la}\Big)\Big|\le 1-\expar_{\theta_n}e^{W_n\land\la},\label{4.f}\\
\Big|R_n\Big|&=& \Big|\expar_{\theta_n+t/{\sqrt n}}\exp\Big\{i
a^TU_n\Big\}{\been}_{[L_n(\theta_n)=-\infty]}\Big|\label{4.g}\\
&&\quad\quad\leq P_{\theta_n+t/{\sqrt
n}}\Big(L_n(\theta)=-\infty\Big).\nonumber
\end{eqnarray}
By contiguity (Corollary \ref{cor:3.2})
 $P_{\theta_n}\Big(L_n(\theta_n)=-\infty\Big)=0$
implies
$$P_{\theta_n+t/{\sqrt n}}\Big(L_n(\theta_n)=-\infty\Big)\to 0$$
and hence
\begin{equation}\label{4.h}
|R_n|\to 0, \quad {\rm as}\  n\to\infty.
\end{equation}
Note that (\ref{3.8}), (\ref{4.c})  and (\ref{3.10}) imply  $\ex e^W=1$.
Consequently, (\ref{4.f}) yields
\begin{equation}\label{4.i}
\lim_{\la\to\infty}\limy
|R_{n\la}|\le\lim_{\la\to\infty}1-\ex e^{W\land\la}=0,
\end{equation}
which together with (\ref{4.h}) shows
\begin{equation}\label{4.j}
\limy \expar_{\theta_n+t/{\sqrt n}}\exp\{i a^TU_n\}
=\lim_{\la\to\infty}\ex\exp\{i a^TU+W\land\la\}=\ex\exp\{i
a^TU+W\}.
\end{equation}

\noindent Combining (\ref{4.d}) and (\ref{4.j}) we obtain
\begin{equation}\label{4.k}
\ex\exp\{ia^TU+t^TV-\halfe t^TI(\theta)t\} =\exp\{ia^T{\dot
q}(\theta)t\} \ex\exp\{ia^TU\}, \ t\in{\mathbb R}^k.
\end{equation}
Since both sides of (\ref{4.k}) are analytic functions (cf. e.g.
Theorem 2.9, p.52 of Lehmann (1959)), (\ref{4.k}) also holds by
analytic continuation (cf. e.g. for   $t^T=-ia^T{\dot
q}(\theta)I^{-1}(\theta)+ib^TA(\theta)$  with
 $b\in{\mathbb R}^m$  arbitrary.
This results in
\begin{eqnarray}\label{4.l}
\lefteqn{\ex\exp\{ia^T(U-{\dot q}(\theta)I^{-1}(\theta)V+ib^TA(\theta)V\}}\\
&=&\ex\exp\{ia^TU+\halfe a^T{\dot q}(\theta)I^{-1}(\theta){\dot
q}^T(\theta)a-\halfe b^TA(\theta)I(\theta)A^T(\theta)b\}\nonumber
\end{eqnarray}
and in particular for  $b=0$,  in
\begin{equation}\label{4.m}
\ex\exp\{ia^T(U-{\dot q}(\theta)I^{-1}(\theta)V)\}
=\ex\exp\{ia^TU+\halfe a^T{\dot q}(\theta)I^{-1}(\theta){\dot
q}^T(\theta)a\}.
\end{equation}
Combining (\ref{4.l}) and (\ref{4.m}), using (\ref{4.b}) and
recalling the notation in (\ref{4.a}) and (\ref{4.9}), we see that
the characteristic function of   $\Big(\Del^T_\theta,
S^T_\theta\Big)^T$ at  $\Big(a^T,b^T\Big)^T$  can be written as
\begin{equation}\label{4.n}
\ex\exp\Big\{ia^T\Del_\theta+ib^TS_\theta\Big\}
=\ex\exp\Big\{ia^T\Del_\theta\Big\} \exp\Big\{-\halfe b^T
A(\theta)I(\theta)A^T(\theta)b\Big\},
\end{equation}
which implies the independence of  $\Del_\theta$  and  $S_\theta$.

Note that (\ref{4.m}) shows that the distribution of $\Del_\theta$
is the same for all subsequences  $\{n''\}$. Since the
distribution of  $S_\theta$  does not depend on the subsequence
either and since  $\Del_\theta$  and  $S_\theta$  are independent
we obtain
\begin{equation}\label{4.o}
{{\sqrt
n}(T_n-\nu(P_{\theta_n})-\ene\nn\tell(X_i;P_{\theta_n}\mid\nu,\Pp))\choose
A(\theta_n)S_n(\theta_n)}{\buildrel\D\over\rightarrow_{\theta_n}}
\ {\Del_\theta\choose S_\theta}.
\end{equation}
However, (\ref{4.m}) also shows that the distribution of
$\Del_\theta$  is continuous in  $\theta$. Combining this with the
continuity of the distribution of  $S_\theta$,  with the
independence of  $\Del_\theta$  and  $S_\theta$  and with
(\ref{4.o}) we obtain (\ref{4.9}). The second part of Theorem
\ref{thm:4.1} has been proved in Remark \ref{rem:4.1}. \hfill
$\Box$

\bigskip

For nonparametric (and semiparametric) models the previous
classical Convolution Theorem can be extended as follows (cf. the
LAS Theorem \ref{thm:new6.2}).

\begin{theorem}[Convolution Theorem]\label{thm:5.2}
 Let  $\nu:\,\Pp\to{\mathbb R}^m$  be pathwise differentiable on  $\Pp$  at  $P_0$
and let  $T_n$  be a locally regular estimator of  $\nu$  (on all
regular parametric submodels $\Qq=\{P_\theta:\,\theta\in{\mathbb
R}^k, \, |\theta|<\varepsilon\}$ of
 $\Pp$  containing  $P_0$), that is,
\begin{equation}\label{5.9}
{\sqrt
n}(T_n-\nu\Big(P_{\theta_n})\Big){\buildrel\D\over\longrightarrow_{\theta_n}
Z}
\end{equation}
whenever  $\theta_n={\O}(n^{-1/2})$. If
\begin{equation}\label{5.10}
[\dun]\subset\odP,
\end{equation}
the closure of the tangent set  $\dP^0$,  then for any  $h\in(\odP)^\ell$
\begin{equation}\label{5.11}
{{\sqrt
n}\Big(T_n-\nu(P_0)-\ene\nn\tel(X_i,P_0\mid\nu,\Pp)\Big)\choose{1\over\sqrt
n}\nn h(X_i)}{\buildrel\D\over\longrightarrow_{P_0}{\Del_0\choose W_0}},
\end{equation}
where  $\Del_0$  and  $W_0$  are independent.
Moreover,  $\Del_0$  is degenerate at  $0$  iff  $T_n$  is regular and
asymptotically linear with influence function  $\tel(\cdot,P_0\mid\nu,\Pp)$
under  $P_0$,  and  $T_n$  is called efficient then.
\end{theorem}

The proof of this Convolution Theorem is similar to the proof of
Convolution Theorem \ref{thm:4.1} and involves an approximation procedure
based on (\ref{5.10}).
It may be found in Section 3.3, pp. 64, 65, of BKRW (1993).

\bigskip

As we have seen, LAS Theorems are valid for any sequence of
estimators and Convolution Theorems for regular sequences only.
For the still more specific regular sequences of asymptotically
linear estimators we obtain via still a different approach
asymptotic optimality as follows.

\begin{theorem}[Regular Linear Estimators]
 Let  $\{T_n\}$  be a regular sequence of asymptotically linear estimators of
 $\nu:\,\Pp\to{\mathbb R}^m$  with influence function  $\psi$  at  $P_0$.
Then  $\nu$  is pathwise differentiable at  $P_0$  with derivative  $\dun$  and
\begin{equation}\label{5.12}
\psi-\dun=\psi-\tel\perp\dP.
\end{equation}
In particular, the limit covariance matrix $<\psi,\psi^T>_0$
equals at least $I^{-1}(P_0\mid\nu,{\cal P})$.
Moreover, $T_n$  is efficient at  $P_0$ in the sense of
\begin{equation}\label{new7.32}
<\psi,\psi^T>_0=I^{-1}(P_0\mid\nu,{\cal P})
\end{equation}
iff
\begin{equation}\label{5.13}
\psi\in\dP^m
\end{equation}
and then  $\psi=\tel$.
\end{theorem}

\noindent{\bf Sketch of Proof.}

Let  $\{P_\eta\}$  be a path through  $P_0$  with score function  $h$  and
choose  $\eta_n=t_n^{-1/2}$.
Regularity and linearity of  $T_n$  imply
\begin{equation}\label{5.2.a}
{\sqrt n}(T_n-\nu(P_{\eta_n})){\buildrel\D\over\longrightarrow_{P_{\eta_n}}}
\Nn(0,\langle\psi,\psi^T\rangle_0\big).
\end{equation}
On the other hand linearity of  $T_n$,  LAN, and Le Cam's third lemma
(e.g. lemma A.9.3 pp. 503, 504 in BKRW (1993))  yield
\begin{equation}\label{5.2.b}
{\sqrt n}(T_n-\nu(P_0)){\buildrel\D\over\longrightarrow_{P_{\eta_n}}}
\Nn(t\langle\psi,h\rangle_0, \ \langle\psi,\psi^T\rangle_0\big).
\end{equation}
Combining (\ref{5.2.a}) and (\ref{5.2.b}) we obtain
\begin{equation}\label{5.2.c}
\nu(P_{\eta_n})=\nu(P_0)+\eta_n\langle\psi,h\rangle_0+o(\eta_n).
\end{equation}
Comparing (\ref{5.2.c}) to (\ref{5.1}) we  arrive at (\ref{5.12}).
Because of $\tilde\ell\in\dP$, (\ref{5.12}) implies that
$<\psi,\psi^T>_0-<\tilde\ell,\tilde\ell^T>_0$ is positive
semidefinite. Finally,
 note that (\ref{5.12}) and (\ref{5.13}) together imply
$\psi=\tel$.

\hfill $\Box$

\chapter{Construction of Estimators}\label{chap:8}

In the preceding two chapters bounds have been constructed on the asymptotic
behavior of estimators of pathwise differentiable parameters
 $\nu:\,\Pp\to{\mathbb R}^m$
in non- and semiparametric models.
A complete theory on construction of efficient estimators in these models,
that is, estimators attaining these bounds in the limit, is lacking.

For regular parametric models Le Cam has proved that efficient
estimators exist indeed. In fact, he gave a construction as
follows. First, he constructed a  $\sqrt n$-consistent estimator
$\tt_n$  of the parameter  $\theta$  in the regular parametric
model
 $\Pp=\{P_\theta:\,\theta\in\Theta\}, \ \Theta\subset{\mathbb R}^k$;  see Theorem 2.5.1 of
BKRW (1993). Then, he discretized  $\tt_n$  to an estimator
$\theta^*_n$  taking its values in a grid in  ${\mathbb R}^k$ with
mesh width  $cn^{-1/2}$,  such that
 $|\theta^*_n-\tt_n|\le Cn^{-1/2}$  a.s.
Motivated by a Newton-Raphson procedure he defined the one-step efficient
estimator
 $\hta_n$  by
\begin{equation}\label{6.1}
\hta_n=\theta^*_n+\ene\nn\tel(X_i;\theta^*_n),
\end{equation}
where the efficient influence function  $\tel(\cdot;\theta)$  is
defined by (cf. (\ref{new4.11}))
\begin{equation}\label{6.2}
\tel(x;\theta)=\tel(x;P_\theta\mid\theta,\Pp)=I^{-1}(\theta)\del(x;\theta),
\quad x\in\Xx.
\end{equation}
\begin{theorem}[Discretization]\label{thm:6.1}
Let  $\Pp=\{P_\theta:\,\theta\in\Th\}, \ \Th\subset{\mathbb R}^k$,
be a regular parametric model and let  $\tt_n$  be a locally
$\sqrt n$-consistent estimator of
 $\theta$  based on  $\Xv$,  that is, for  every  $\theta\in{\mathbb R}^k$  and every sequence
 $\{\theta_n\}$  with  $\theta_n=\theta+{\O}(n^{-1/2})$
\begin{equation}\label{6.3}
\lim\limits_{M\to\infty}\limsup_{n\to\infty}P_{\theta_n}\Big(|{\sqrt
n}(\tt_n-\theta_n)|>M\Big)=0.
\end{equation}
If  $\theta^*_n$  is a discretized version of  $\tt_n$,  then
$\hta_n$  as defined by {\rm(\ref{6.1})} and {\rm(\ref{6.2})} is
locally efficient in the sense (cf. {\rm Definition
\ref{newdefn:4.3}} and {\rm Remark \ref{rem:4.1}})
\begin{equation}\label{6.4}
\limy P_{\theta_n}\Big(|{\sqrt
n}(\hta_n-\{\theta_n+\ene\nn\tel(X_i;\theta_n)\})|>\varepsilon\Big)=0,
\quad \varepsilon>0.
\end{equation}
\end{theorem}

\noindent
{\bf Proof.}

Fix  $\varepsilon>0$  and  $\delta>0$. In view of (\ref{6.3}) and the
definition of  $\theta^*_n$  we may choose  $M$ sufficiently large
such that for  $n$  large
\begin{equation}\label{6.a}
P_{\theta_n}\left(|\sqn(\theta^*_n-\theta_n)|>M\right)<\delta.
\end{equation}
Let  $B_n$  be the collection of grid points at distance at most
$Mn^{-1/2}$ from  $\theta_n$. By (\ref{6.a}) we have
\begin{eqnarray}
\lefteqn{P_{\theta_n}\Big(|\sqn(\hta_n-\{\theta_n+\ene\nn\tel(X_i;\theta_n)\})|
>\varepsilon\Big)\le P_{\theta_n}(\theta^*_n\not\in B_n)}\nonumber\\
&+&P_{\theta_n}\Big(|\sqn(\theta^*_n-\theta_n+\ene\nn\tel(X_i;\theta^*_n)
 -\ene\nn\tel(X_i;\theta_n))|>\varepsilon,   \theta^*_n\in B_n\Big)\label{6.b}\\
&\le&\delta+\sum_{\theta'_n\in
B_n}P_{\theta_n}\Big(|\sqn(\theta'_n-\theta_n+\ene\nn\tel(X_i;\theta'_n)-\ene\nn\tel(X_i
;\theta_n))|>\varepsilon \Big).\nonumber
\end{eqnarray}
Since the number of points in  $B_n$  is bounded (uniformly in  $n$) and since
 $\delta$  is arbitrary, it suffices to prove that for every sequence
 $\{\theta'_n\}$  with  $\theta'_n=\theta+{\O}(n^{-1/2})$  the probability at the
right-hand side of (\ref{6.b}) converges to  $0$. However, this
holds in view of the smoothness property (\ref{3.9}), the
continuity of the Fisher information  $I(\cdot)$  and the
definition of
 $\tel(\cdot;\cdot)$  in (\ref{6.2}). \hfill$\Box$

\bigskip

The estimator  $\tt_n$  is called a preliminary estimator.
The proof makes clear that the discretization is used to force (a kind of)
independence between the preliminary estimator and the observations.
Another technique to force this is by splitting the sample into two parts.
One part is used for the preliminary estimator and  an average is constructed
of
the efficient influence function at the observations from the other part.
Interchanging the roles of these two parts we obtain an estimator that is
asymptotically linear in the observations of the first part.
An appropriate convex combination of these  two estimators is efficient.

\begin{theorem}[Sample Splitting]\label{thm:6.2}
Let  $\Pp=\{P_\theta:\,\theta\in\Th\}$  be a regular parametric
model and let
 $\tt_n=t_n(\Xv)$ be a locally  ${\sqrt n}$-consistent estimator of  $\theta$.
With  $\la_n\in{\mathbb N}$  we define
$$\tt_{n1}=t_{\la_n}(X_1,\dots,X_{\la_n}), \
\tt_{n2}=t_{n-\la_n}(X_{\la_n+1},\dots,X_n).$$
If
\begin{equation}\label{6.5}
\la_n/n\to\la\in(0,1),
\end{equation}
then the estimator
\begin{equation}\label{6.6a}
\hta_n={\la_n\over
n}\Big(\tt_{n2}+{1\over\la_n}\sum_{i=1}^{\la_n}\tel(X_i;\tt_{n2})\Big)
+
{n-\la_n\over n}\Big(\tt_{n1}+{1\over
n-\la_n}\sum_{i=\la_n+1}^n\tel(X_i;\tt_{n1})\Big)
\end{equation}
is locally efficient.
\end{theorem}

\noindent
{\bf Proof.}

First note that for (\ref{6.4}) to hold with  $\hta_n$  as in (\ref{6.6a}),
it suffices
to prove
\begin{equation}\label{6.2.a}
\limy
P_{\theta_n}\Big(|{\sqrt{\la_n}}(\tt_{n2}-\theta_n+{1\over\la_n}
\sum_{i=1}^{\la_n}\tel(X_i;\tt_{n2})
 -\
{1\over\la_n}\sum_{i=1}^{\la_n}\tel(X_i;\theta_n))|>\varepsilon\Big)=0,
\ \varepsilon>0.\nonumber
\end{equation}
Let  $M$  be sufficiently large and let  $B_n$  be the ball of radius
 $M\la_n^{-1/2}$  about  $\theta_n$.
Then the probability in (\ref{6.2.a}) is bounded from above by
\overfullrule=0pt
\begin{eqnarray*}
\lefteqn{P_{\theta_n}(\tt_{n2}\not\in B_n)+
E_{\theta_n}\Big({\been}_{B_n}(\tt_{n2})
P_{\theta_n}(|{\sqrt{\la_n}}
(\tt_{n2}-\theta_n+}\\
&&\quad\quad+
{1\over\la_n}\sum_{i=1}^{\la_n}\tel(X_i;\tt_{n2})-{1\over\la_n}
\sum_{i=1}^{\la_n}\tel(X_i;\theta_n))|>\varepsilon\mid\tt_{n2}\Big)).
\end{eqnarray*}
But, the conditional probability converges to  $0$,  for  $\tt_{n2}\in
B_n$,  by the
same argument as at the end of the proof of Theorem \ref{thm:6.1}
 in view of (\ref{6.5}) and the
independence of  $\tt_{n2}$  and  $(X_1,\dots,X_{\la_n})$.
\hfill $\Box$

Of course, both discretization and sample splitting are artificial
techniques used only to make proofs work. The estimator  $\hta_n$
from (\ref{6.1}) with  $\theta^*_n$  replaced by $\tt_n$ itself
would be more natural, since it is very close to a one-step
Newton-Raphson approximation of the maximum likelihood equation
\begin{equation}\label{6.6}
\nn\del(X_i;\theta)=0.
\end{equation}
Indeed, for  $k=1$  with  $\tt_n$  close to a root of this equation one-step
Newton-Raphson results in
\begin{equation}\label{6.7}
\tt_n-\ene\nn\del(X_i;\tt_n)\Big\{\ene\nn{\partial\over\partial\theta}
\del(X_i;\tt_n)\Big\}^{-1}.
\end{equation}
By the law of large numbers and under extra smoothness conditions on
 $\del(\cdot;\theta)$  the denominator should converge to
\begin{eqnarray}
\expar_\theta{\partial\over\partial\theta}\del(X;\theta)&=&\int\Big\{{\partial
\over\partial\theta}
(\del(x;\theta)p(x;\theta))-\del(x;\theta){\partial\over\partial\theta}p(x;\theta)\Big\}
d\mu(x)\quad
\label{6.8}\\
&=&{\partial\over\partial\theta}
\expar_\theta\del(X;\theta)-\int\del^2(x;\theta)p(x;\theta)d\mu(x)=-I(\theta)\nonumber
\end{eqnarray}
and hence (\ref{6.7}) resembles (\ref{6.1}). Under additional
regularity conditions it can be shown that the maximum likelihood
estimator itself is efficient; see pp. 500--504 of Cram\'er (1946)
and Le Cam (1970).

\bigskip

Here we will apply the sample splitting technique to construct
efficient estimators of the Euclidean parameter  $\theta$  in the
semiparametric model
\begin{equation}\label{6.9}
\Pp=\{P_{\theta,G}:\,\theta\in\Th, \ G\in\Gg\},
\end{equation}
where  $\Th\subset{\mathbb R}^k$  is open and  $\Gg$  is general.
Under the conditions of Theorem \ref{thm:5.3} with  $\nu=\theta$
the efficient influence function is well-defined and we will write
(cf. (\ref{5.18}))
\begin{equation}\label{6.10}
\tel(x;\theta;G)=(\ex{\stel_1}{\stel_1}^T)^{-1}{\stel_1}.
\end{equation}
Since  $\tel$  depends on  $G$  and since  $G$  is unknown, the estimator
(\ref{6.6a})
cannot be used here.
We will introduce an extra splitting step to handle this problem under the
assumption that there exists an estimator
 $\tel_n(\cdot;\theta;\uX)=\tel_n(\cdot;\theta;\,\Xv)$  of the efficient influence
function  $\tel(\cdot;\theta;G)$  satisfying ({\bf consistency})
\begin{equation}\label{6.11}
\int\Big|\tel_n(x;\theta_n;\uX)-\tel(x;\theta_n;G)\Big|^2dP_{\theta_n,G}(x)=o_{\theta_
n,G}(1)
\end{equation}
and ({\bf $\sqn$-unbiasedness})
\begin{equation}\label{6.12}
\sqn\int\tel_n(x;\theta_n;\uX)dP_{\theta_n,G}(x)=o_{\theta_n,G}(1),
\end{equation}
for all  $(\theta,G)\in\Th\times\Gg$  and all sequences
$\{\theta_n\}$  with
 $\theta_n=\theta+{\O}(n^{-1/2})$.
To interpret (\ref{6.12}) recall that
$$\int\tel(x;\theta_n;G)dP_{\theta_n,G}(x)=0.$$
Furthermore, we assume the existence of a  $\sqn$-consistent
preliminary estimator of  $\theta$ denoted by $\tt_n=t_n(\Xv)$
satisfying (\ref{6.3}) under $P_{\theta_n,G}$ for every
$G\in\Gg$. Note that in the parametric case such a
$\sqn$-consistent estimator always exists according to Le Cam (cf.
Exercise \ref{exer:6.1}), but in a semiparametric model we do not
know this beforehand (cf. Exercise \ref{exer:6.3}).

Let  $\{\la_n\}, \{\mu_n\}$,  and  $\{\nu_n\}$  be sequences of integers with
$0<\la_n<\mu_n<\nu_n<n$  and
\begin{equation}\label{6.13}
{\la_n\over n}\to\la, \ {\mu_n\over n}\to\mu, \ {\nu_n\over n}\to\nu,
\end{equation}
with  $0<\la<\mu<\nu<1$.
We define (cf. Theorem 6.2)
\begin{equation}\label{6.14}
\tt_{n1}=t_{\la_n}(X_1,\dots,X_{\la_n}), \
\tt_{n2}=t_{\nu_n-\mu_n}(X_{\mu_n+1},\dots,X_{\nu_n}),
\end{equation}
\begin{eqnarray}
\tel_{n1}(x;\theta)&=&\tel_{\mu_n-\la_n}(x;\theta;X_{\la_n+1},\dots,
X_{\mu_n}),\label{6.15}\\
\tel_{n2}(x;\theta)&=&\tel_{n-\nu_n}(x;\theta;X_{\nu_n+1},\dots,X_n),
\nonumber
\end{eqnarray}
and  we notice that these four estimators are independent, because they are
based on the four independent blocks of observations into which our sample of
 $n$  observations is split up.
Analogously to (\ref{6.6a}) we define
\begin{equation}\label{6.16}
\hta_n={\mu_n\over
n}\Big(\tt_{n2}+{1\over\mu_n}\sum_{i=1}^{\mu_n}\tel_{n2}(X_i;\tt_{n2})\Big)
+{n-\mu_n\over n}\,\Big(\tt_{n1}+{1\over
n-\mu_n}\sum_{i=\mu_n+1}^n\tel_{n1}(X_i;\tt_{n1})\Big).
\end{equation}
To prove efficiency of this estimator  $\hta_n$  in the semiparametric
model (\ref{6.9}) we assume smoothness of the efficient influence functions as
follows,
\begin{equation}\label{6.17}
\sqn\Big(\theta_n-\theta+\ene\nn\tel(X_i;\theta_n;G)-\ene\nn\tel(X_i;\theta;G)\Big)=
o_{\theta,G}(1),
\end{equation}
for every  $(\theta,G)\in\Th\times\Gg$  and for all sequences
$\{\theta_n\}$  with
 $\theta_n+{\O}(n^{-1/2})$.
Note that in a least favorable regular parametric submodel this follows from
 the
LAN-Theorem \ref{thm:3.1} via (\ref{3.9}) as in the proof of Theorem
\ref{thm:6.1}.

\begin{theorem}[Semiparametric Sample Splitting]\label{thm:6.3}
In the semiparametric model {\rm(\ref{6.9})} let
$\Pp_1(G)=\{P_{\theta,G}:\, \theta\in\Th\}$ be a regular
parametric submodel for every  $G\in\Gg$  and let  $\tel$  be
well-defined by {\rm(\ref{6.10})}. Let  $\tt_n$  be a
$\sqn$-consistent estimator of  $\theta$  and let
 $\tel_n(\cdot;\cdot;\uX)$  satisfy {\rm(\ref{6.11})} and {\rm(\ref{6.12})}.
If {\rm(\ref{6.17})} holds, then the estimator  $\hta_n$  from
{\rm(\ref{6.16})} is locally uniformly efficient, i.e.
\begin{equation}\label{6.18}
\sqn\Big(\hta_n-\theta_n-\ene\nn\tel(X_i;\theta_n;G)\Big)=o_{\theta_n,G}(1),
\end{equation}
for every  $(\theta,G)\in\Th\times\Gg$  and every sequence
$\{\theta_n\}$  with
 $\theta_n=\theta+{\O}(n^{-1/2})$.
\end{theorem}

\noindent {\bf Proof.} Fix  $\theta\in\Th$  and  $G\in\Gg$. In
view of Theorem \ref{thm:3.1} and Corollary \ref{cor:3.2}
 the regularity of  $\Pp_1(G)$
implies the mutual contiguity of  $\{P_{\theta,G}^n\}$  and
$\{P_{\theta_n,G}^n\}$  for every sequence  $\{\theta_n\}$  with
$\theta_n=\theta+{\O}(n^{-1/2})$. Together with the smoothness
condition (\ref{6.17}) this shows that it suffices to prove
(\ref{6.18}) with  $\theta_n=\theta$. We will just prove
\begin{equation}\label{6.3.a}
P_{\theta,G}\Big(\Big|{\sqrt\mu_n}\Big(\tt_{n2}-\theta+{1\over\mu_n}
\sum_{i=1}^{\mu_n}\tel_{n2}(X_i;\tt_{n2})-
{1\over\mu_n}\sum_{i=1}^{\mu_n}\tel(X_i;\theta;G)\Big)\Big|>\varepsilon\Big)\to
0, \ \varepsilon>0.
\end{equation}
As in the proof of Theorem \ref{thm:6.2}, let  $M$  be
sufficiently large and let  $B_n$  be the ball of radius
$Mn^{-1/2}$  about  $\theta$. Then the probability in
(\ref{6.3.a}) is bounded from above by
\begin{eqnarray*}
&&P_{\theta,G}(\tt_{n2}\not\in B_n)
+E_{\theta,G}\Big({\been}_{B_n}(\tt_{n2})P_{\theta,G}\Big(\Big|
{\sqrt\mu_n}\Big(\tt_{n2}-\theta+{1\over\mu_n}\sum_{i=1}^{\mu_n}\tel_{n2}(X_i;
\tt_{n2})\\
&&\quad\quad\quad\quad\quad\quad-{1\over\mu_n}\sum_{i=1}^ {\mu_n}
\tel(X_i;\theta;G)\Big)\Big|>\varepsilon\mid\tt_{n2}\Big)\Big),
\end{eqnarray*}
and by the  $\sqn$-consistency of  $\tt_{n2}$,  by the independence of
 $\tt_{n2}$   and  $(X_1,\dots,$ $X_{\mu_n},\ X_{\nu_n+1},\dots,X_n)$  and
by the dominated convergence theorem it suffices to prove that for
$\theta_n=\theta+{\O}(n^{-1/2})$
\begin{equation}\label{6.3.b}
P_{\theta,G}\Big(\Big|{\sqrt\mu_n}\Big(\theta_n-\theta+{1\over\mu_n}
\sum_{i=1}^{\mu_n}\tel_{n2}(X_i;\theta_n)-{1\over\mu_n}
\sum_{i=1}^{\mu_n}\tel(X_i;\theta,G)\Big)\Big|>\varepsilon\Big)\to 0.
\end{equation}
In view of (\ref{6.17}) and (\ref{6.13}), and because of
contiguity this holds if for all $a\in{\mathbb R}^k$
\begin{equation}\label{6.3.c}
{1\over {\sqrt\mu_n}}\sum_{i=1}^{\mu_n}
a^T\Big(\tel_{n2}(X_i;\theta_n)-
\tel(X_i;\theta_n,G)\Big)=o_{\theta_n,G}(1).
\end{equation}
Computing the conditional expectation and conditional variance of
the left-hand side of (\ref{6.3.c}) given
$(X_{\nu_n+1},\dots,X_n)$,  and  using (\ref{6.12}) and
(\ref{6.11}) we see that (\ref{6.3.c}) holds indeed.\hfill $\Box$

\bigskip

\begin{remark}[Necessity]\label{rem:6.1}
{\rm
Theorem \ref{thm:6.3} may be extended to general linear estimators as in
Klaassen (1987).
Our proof is taken from this paper.
The proof as given for Theorem 7.8.1 of BKRW (1993) has been messed up at
the type setting stage but has been corrected in the 1998 edition.
In Klaassen (1987) it has also been proved under a uniform integrability
condition on the efficient influence functions that existence of an estimator
of
the efficient influence function satisfying (\ref{6.11}) and (\ref{6.12}),
 is necessary for
existence of an efficient estimator of  $\theta$. }
\end{remark}

\begin{remark}\label{rem:6.2}
{\rm
For a broad class of semiparametric models estimation of the efficient
 influence
function reduces to estimation of the score function for location
$-g'/g$ for a density  $g$. A method for this has been given by
Bickel; see Bickel (1982), Bickel and Klaassen (1986) and
Proposition 7.8.1, p. 400, of BKRW (1993). For other methods of
estimation of  $\theta$,  efficiently or just
 $\sqn$-consistently, we refer to Chapter 7 of BKRW (1993).
}
\end{remark}

\section{Exercises Chapter {\ref{chap:8}}}

\begin{exercise}[$\sqn$-consistency]\label{exer:6.1}
{\rm
In a regular parametric model  $\tt_n$  satisfies (\ref{6.3}) iff for every
 $\theta\in{\mathbb R}^k$
\begin{equation}\label{6.19}
\lim_{M\to\infty}\limsup_{n\to\infty}
P_\theta\Big(\Big|\sqn(\tt_n-\theta)\Big|>M\Big)=0.
\end{equation}
}
\end{exercise}

\begin{exercise}[Smoothness]\label{exer:6.2}
{\rm If  $\Pp_1(G)$  is regular, (\ref{6.17}) is necessary for the
existence of efficient estimators of  $\theta$. }
\end{exercise}

\begin{exercise}[Symmetric Location]\label{exer:6.3}
{\rm
Verify that Theorem \ref{thm:6.3} may be applied to the symmetric location
case.
To construct a  $\sqn$-consistent estimator of the location parameter, consider
 $\psi:\,{\mathbb R}\to{\mathbb R}$,  strictly increasing, uneven, bounded and with bounded
first and
second
derivatives.
Let  $\tt_n$  be the  $M$-estimator solving
\begin{equation}\label{6.20}
\nn\psi(X_i-\theta)=0.
\end{equation}
Note that
\begin{equation}\label{6.21}
P_{\theta,G}\Big(\sqn(\tt_n-\theta)\le
y\Big)=P_{\theta,G}\Big({1\over\sqn}\nn\psi(X_i-
\theta-{y\over\sqn})\le 0\Big)
\end{equation}
and show that  $\tt_n$  is  $\sqn$-consistent.
}
\end{exercise}

\begin{exercise}[Some Semiparametric Models]\label{exer:6.4}
{\rm We have to cope with the following estimation problem. $X_1,
\ldots, X_n$ are independent and identically distributed random
quantities with distribution $P_{\nu,g}$\,, where $\nu \in
{\mathbb R}$ and $g\in {\cal G}$ are unknown. We want to estimate
$\nu$ in the presence of the nuisance parameter $g$ and based on
$X_1, \ldots, X_n$. The following questions (A. through F.) have
to be answered for the estimation problems from the list following
the list of questions. The degree of mathematical detail in your
answers may vary and typically decreases when you proceed from
question A to F. For example, to answer question B you have to
determine tangent spaces on which to project. This is quite hard
to do exactly, but your efforts will be appreciated! Furthermore,
in answering question F it suffices to indicate the lines along
which an efficient estimator could be constructed: recall Section
8. Answering question F in full mathematical rigor would be
equivalent to writing a scientific paper about the model.

\begin{description}

\item
A. Fix $g_0$. If $g=g_0$ is known, then what does the Local
Asymptotic Spread Theorem tell us?

\item
B. Determine a lower bound on the asymptotic performance of
estimator sequences of $\nu$ in case $g$ is unknown.

\item
C. What is the information loss due to $g$ being unknown?

\item
D. Construct an efficient estimator of $\nu$ in case $g$ is known.

\item
E. Can you determine a $\sqrt n$-consistent estimator of $\nu$
when $g$ is unknown?

\item
F. (When possible!) Construct an efficient estimator of $\nu$ if
$g$ is unknown.

\end{description}

Formulate the regularity conditions that you need in your answers
and solutions; you are free to choose them yourself.

\begin{description}

\item
1. {\bf Symmetric Location}\\
\newblock Let $\cal G$ be the class of densities $g$ on ${\mathbb R}$ with respect
to Lebesgue measure $\lambda$\,, that are symmetric about 0, that
are absolutely continuous with Radon-Nikodym derivative $g'$, and
that have finite Fisher information for location
$$I(g) = \int \bigl({{g'} \over {g}}\bigr)^2 g\,. $$
$P_{\nu,g}$ has density $g(\cdot -\nu)$ on $({\mathbb R}, {\cal
B}, \lambda)$.

\item
2. {\bf Linear Regression}\\
\newblock The observation $X$ is defined by
$$ X=(Y,Z)\,,\, Y=\mu+\nu Z+\varepsilon.$$
Here, $\varepsilon$ and $Z$ are independent random variables,
$\varepsilon$ has unknown density $g_1$ on ${\mathbb R}$ with respect
to Lebesgue measure $\lambda$, that is absolutely continuous with
derivative $g'_1$ and that has finite Fisher information for
location
$$I(g_1) = \int \bigl({{g'_1} \over {g_1}}\bigr)^2 g_1\,, $$
$Z$ has unknown density $g_2$ on $({\mathbb R}, {\cal
B},\lambda)$, and $\mu \in {\mathbb R}$ is unknown as well. Of
course, $\cal G$ is the class of triplets $(\mu,g_1,g_2)$, where
$\mu, g_1$, and $g_2$ satisfy the above description.

\item
3. {\bf Heteroscedasticity}\\
\newblock The generic observation $X$ is defined by
$$ X=(Y,Z)\,,\, Y=\mu+e^{\nu Z}\varepsilon.$$
Here, $\varepsilon$ and $Z$ are independent random variables,
$\varepsilon$ has the standard normal distribution, $Z$ has finite
second moment and an unknown density $g$ on $({\mathbb R}, {\cal
B},\lambda)$, and $\mu \in {\mathbb R}$ is unknown as well.
Consequently, $\cal G$ is the class of pairs $(\mu,g)$.

\item
4. {\bf Partial Splines}\\
\newblock The generic observation $X$ is defined by
$$ X=(Y,Z)\,,\, Y=\nu+g(Z)+\varepsilon.$$
Here, $\varepsilon$ and $Z$ are independent, $\varepsilon$ is standard
normally distributed, and $Z$ is uniformly distributed on $(0,1)$.
Furthermore, $\cal G$ is the class of functions $g:(0,1)\to
{\mathbb R}$ with the properties
$$\int_0^1 g(z)dz=0, \int_0^1g^2(z)dz < \infty.$$
A model that comes closer to the one of Engle, Granger, Rice, and
Weiss (1986), is the following. $X=(W,Y,Z),
Y=\mu+\nu^TW+g(Z)+\varepsilon$, with $W,Z$, and $\varepsilon$
independent, $EW=0$, $\mu\in{\mathbb R}$ unknown, and $\cal G$ as
above. When possible, treat this problem as well.

\item
5. {\bf Projection Pursuit Regression}\\
\newblock The observation $X$ is defined by
$$ X=(Y,Z_1,Z_2)\,,\, Y=g(Z_1+\nu Z_2) +\varepsilon.$$
Here, $\varepsilon, Z_1,$ and $Z_2$ are independent standard normal
random variables. Furthermore, $g\,:\,{\mathbb R}\,\to\,{\mathbb
R}$ is an unknown differentiable function, and hence ${\cal
G}=\bigl\{g\,:\,{\mathbb R}\,\to\,{\mathbb R}\,\, {\rm
differentiable} \bigr\}$.

\item
6. {\bf Logistic Regression}\\
\newblock The generic observation $X$ is defined by
$ X=(W,Y,Z)$, where the random variable $Y$ has a Bernoulli
distribution with success probability
$$1 \over {1+e^{\nu Z+g(W)}}.$$
Here, $W$ and $Z$ are independent standard normal random variables
and ${\cal G}=\bigl\{g\,:\,{\mathbb R}\,\to\,{\mathbb R}\,\, {\rm
measurable}\bigr\}$.

\item
7. {\bf Errors in variables}\\
\newblock The observation $X$ is defined by
$$ X=(Y,Z)\,,\, Z=Z'+\varepsilon_1\,,\,Y=\nu Z' +\varepsilon_2.$$
Here, $Z',\varepsilon_1,$ and $\varepsilon_2$ are independent random
variables, $\varepsilon_1$ and $\varepsilon_2$ have a standard normal
distribution, and
$Z'$ has unknown density $g$ on $({\mathbb R}, {\cal B},\lambda)$.\\
{\em Hint}. $\dot{\cal P}_2$ exists of all random variables with
mean 0 and finite variance that can be written as a function of
$Z+\nu Y$.

\end{description}
}
\end{exercise}

\chapter{Time Series}\label{chap:9}

In the preceding Chapters we have studied asymptotically efficient estimation of
Euclidean parameters in semiparametric models for i.i.d. random variables.
Crucial in our set-up has been Local Asymptotic Normality of regular parametric
models and submodels as formulated in the LAN-Theorem \ref{thm:3.1}.
Via the LAS Theorems \ref{thm:new4.1} and \ref{thm:new6.2}
and the Convolution Theorems \ref{thm:4.1} and \ref{thm:5.2} LAN has provided a lower bound on the
asymptotic performance of estimators both in parametric and in nonparametric
models, and via the smoothness properties (\ref{3.9}) and (\ref{6.17}) it has been essential
in the proof of asymptotic linearity in the efficient influence  function of the
one-step estimators we have constructed; cf. the proof of Theorem \ref{thm:6.1}, the
sentence after (\ref{6.17}), and Theorem \ref{6.3}.

\bigskip

In many applications time-series models play an important part, e.g. in
econometrics.
They constitute a large group of models for non-i.i.d. observations,
 $Y_1,\dots,Y_n$.
Typically, an i.i.d. structure is hidden in a time-series  as follows.
Let $\varepsilon_1,\dots,\varepsilon_n$ be i.i.d. observations from a distribution with
density $g$, which are independent of the random vector $X_n$.
Assume that the $\si$-fields $\Ff^n_0=\Ff(X_n), \
\Ff^n_t=\Ff^n_{t-1}\vee\Ff(\varepsilon_t), \ t=1,\dots,n$, define a filtration.
One observes $X_n\in\Ff^n_0$ and $Y_1,\dots, Y_n$ with
\begin{equation}\label{7.1}
Y_t=\mu_t(\theta)+\si_t(\theta)\varepsilon_t,\quad t=1,\dots,n,
\end{equation}
where the time-dependent location-scale parameter $(\mu_t(\theta),
\ \si_t(\theta))^T\in(\Ff^n_{t-1})^2$ is supposed to depend on
$\theta$, the observed starting values $X_n\in\Ff^n_0$, and the
first $t-1$ observations $Y_1,\dots,Y_{t-1}$. So, $Y_t$ depends on
the ``past'' and the new independent innovation $\varepsilon_t$.

As in the i.i.d. case it is crucial to obtain LAN. Fix $g$ and
$\theta\in\Th\subset{\mathbb R}^k$. Let $\theta_n$ denote the true
parameter point, suppose $\theta_n\to\theta$, and let $\tt_n$ be
such that $\sqn(\tt_n-\theta_n)\to\la$. We denote the
log-likelihood ratio statistic of the observations $X_n,\,
Y_1,\dots,Y_n$ for $\tt_n$ with respect to $\theta_n$ by $\La_n$.
Let
\begin{equation}\label{7.2}
\Qq=\{Q_{\mu,\si}:\, \mu\in{\mathbb R}, \si>0\}
\end{equation}
with Lebesgue densities $q(\zeta)=\si^{-1}g(\si^{-1}(\cdot-\mu)), \
\zeta=(\mu,\si)^T$, of $Q_{\mu,\si}$, be the location-scale family of $g$
and denote
\begin{equation}\label{7.3}
\ell(\zeta)=\log q(\zeta).
\end{equation}
We assume that the innovations can be reconstructed  from the
observations if $\theta$ is known and we write
\begin{equation}\label{7.4}
\varepsilon_t(\theta)=\varepsilon(X_n,\, Y_1,\dots,Y_t,\,\theta).
\end{equation}
Under $\theta$ we have $\varepsilon_t(\theta)=\varepsilon_t, \ t=1,\dots,n$.
We define $k\times 2$-matrices $W_{nt}$ by their rows,
$j=1,\dots,k$,
\renewcommand{\arraystretch}{2}
\begin{equation}\label{7.5}
\quad W_{ntj}=
\left\{
\begin{array}{ll}
\si^{-1}_t(\theta_n)(\tt_{nj}-\theta_{nj})^{-1}
\Big(\mu_t(\tt_n^j)-\mu_t(\tt_n^{j-1}),\si_t(\tt_n^j)-\si_t(\tt_n^{j-1})\Big)
&\mbox{if}\quad \tt_{nj}\ne \theta_{nj}\\
\si^{-1}_t(\theta_n){\partial\over\partial\theta_j}(\mu_t(\theta),
\, \si_t(\theta)\Big)|_{\theta=\tt_n^j} &\mbox{if}\quad
\tt_{nj}=\theta_{nj},
\end{array}
\right.
\end{equation}
\renewcommand{\arraystretch}{1.5}
where ~$\tt_n^j=(\tt_{n1},\dots,\tt_{nj}, \
\theta_{n,j+1},\dots,\theta_{nk})^T$. Note that
\begin{equation}\label{7.6}
W^T_{nt}\Big(\tt_n-\theta_n\Big)=\si^{-1}_t(\theta_n)\Big(\mu_t(\tt_n)-\mu_t(\theta
_n),\, \si_t(\tt_n)-\si_t(\theta_n)\Big)^T.
\end{equation}
Now, the log-likelihood ratio statistic $\La_n$ may be written with
$\zeta_0=(0,1)^T$ as
\begin{equation}\label{7.7}
\La_n=\La^s_n+\sum
_{t=1}^n\Big\{\ell\Big(\zeta_0+W_{nt}^T(\tt_n-\theta_n)\Big)
\Big(\varepsilon_t(\theta_n)\Big)-\ell(\zeta_0)\Big(\varepsilon_t(\theta_n)\Big)\Big\}.
\end{equation}
Here $\La_n^s\in\Ff_0^n$ depends on the starting observations $X_n$
only, and we assume that they have a negligible effect, i.e.
\begin{equation}\label{7.8}
\La^s_n=o_{\theta_n,g}(1).
\end{equation}
Furthermore, we assume that $\Qq$ is a regular parametric model with
score function
\begin{equation}\label{7.9}
\del(x)=
\begin{pmatrix}
   -\,g'/g(x) \\
  -1-xg'/g(x)
\end{pmatrix},\quad x\in{\mathbb R},
\end{equation}
and Fisher information matrix
\begin{equation}\label{7.10}
J=\ex_g\del\del^T(X)
\end{equation}
at $\zeta_0=(0,1)^T$. Finally, we assume that there exists a
continuous nonsingular Fisher information matrix $I(\theta)$ and
square integrable $k\times 2$-matrices
$W_t(\theta_n)\in\Ff^n_{t-1}$ satisfying
\begin{equation}\label{7.11}
\ene\nn W_t(\theta_n)JW_t^T(\theta_n){\buildrel
P\over\longrightarrow}\, I(\theta),
\end{equation}
and
\begin{equation}\label{7.12}
\ene\sum_{t=1}^n\Big|W_t(\theta_n)\Big|^2{\been}_{[|W_t(\theta_n)|>\delta\sqn]}{
\buildrel P\over\longrightarrow}\, 0
\end{equation}
under $\theta_n$~ for all $\delta>0$, such that
\begin{equation}\label{7.13}
\sum_{t=1}^n\Big|(W_{nt}-W_t(\theta_n))^T(\tt_n-\theta_n)\Big|^2{\buildrel
P\over\longrightarrow}\, 0
\end{equation}
under $\theta_n$. We will denote the score functions of our
time-series model by
\begin{equation}\label{7.14}
\del_t(\theta)=W_t(\theta)\del(\varepsilon_t(\theta)), \ t=1,\dots, n.
\end{equation}

\begin{theorem} [LAN]\label{thm:7.1}
Under the above conditions, write
\begin{equation}\label{7.15}
\La_n=\la^T{1\over\sqn}\sum_{t=1}^n\del_t(\theta_n)-{1\over
2n}\sum_{t=1}^n \Big|\la^T\del_t(\theta_n)\Big|^2+R_n.
\end{equation}
Under $\theta_n$ and $g$
\begin{equation}\label{7.16}
R_n{\buildrel P\over\longrightarrow}\, 0, \
\La_n{\buildrel\D\over\longrightarrow}\,
\Nn\left(-\halfe\la^TI(\theta)\la, \ \la^TI(\theta)\la \right).
\end{equation}
The distributions of $(X_n,\, Y_1,\dots,Y_n)$ under $\theta_n$ and
$\tt_n$ are contiguous. In case $\theta_n=\theta+{\O}(n^{-1/2})$,
the smoothness condition holds:
\begin{equation}\label{7.17}
\sqn\Big(\tt_n-\theta_n+\ene\sum_{t=1}^n
I^{-1}(\tt_n)\del_t(\tt_n)-\ene\sum_{t=1}^n
I^{-1}(\theta_n)\del_t(\theta_n)\Big){\buildrel
P\over\longrightarrow}\, 0.
\end{equation}
\end{theorem}
{\bf Proof.}

A proof is given in Drost, Klaassen and Werker (1997). It is along
similar lines as the proof of Theorem \ref{thm:3.1}, to wit Taylor
expansion of (\ref{7.7}), but with ergodic theorems and martingale
difference central limit theorems replacing the law of large
numbers and the ``ordinary'' central limit theorem. \hfill $\Box$

\bigskip

The classical convolution theorem of H\'ajek (1970) shows that LAN
suffices for the convolution structure of the limit distribution
of regular estimators. (Note that in the proof of Theorem
\ref{thm:4.1} we only used (\ref{3.7}), (\ref{3.8}) and
contiguity, hence only LAN, and not anything else about regular
parametric models). Therefore, in our time-series models,
estimators $T_n$ of $q(\theta)$ are asymptotically efficient if
they satisfy
\begin{equation}\label{7.18}
\sqn\Big(T_n-q(\theta_n)-\ene\sum_{t=1}^n\dq(\theta_n)I^{-1}(\theta_n)\del_t(\theta_n)
\Big) {\buildrel P\over\longrightarrow}\, 0
\end{equation}
under $\theta_n$ and fixed $g$ for
$\theta_n=\theta+{\O}(n^{-1/2})$. Note that we have considered a
parametric model, in fact, since we kept $g$ fixed. Of course, the
convolution lower bound is still valid in the semiparametric model
with $g$ unknown. Often, under the ``right'' parametrization and
for the ``right'' function $q(\cdot)$ estimators  of $q(\theta)$
may be given for the semiparametric model that attain this bound
from the parametric submodel with $g$ fixed. We call such
estimators adaptive. The geometry of this situation is studied in
Drost, Klaassen and Werker (1994).

In Drost, Klaassen and Werker (1997), a general method for the construction
of such adaptive estimators is given based on sample splitting {\sl and}
discretization for the present time-series models.
We present two examples here.

\begin{example}[ARMA $(p,q)$]\label{exam:7.1}
{\rm
In the Auto Regressive Moving  Average model of orders $p$ and $q$, we
observe
realizations of $Y_1,\dots,Y_n$ with
\begin{equation}\label{7.19}
Y_t=\rho_1Y_{t-1}+\cdots+\rho_pY_{t-p}+\varf_1\varepsilon_{t-1}+
\cdots+\varf_q\varepsilon_{t-q}+\varepsilon_t.
\end{equation}

As in the general model (\ref{7.1}) the innovations $\varepsilon_t$ are
i.i.d. with unknown density $g$. Furthermore,
$\theta=(\rho_1,\dots,\rho_p, \, \varf_1,\dots,\varf_q)^T$.
Comparing (\ref{7.19}) to (\ref{7.1}) we note that
$\si_t(\theta)=1$ here. Consequently, we do not  need the full
location-scale model (\ref{7.2}), but only its restriction to
location. Therefore, we will assume that $g$ has finite Fisher
information for location $I_\ell(g)=\int(g'/g)^2g$, thus obtaining
regularity of the location model; cf. Exercises \ref{exer:3.1} and
\ref{exer:7.1}.

We will also assume $\expar_g\varepsilon_t=0, \ \expar_g\varepsilon_t^2<\infty$.
The starting observations are collected into the vector
$X_n=(Y_0,\dots,Y_{1-p}, \ \varepsilon_0,\dots,\varepsilon_{1-q})^T$ and we assume
(partially unrealistically) that $X_n$ will be observed too.
It may be checked straightforwardly that the innovations can be recovered from
the observations $X_n,\ Y_1,\dots,Y_n$ as in (\ref{7.4}).

Roughly speaking, an adaptive estimator of $\theta$ can be
constructed as follows. Via the first $p+q+1$ sample
autocovariances $\theta$ is estimated $\sqn$-consistently and this
estimator is discretized. Given this discretized estimator the
innovations can be recovered approximately from the  observations;
usually, one calls these ``estimated'' innovations the residuals.
Via these residuals the score function for location can be
estimated by the same methods as mentioned in Remark
\ref{rem:6.2}. Based on these estimators for $\theta$ and the
score function an efficient, adaptive estimator of $\theta$ can be
constructed based on sample splitting as in (\ref{6.16}). This
sample splitting and the discretization yield enough
``independence'' to be able to prove adaptivity (and hence
efficiency) of this estimator of $\theta$. For technical details
we refer to Example 4.2, pp. 809--811, of Drost, Klaassen and
Werker (1997). \hfill $\Box$ }
\end{example}

\begin{example}[GARCH $(p,q)$]\label{exam:7.2}
{\rm
Consider the time series model
\begin{equation}\label{7.20}
Y_t=\mu h^{1/2}_t+\si h^{1/2}_t\varepsilon_t
\end{equation}
with
\begin{equation}\label{7.21}
h_t=1+\be_1h_{t-1}+\cdots+\be_ph_{t-p}+\al_1Y^2_{t-1}+
\cdots+\al_qY^2_{t-q}.
\end{equation}
Originally, this model was introduced with $\al_1=\cdots=\al_q=0$; since
(\ref{7.21}) exhibits an autoregression structure for the $h_t$'s and since the
innovations are multiplied by a random factor, it was called an AutoRegressive
Conditional  Heteroskedasticity model: ARCH.
Therefore, the generalized model above is called GARCH.
Note that it fits into the framework of (\ref{7.1}) via
\begin{equation}\label{7.22}
\mu_t(\theta)=\mu h^{1/2}_t,  \si_t(\theta)=\si h^{1/2}_t,\quad
\theta=\Big(\al_1,\dots\al_q,  \be_1,\dots,\be_p, \mu,\si\Big)^T.
\end{equation}

Under appropriate regularity conditions and via a construction as
in Example \ref{exam:7.1}, an estimator of
$\nu(\theta)=(\al_1,\dots,\al_q, \ \be_1,\dots,\be_p)^T$ can be
constructed  which is efficient in the semiparametric model with
the density $g$ of the innovations $\varepsilon_t$ unknown. This
estimator is in fact adaptive in the presence of the nuisance
parameters $\mu$ and $\si$. This means that, with $\mu$ and $\si$
as nuisance parameters, it is performing asymptotically as well
under ~$g$~ unknown as the best estimator of $\nu(\theta)$ under
$g$ known (but still with $\mu$ and $\si$ unknown nuisance
parameters). Details may be found in Drost and Klaassen (1997).
\hfill $\Box$ }
\end{example}

\section{Exercises Chapter {\ref{chap:9}}}

\begin{exercise}[Location-Scale]\label{exer:7.1}
{\rm Let $g$ be a density on ${\mathbb R}$ which is absolutely
continuous with respect to Lebesgue measure with derivative $g'$.
Define the Fisher information for location and for scale
respectively, by
\begin{equation}\label{7.23}
I_\ell(g)=\int\Big({g'\over g}\Big)^2g, \ I_s(g)=\int
\Big(1+x{g'\over g}(x)\Big)^2g(x)dx.
\end{equation}
Use the technique of pp. 211--212 and p. 214 of H\'ajek and
$\breve {\rm S}$id\'ak (1967) to show that finiteness of
$I_\ell(g)$ and $I_s(g)$ imply that the location-scale family
$\Qq$ from (\ref{7.2}) is a regular parametric model. \hfill
$\Box$ }
\end{exercise}

\chapter{Banach Parameters}\label{chap:10}

Up to now, we have considered only Euclidean (finite dimensional)
parameters, i.e. parameterfunctions ~$\nu:\, \Pp\to{\mathbb R}^m$.
Here we will discuss generalization to Banach (infinite
dimensional) parameters, i.e. to parameterfunctions  $\nu:\,
\Pp\to\B$ with $\B$ a Banach space with norm $\|\cdot\|_\B$. The
connection between Euclidean and Banach parameters is made via the
dual space $\B^*$. Recall that $\B^*$ is the Banach space of real
valued bounded linear functions $b^*:\B\to{\mathbb R}$. So, for
any $b^*\in\B^*$~ the parameterfunction $b^*\nu:\,\Pp\to{\mathbb
R}$ is Euclidean. Not surprisingly the convolution theorem for
estimators of Euclidean parameters thus gives rise to convolution
theorems for Banach parameters. To formulate them accurately we
need some definitions similar to those in Chapter \ref{chap:6}.

\begin{definition}[Pathwise
(Weak-)Differentiability]\label{defn:8.1}
{\rm
The parameter $\nu:\,\Pp\to\B$ is {\sl pathwise differentiable} on $\Pp$ at
$P_0\in\Pp$, if there exists a bounded linear operator
$\dun(P_0)=\dun:\,\dP\to\B$ such that for all one-dimensional regular
parametric
 submodels $\{P_\eta:|\eta|<\varepsilon\}, \ \varepsilon>0$, with score
function $h$ at $P_0$
\begin{equation}\label{8.1}
\|\nu(P_\eta)-\nu(P_0)-\eta\dun(h)\|_\B=o(\eta), \ {\rm as}~ \eta\to 0.
\end{equation}
The parameter $\nu$ is {\sl pathwise weak-differentiable} if
\begin{equation}\label{8.2}
b^*\nu(P_\eta)=b^*\nu(P_0)+\eta b^*\dun(h)+o(\eta),
\end{equation}
for all ~$b^*\in\B^*$.
}
\end{definition}

The transpose  of $\dun:\,\dP\to\B$ is the map $\dun^T:\,\B^*\to\dP$
defined by
\begin{equation}\label{8.3}
\langle\dun^Tb^*,\,h\rangle_0=b^*\dun(h), \quad h\in\dP.
\end{equation}
By the Riesz representation theorem (cf. (\ref{5.2})) there exists a unique
$\dun_{b^*}\in\dP$ with
\begin{equation}\label{8.4}
b^*\dun(h)=\langle\dun_{b^*},\, h\rangle_0,
\end{equation}
and hence we have
\begin{equation}\label{8.5}
\dun_{b^*}=\dun^Tb^*, \quad b^*\in\B^*.
\end{equation}

\begin{definition}[Efficient Influence Operator and
Function]\label{defn:8.2}
{\rm
If $\nu$ is pathwise (weak-) differentiable with derivative operator
$\dun:\,\dP\to\B$, then the {\sl efficient influence operator}
$\tel(P_0\mid\nu,\Pp)=\tel_\nu:\,\B^*\to\dP$ is defined as its transpose
\begin{equation}\label{8.6}
\tel_\nu(b^*)=\dun^Tb^*=\dun_{b*},
\end{equation}
and the {\sl inverse information covariance functional}
$I^{-1}(P_0\mid\nu,\Pp)=I_\nu^{-1}:\B^*\times\B^*\to{\mathbb R}$
by
\begin{equation}\label{8.7}
I^{-1}_\nu(b^*_1,b^*_2)=\langle\tel_\nu(b^*_1),\,\tel_\nu(b^*_2)\rangle_0.
\end{equation}
If there exists a map $\tel:\Xx\to\B$ such that for all $b^*\in\B^*$
\begin{equation}\label{8.8}
b^*\tel(\cdot)=\dun_{b^*}(\cdot)=\tel_\nu(b^*)(\cdot),
\end{equation}
then we call $\tel$ the {\sl efficient influence function}.
}
\end{definition}

\begin{definition}[Weak Regularity]\label{defn:8.3}
{\rm The estimator sequence $\{T_n\}$ of $\nu(P)$ is said to be
{\sl weakly regular} at $P_0\in\Pp$ if there exists a process
$\{b^*\Z:\, b^*\in\B^*\}$ on $({\mathbb R}^{\B^*},\Bb^{\B^*})$
such that for all one-dimensional regular parametric submodels
$\{P_\eta:\, |\eta|<\varepsilon\}, \ \varepsilon>0$, sequences
$\eta_n={\O}(n^{-1/2})$, and $b^*\in\B^*$
\begin{equation}\label{8.9}
\sqn\Big(b^*T_n-b^*\nu(P_{\eta_n})\Big)\to b^*\Z, ~{\rm as}~ n\to\infty.
\end{equation}
}
\end{definition}

\bigskip

With these definitions we may formulate a generalization of Theorem \ref{thm:5.2}.

\begin{theorem}[Convolution Theorem]\label{thm:8.1}
Let $\nu:\,\Pp\to\B$ be pathwise weak-differentiable at
$P_0\in\Pp$, let $\{T_n\}$ be weakly regular with limit $\Z$. If
the tangent  set $\dP^0$ is linear, then there exist processes
$\{b^*\Z_0:b^*\in\B^*\}$ and $\{b^*\Del_0:b^*\in\B^*\}$ on
$\{{\mathbb R}^{\B^*},\Bb^{\B^*}\}$ such that:
\begin{equation}\label{8.10}
b^*\Z\,{\mathop=^{\rm D}}\, b^*\Z_0+b^*\Del_0, \quad b^*\in\B^*,
\end{equation}
$\Z_0$ and $\Del_0$ are independent, $\Z_0$  is Gaussian with mean
$0$ and
\begin{equation}\label{8.11}
\Cov(b^*_1\Z_0,\, b^*_2\Z_0)=I^{-1}_\nu(b^*_1,b^*_2),
\end{equation}
and for every $h\in\dP$ and $b^*\in\B^*$
\begin{equation}\label{8.12}
{\sqn\Big(b^*T_n-b^*\nu(P_0)-\ene\nn\tel_\nu(b^*)(X_i)\Big)\choose{1\over
\sqn} \nn
h(X_i)}{\buildrel D\over \longrightarrow_{P_0}{b^*\Del_0\choose
W_0}}
\end{equation}
with $\Del_0$ and $W_0$ independent. Furthermore, $\Del_0$ is
degenerate at $0\in{\mathbb R}^{\B^*}$ iff for every $b^*\in\B^*,
\ b^*T_n$ is asymptotically linear with influence function
$\tel_\nu(b^*)$, and $\{T_n\}$ is called weakly efficient at $P_0$
in this case.
\end{theorem}

This theorem is proved in Section 5.2 of BKRW (1993) together with a further
generalization for regular estimators.

Consider a semiparametric model in which we have an efficient estimator of the Euclidean parameter $\theta$, and for which,
given the value of $\theta$, there exists an estimator of the Banach parameter, which depends on $\theta$ and
is efficient within this restricted model (with $\theta$ known).

Substituting the efficient estimator of $\theta$
for the value of this parameter in the estimator of the Banach
parameter, one can obtain by a sample splitting technique as in Theorem \ref{thm:6.2}, an efficient estimator of the Banach
parameter for the full semiparametric model with the Euclidean
parameter unknown; see Klaassen and Putter (2005).
This heredity property of efficiency completes
estimation in semiparametric models in which the Euclidean
parameter has been estimated efficiently. Typically, estimation
of both the Euclidean and the Banach parameter is necessary in
order to describe the random phenomenon under study to a sufficient
extent.

\begin{example}[Estimation of Distribution Functions]\label{exam:8.1}
{\rm Let $\Pp$ be the collection of all
probability measures on ${\mathbb R}$ and $\nu(P)$ the
distribution function $F$ of $P$. With $\B$ the space of cadlag
functions with the supnorm, $\nu:\,\Pp\to\B$ is pathwise
differentiable at $P_0$ with $\dun:\,\dP\to\B$ given by
\begin{equation}\label{8.13}
\dun(h)(t)=\int\Big({\been}_{[x\le t]}-F_0(t)\Big)h(x)dP_0(x),
\quad t\in{\mathbb R}.
\end{equation}
If $b^*:\,\Bb\to{\mathbb R}$ is such that $b^*(f)=f(t)$, then
\begin{equation}\label{8.14}
\tel_\nu(b^*)(x)={\been}_{[x\le t]}-F_0(t),
\end{equation}
and from the convolution theorem we conclude that the empirical
distribution function
\begin{equation}\label{8.15}
{\hat F}_n(t)=\ene\nn{\been}_{[X_i\le t]}, \quad t\in{\mathbb R},
\end{equation}
is a weakly efficient estimator of $\nu(P)=F(\cdot),$ with
efficient influence function $\tel:\Xx\to\B$ satifying $\tel(x): t
\mapsto {\been}_{[x\le t]}-F_0(t).$
\hfill $\Box$ }
\end{example}

Next to its asymptotic efficiency the empirical distribution function is efficient for finite samples too, as we show by proving
a Cram\'er-Rao inequality for unbiased estimators of distribution functions.
\begin{theorem}[Cram\'er-Rao Inequality for Distribution Functions]\label{thm:10.2}
Let $\Pp$ be the collection of all probability measures on
${\mathbb R}$ and $\nu(P)=F$ the distribution function of $P$.
For every unbiased estimator ${\bar F}_n$ of
$\nu(P)=F$ and every $P \in \Pp$ with distribution
function $F$ we have
\begin{equation}\label{8.16}
\var_P \left( {\bar F}_n(t) \right) \geq \var_P \left( {\hat F}_n(t) \right)
= \frac 1n F(t)\left(1-F(t)\right), \quad t \in {\mathbb R},
\end{equation}
with ${\hat F}_n$ as in (\ref{8.15}).
\end{theorem}
{\bf Proof.}
Fix  $t \in {\mathbb R}$ and the distribution function $F$.
For $0\leq \varepsilon \leq 1$ we define
\begin{equation}\label{8.17}
F_\varepsilon(x) = F(x) - \varepsilon \left(F(x \wedge t) - F(x) F(t) \right), \quad x \in {\mathbb R}.
\end{equation}
Note that $F_\varepsilon$ is a distribution function with $F_0=F$ and
\begin{equation}\label{8.18}
dF_\varepsilon(x) = \left( 1 - \varepsilon \left( {\bf 1}_{[x \leq t]} - F(t) \right) \right)\, dF(x).
\end{equation}
Unbiasedness of an estimator ${\bar F}_n(t) = {\bar F}_n(t; X_1,\dots,X_n)$ means
\begin{eqnarray}\label{8.19}
\lefteqn{ F_\varepsilon (t) = \expar_{F_\varepsilon} {\bar F}_n(t) } \\
&& = \idotsint_{{\mathbb R}^n} {\bar F}_n(t; x_1,\dots,x_n) \prod_{j=1}^n \left( 1 - \varepsilon \left( {\bf 1}_{[x_j \leq t]} - F(t) \right) \right)\, dF(x_j). \nonumber
\end{eqnarray}
Taking the right hand derivative of (\ref{8.19}) with respect to $\varepsilon$ at 0, we obtain
\begin{eqnarray}\label{8.20}
\lefteqn{ F(t) (1-F(t)) = \idotsint_{{\mathbb R}^n} {\bar F}_n(t; x_1,\dots,x_n) \sum_{i=1}^n \left( {\bf 1}_{[x_i \leq t]} - F(t) \right)\, \prod_{j=1}^n dF(x_j) } \\
&& \hspace{10em} = n \ \expar_F \left( {\bar F}_n(t) \left( {\hat F}_n(t) -F(t) \right) \right), \nonumber
\end{eqnarray}
which implies
\begin{equation}\label{8.21}
\Cov \left( {\bar F}_n(t),{\hat F}_n(t) \right) = \var \left( {\hat F}_n(t) \right).
\end{equation}
Applying the Cauchy-Schwarz inequality to (\ref{8.21}), we obtain (\ref{8.16}).
\hfill $\Box$

\section{Exercises Chapter {\ref{chap:10}}}

\begin{exercise}[Interchange of differentiation and integration]\label{exer:10.1}
{\rm Prove the validity of (\ref{8.20}). {\em Hint}: Write the difference quotient of the integrand in the multiple integral from (\ref{8.19})
as a telescoping sum and apply dominated convergence.
\hfill $\Box$ }
\end{exercise}

\appendix

\chapter{Additions}

\section{Measure Theory and Probability}

\subsection{Vitali's theorem}
Pointwise convergence together with convergence of the norms
implies convergence {\em in} norm; more precisely:
\begin{theorem}[Vitali]\label{vitali}
If $f_n\to f$, $\mu$-almost everywhere, and
$$\limsup_{n\to\infty} \int |f_n|^pd\mu \leq \int |f|^pd\mu,\quad0<p<\infty,$$
then we have
\begin{equation}\label{vitaly:1}
\int|f_n-f|^pd\mu \to 0.
\end{equation}
\end{theorem}
\medskip
\noindent{\bf Proof} Note that we have $(a+b)^p\leq 2^p(a^p+b^p)$ for all
$a,b \geq 0$. So
\begin{equation}\label{pos}
|f_n-f|^p \leq (|f_n|+|f|)^p\leq 2^p(|f_n|^p+|f|^p).
\end{equation}
The  convergence of $f_n$ to $f$  implies
$$2^{p+1}|f|^p=\lim_{n\to\infty}\Big(2^p(|f_n|^p+|f|^p) - |f_n-f|^p\Big),$$
$\mu$-almost everywhere, and by Fatou's lemma
\begin{eqnarray}
\lefteqn{\int 2^{p+1}|f|^pd\mu}\nonumber\\
&=&\int \lim_{n\to\infty}\Big(2^p(|f_n|^p+|f|^p) - |f_n-f|^p\Big)d\mu\nonumber\\
&\leq&\liminf_{n\to\infty}\int\Big(2^p(|f_n|^p+|f|^p) - |f_n-f|^p\Big)d\mu\label{newA.3}\\
&\leq&\limsup_{n\to\infty}2^p\int|f_n|^pd\mu + 2^p\int|f|^pd\mu -
\limsup_{n\to\infty}\int|f_n-f|^pd\mu\nonumber\\
&\leq&\int 2^{p+1}|f|^pd\mu -
\limsup_{n\to\infty}\int|f_n-f|^pd\mu.\nonumber
\end{eqnarray}
This implies (\ref{vitaly:1}). Note that by (\ref{pos}) the second
integrand in (\ref{newA.3}) is positive which is essential in an
application of Fatou's lemma. \hfill$\Box$

\subsection{Total variation distance}\label{tv}

Consider a set of measures ${\cal P}$ on a separable metric space ${\cal X}$
with
Borel sets ${\cal B}$. We define the {\em total variation
distance} between two measures $P$ and $Q$ as
\begin{equation}\label{totvar}
d(P,Q) = \sup_{A\in {\cal B}} |P(A) - Q(A)|.
\end{equation}
If both $P$ and $Q$ are dominated by $\mu$ then $d(P,Q)$ is half the
${\cal L}_1(\mu)$ distance between the densities
\begin{equation}\label{l1}
d(P,Q) = \int_{A_0}|p-q|d\mu = {1\over 2}\int_{{\cal X}}|p-q|d\mu,
\end{equation}
where $A_0=\{x\in {\cal X}:p(x)>q(x)\}$.

Examples of sets $K$ of measures that are compact in this metric,
are for instance sets existing of finitely many elements, or sets
$\{P,P_1,P_2,\dots\}$ where $(P_n)_{n=1}^\infty$ is a sequence of
measures converging to $P$ in total variation. 
If ${\cal P}=\{P_\theta:\theta\in \Theta\}$ is a regular parametric model
then for $\theta_n\to\theta$ and $t_n\to t$ the sequence
$(P_{\theta_n+ t_n/\sqrt{n}})_{n=1}^\infty$ converges in total
variation to $P_\theta$. This follows from the fact that the
Fr\'echet differentiability implies $\L_2$ convergence of the
densities, and hence $\L_1$ convergence too. Hence
$K=\{P_{\theta_n+t_n/\sqrt{n}},n=1,2,\dots\}$ is also a
compact set.

\subsection{Tightness and Prohorov's theorem}\label{proh}

We need the concept of {\em tightness} of a sequence of distributions.

\begin{definition}[Tigthness]
{\em A sequence of probability measures $(P_n)_{n=1}^\infty$ on a
measurable space $({\cal X},{\cal B})$ is called tight if for every
$\varepsilon>0$ there exists a compact set $K_\varepsilon$ such that
$P_n(K_\varepsilon)>1-\varepsilon$, for all $n$.
}
\end{definition}

The next theorem establishes a correspondence between tightness and weak
convergence of subsequences. For a proof we refer to Billingsley (1968)
page 37.

\begin{theorem}[Prohorov]
Let  $(P_n)_{n=1}^\infty$ denote a sequence of probability measures
on a complete separable measurable space $({\cal X},{\cal B})$.
Then $(P_n)_{n=1}^\infty$ is tight if and only if for each subsequence
there exists a further subsequence that converges weakly to a probability
measure $P$. This limit measure $P$ may depend on the specific subsequences.
\end{theorem}

Weak convergence of the whole sequence can also be characterized in a
similar manner. The proof is straightforward.

\begin{theorem}
Let  $(P_n)_{n=1}^\infty$ denote a sequence of probability measures
on a  measurable space $({\cal X},{\cal B})$.
Then $(P_n)_{n=1}^\infty$ converges weakly to a probability measure $P$
if and only if for each subsequence
there exists a further subsequence that converges weakly to  $P$.
\end{theorem}

\chapter{Notes}

\section{Notes Chapter \ref{chap:2}}\label{noteschap2}

\begin{Note}\label{note:2.1}

\noindent {\bf Proof of the "if" part of Theorem \ref{thm:2.2}}.

Under (\ref{new2.16}) equality in (\ref{2.10}) can be violated only if
$G^{-1}(v)-G^{-1}(u) >\int_u^v 1/g(G^{-1}(s))ds$.
This can happen only if $G^{-1}$ has jumps, that is, if $g$ vanishes on
$[y_0,y_1]$, say, with $y_0<y_1$ and with $g$ positive at some points
$y<y_0$ and some $y>y_1$. However, by (\ref{2.4})
\begin{eqnarray*}
\lefteqn{g(y)=\ex SI_{[T-\vart>y]}}\\
&=&
\ex SI_{[G^{-1}(H(S))>y]}\\
&=&
\ex (S|G^{-1}(H(S))>y)P(G^{-1}(H(S))>y)
\end{eqnarray*}
holds and the first factor at the right hand side is nonnegative and nondecreasing in
$y$.
Now, assume $g(y_2)=0$ for some $y_2$. If $P(G^{-1}(H(S))>y_2)=0$
then $P(G^{-1}(H(S))>y)=0$ and hence $g(y)$
vanishes for all $y\geq y_2$.
If $\ex(S| G^{-1}(H(S))>y_2)=0$, then $\ex(S| G^{-1}(H(S))>y)$
and hence $g(y)$ vanish for all $y\leq y_2$.
It follows that $g$ cannot vanish on an
interval strictly within its support. Consequently, $G^{-1}$
cannot have jumps.
\hfill$\Box$
\end{Note}

\newpage

\section{Notes Chapter \ref{chap:3}}\label{noteschap3}

\begin{Note}\label{note:1}

Recall that the directional derivative of a function $g: {\mathbb
R}^k\to{\mathbb R}$ on the line segment from $a\in {\mathbb R}^k$
to $b\in {\mathbb R}^k$ is given by
$$
{d\over d\lambda}g(a+\lambda(b-a)) =
(b-a)^Tg'(a+\lambda(b-a)),\quad \lambda\in [0,1],
$$
where $g'$ is the gradient vector of partial derivatives of $g$.

\noindent In the context of the proof of Proposition \ref{prop:3.1}
we get, since ${\dot s}={1\over 2}{\del s}$,
\begin{eqnarray*}
\lefteqn{ s(x;\theta)-s(x;\theta_0)}\\
&=& \int_0^1 {d\over d\lambda} s(\theta_0+\lambda (\theta-\theta_0))\, d\lambda\\
&=& \int_0^1 (\theta-\theta_0)^T{\dot s}(\theta_0+\lambda (\theta-\theta_0))\, d\lambda\\
&=& \int_0^1 {1\over 2}(\theta-\theta_0)^T\del s(\theta_0+\lambda (\theta-\theta_0))\, d\lambda.
\end{eqnarray*}
\end{Note}

\begin{Note}\label{note:2}

Writing
$$g(\lambda)={1\over 2}(\theta-\theta_0)^T\del s(\theta_0+\lambda
(\theta-\theta_0))
$$
use the inequality (Jensen or Cauchy-Schwarz)
$$
\Big(\int_0^1 g(\lambda)d\lambda\Big)^2 \leq \int_0^1 g(\lambda)^2d\lambda .
$$
\end{Note}

\begin{Note}\label{note:3}

Let $a$ and $b$ be $k$-vectors. The following equality will be
used repeatedly
\begin{equation}\label{neq:0}
(a^Tb)^2= a^T b(a^Tb) = a^T b (a^T b)^T = a^T b(b^T a) = a^T(bb^T)a.
\end{equation}
By this equality we obtain
\begin{eqnarray*}
\lefteqn{\int (t^T\del s(x;\theta))^2d\mu(x) =
\int (t^T\del(x;\theta))^2 s(x;\theta)^2d\mu(x)}\\
&=&\int (t^T\del\del^T(x;\theta)t)p(x;\theta)d\mu(x)=
t^T\int \del\del^T(x;\theta)p(x;\theta)d\mu(x)\ t\\
&=&t^TI(\theta)t.
\end{eqnarray*}
Now substitute $\theta-\theta_0$ for $t$ and
$\theta_0+\lambda(\theta-\theta_0)$ for $\theta$.
\end{Note}

\begin{Note}\label{note:4}

Suppose that (\ref{d}) is not true. Then there exists a sequence
$\{\theta_n\}$ such that
\begin{equation}\label{ong:1}
{1\over{|\theta_n-\theta_0|^2}} \int\limits_{p(x;\theta_0)>0}
|s(x;\theta_n)-s(x;\theta_0)-\halfe(\theta_n- \theta_0)^T\del
s(x;\theta_0)|^2d\mu(x)
\end{equation}
is at least at distance $\varepsilon$ from $0$ for some $\varepsilon>0$.
Now consider the sequence
$\{(\theta_n-\theta_0)/|\theta_n-\theta_0|\}$. This is a sequence
of points on the unit ball of ${\mathbb R}^k$. So there is a
convergent subsequence for which (\ref{ong:1}) does not converge.
Hence it suffices to prove (\ref{d}) for all sequences
$\{\theta_n\}$ for which $(\theta_n-\theta_0)/|\theta_n-\theta_0|$
converges.

Now consider an arbitrary sequence  $\{\theta_n\}$ with
$(\theta_n-\theta_0)/|\theta_n-\theta_0| \to c \in {\mathbb R}^k$.
Then (\ref{ong:1}) equals
\begin{equation}\label{ong:2}
\int\limits_{p(x;\theta_0)>0} |{s(x;\theta_n)-s(x;\theta_0)
\over{|\theta_n-\theta_0|}}-\halfe c^T  \del
s(x;\theta_0)|^2d\mu(x)+o(1)
\end{equation}
 From (\ref{b}) we get
$$
{s(x;\theta_n)-s(x;\theta_0) \over{|\theta_n-\theta_0|}} \to
\halfe c \del s(x;\theta_0)
$$
and by  (\ref{c}) we have (use $c^T\del\del^Tc=(c^T\del)^2$)
$$
\limsup_{n\to \infty} \int\limits_{p(x;\theta_0)>0}
 \Big|{s(x;\theta_n)-s(x;\theta_0) \over{|\theta_n-\theta_0|}}
\Big|^2 d\mu(x) \leq \int\limits_{p(x;\theta_0)>0} \Big| \halfe
c^T \del s(x;\theta_0)\Big|^2d\mu(x).
$$
By Vitali's theorem \ref{vitali} we can now conclude that (\ref{ong:1})
converges to zero. This proves (\ref{d}).
\end{Note}

\begin{Note}\label{note:5}

\begin{lemma}\label{expscore}
If ${\cal P} = \{P_\theta :\theta\in\Theta\}$ is a regular
parametric model then
\begin{equation}
\int \del(x;\theta)dP_\theta(x) = 0.
\end{equation}
\end{lemma}

\medskip
\noindent{\bf Proof}
By the Fr\'echet differentiability we have
$$\|s(\theta+h) - s(\theta)\|\leq \|s(\theta+h) - s(\theta)-{1\over2}h^T \del(\theta) s(\theta)\|+ \|{1\over2}h^T \del s\| ={\cal O}(|h|).$$
Now note that, working in ${\cal L}_2(\mu)$, we have
$$<s(\theta),s(\theta)> = \int s(\theta)s(\theta)d\mu = \int p(\theta)d\mu= 1.$$
This implies
\begin{eqnarray*}
\lefteqn{|<h^T\del(\theta)s(\theta),s(\theta)>|}\\
&& = |<s(\theta+h),s(\theta+h)> - <s(\theta),s(\theta)> - <{1\over 2}h^T\del(\theta)s(\theta), 2 s(\theta)>|\\
&& = |<s(\theta+h)-s(\theta)-{1\over 2}h^T\del(\theta)s(\theta),2s(\theta)> +<s(\theta+h)-s(\theta),s(\theta+h)-s(\theta)>|\\
&& \leq |<s(\theta+h)-s(\theta)-{1\over2}h^T\del(\theta)s(\theta),2s(\theta)>| + |<s(\theta+h)-s(\theta),s(\theta+h)-s(\theta)>|\\
&& \leq 2\|s(\theta+h)-s(\theta)-{1\over 2}h^T\del(\theta)s(\theta)\| + \|s(\theta+h)-s(\theta)\|^2\\
&& = {\cal O}(|h|),\quad\mbox{as}\ h\to 0.
\end{eqnarray*}
But, by taking $h=(0,\dots,0,h_i,0,\dots,0)^T$ with $h_i\to 0$,
this implies
$$\int \del_i(x;\theta)dP_\theta(x) = <\del(\theta)s(\theta),s(\theta)>=0, \quad i=1,\dots,k.$$
\hfill$\Box$
\end{Note}

\newpage
\begin{Note}\label{note:6}

\begin{lemma}
If the parametrization  $\theta\longrightarrow P_{\theta}$ is
regular and $T_{n}$ is defined by
\[
T_{n}=T_{n}(X,t)= 2 \left( {{s(X;\theta
+{{t}\over{\sqrt{n}}})-s(X;\theta)}\over{s(X;\theta)}} \right),
\]
then we have uniformly for $\theta\in K$ and $|t|\leq M$
\begin{enumerate}
\item
${\rm E}_{\theta}T_{n}= -{1\over{4n}}
t^TI(\theta)t+o({{1}\over{n}})$
\item
${\rm E}_{\theta}T_{n}^{2}= {1\over{n}}
t^TI(\theta)t+o({{1}\over{n}})$
\item
${\rm E}_{\theta} \left(
T_{n}-{{t^T}\over{\sqrt{n}}}\dot{\ell}(\theta) \right)^{2}
=o({{1}\over{n}})$
\item
$P_{\theta} \left( \left\{ |T_{n}|>\varepsilon \right\} \right)
=o({{1}\over{n}}) \;\;\;{\rm voor\;alle\;}\varepsilon >0$.
\end{enumerate}
\end{lemma}

\noindent{\bf Proof}

We give the proof for $\theta$ and $t$ fixed. The generalization
to convergent sequences $\theta_n$ and $t_n$, to prove uniformity,
is straightforward.

\noindent Write $A(\theta)=\{x:\;s(x;\theta)=0\}$. We need the
following bound
\begin{equation}\label{neq:1}
\int_{A(\theta)}s^{2}(x;\theta+h)d\mu(x)=o(|h|^{2}),\quad\mbox{if}\
h\rightarrow 0.
\end{equation}
This follows from the regularity of the parametrization, since we have
\begin{eqnarray*}
\lefteqn{
\int_{A(\theta)}s^{2}(x;\theta+h)d\mu(x)}\\
&\leq&
\int_{A(\theta)}(s(x;\theta+h)+s(x;\theta-h))^{2}d\mu(x)\\
&=& \int_{A(\theta)}(s(x;\theta+h)-{1\over 2}h^T\del s(x;\theta)+
{1\over 2}h^T\del s(x;\theta)
+s(x;\theta-h))^{2}d\mu(x)\\
&\leq&
2\int_{A(\theta)}(s(x;\theta+h)-{1\over 2}h^T\del s(x;\theta))^{2}d\mu(x)\\
&&+\quad2\int_{A(\theta)}(s(x;\theta-h)+{1\over 2}h^T\del
s(x;\theta))^{2}
d\mu(x)\\
&=&
o(|h|^{2})
\quad\mbox{as}\ h\rightarrow 0.
\end{eqnarray*}
Here we have used the inequality $(a+b)^2\leq 2(a^2+b^2)$. \\
Next we prove ${\rm E}_{\theta}T_{n}^{2}=-4{\rm
E}_{\theta}T_{n}+o({1/n})$. Define the set
$B(\theta)=A(\theta)^c$$=\{x:\;s(x;\theta)>0\}$, then, using $\int
s(x;\theta)^2d\mu = \int s(x;\theta+{t\over \sqrt{n}})^2d\mu = 1$,
we get
\begin{eqnarray*}
\lefteqn{{\rm E}_{\theta}T_{n}^{2} = 4 \int_{\bf X} \left( {{s(x;\theta +{{t}\over{\sqrt{n}}})-s(x;\theta)}\over{s(x;\theta)}} \right)^{2}
dP_{\theta}(x)
= 4 \int_{B(\theta)} \left( {{s(x;\theta+{{t}\over{\sqrt{n}}})-s(x;\theta)}\over{s(x;\theta)}} \right)^{2}
dP_{\theta}(x) } \\
&&= 4 \int_{B(\theta)} \left( s(x;\theta+{{t}\over{\sqrt{n}}})-s(x;\theta) \right)^{2}d\mu(x)\\
&&= 4 \int_{B(\theta)} \Big(s^{2}(x;\theta +{{t}\over{\sqrt{n}}})-2s(x;\theta +{{t}\over{\sqrt{n}}})s(x;\theta)+s^{2}(x;\theta)\Big)d\mu(x)\\
&&= -8 \int_{B(\theta)} s(x;\theta) \left( s(x;\theta+{{t}\over{\sqrt{n}}})-s(x;\theta) \right)d\mu(x)- 4 \int_{A(\theta)} s^{2}(x;\theta +{{t}\over{\sqrt{n}}})d\mu(x)\\
&&\stackrel{(\ref{neq:1})}{=} -8 \int_{B(\theta)} \left({{s(x;\theta+{{t}\over{\sqrt{n}}})-s(x;\theta)}\over{s(x;\theta)}} \right)
dP_{\theta}(x) + o({{1}\over{n}})\\
&& = -4{\rm E}_{\theta}T_{n} +o({{1}\over{n}}).
\end{eqnarray*}
Next we prove { 2.}, which in its turn proves { 1.}
We have
\begin{eqnarray*}
\lefteqn{{\rm E}_{\theta}T_{n}^{2} = 4 \int_{B(\theta)} \left(
s(x;\theta +{{t}\over{\sqrt{n}}})-s(x;\theta) \right)^{2}
d\mu(x)}\\
&=& 4 \int_{B(\theta)} \left( s(x;\theta
+{{t}\over{\sqrt{n}}})-s(x;\theta)-{1\over
2}{{t^T}\over{\sqrt{n}}} \del s(x;\theta)+{1\over
2}{{t^T}\over{\sqrt{n}}} \del s(x;\theta) \right)^{2}
d\mu(x)\\
&=& 4 \int_{B(\theta)} \left( s(x;\theta
+{{t}\over{\sqrt{n}}})-s(x;\theta)-{1\over
2}{{t^T}\over{\sqrt{n}}} \del s(x;\theta) \right)^{2}
d\mu(x)\\
&&\quad + {1\over n}\int_{B(\theta)} (t^T\del s (x;\theta))^2
d\mu(x)\\
&&\quad+ \int_{B(\theta)} {{4t^T}\over{\sqrt{n}}}\del s(x;\theta)
\left( s(x;\theta +{{t}\over{\sqrt{n}}})-s(x;\theta)-{1\over
2}{{t}\over{\sqrt{n}}} \del s(x;\theta) \right)
d\mu(x)\\
&\stackrel{(\ref{neq:0})}{=}& {1\over{n}} t^TI(\theta)t
+o({{1}\over{n}}),
\end{eqnarray*}
since by the Cauchy-Schwarz inequality
\begin{eqnarray*}
\lefteqn{ \left( \int_{B(\theta)} t^T\del s(x;\theta) \left(
s(x;\theta +{{t}\over{\sqrt{n}}})-s(x;\theta)-{1\over
2}{{t^T}\over{\sqrt{n}}} \del s(x;\theta) \right) d\mu(x)
\right)^{2}}\\
&\leq& \left( \int_{B(\theta)} \left( s(x;\theta
+{{t}\over{\sqrt{n}}})-s(x;\theta)-{1\over
2}{{t^T}\over{\sqrt{n}}} \del s(x;\theta) \right)^{2} d\mu(x)
\right)\\
&&\quad\quad\times \left( \int_{B(\theta)} (t^T\del s(x;\theta))^2
d\mu(x)
\right)\\
&=&
o({{1}\over{n}}).
\end{eqnarray*}
Statement {3.} follows from
\begin{eqnarray*}
\lefteqn{{\rm E}_{\theta} \left(
T_{n}-{{t^T}\over{\sqrt{n}}}\del(\theta)
\right)^{2}}\\
&=& 4\int_{B(\theta)} \left( s(x;\theta
+{{t}\over{\sqrt{n}}})-s(x;\theta)-{1\over
2}{{t^T}\over{\sqrt{n}}} \del s(x;\theta) \right)^{2}
d\mu(x)\\
&=&
o({{1}\over{n}}).
\end{eqnarray*}
Finally we prove { 4.} We have
\begin{eqnarray*}
\lefteqn{P_{\theta} \left( \left\{ \left| T_{n} \right| \geq
\varepsilon \right\}
\right)}\\
&\leq& P_{\theta} \left( \left\{\left|
T_{n}-{{t^T}\over{\sqrt{n}}}\dot{\ell}(\theta) \right|
\geq{{\varepsilon}\over{2}} \right\} \right) + P_{\theta} \left(
\left\{\left| {{t^T}\over{\sqrt{n}}}\dot{\ell}(\theta) \right|
\geq{{\varepsilon}\over{2}} \right\}
\right)\\
&\leq& {{4}\over{\varepsilon^{2}}} {\rm E}_{\theta} \left(
T_{n}-{{t^T}\over{\sqrt{n}}}\dot{\ell}(\theta) \right)^{2} +
{{4}\over{\varepsilon^{2}}} {\rm E}_{\theta} \left(
{{t^T}\over{\sqrt{n}}}\dot{\ell}(\theta) \right)^{2}
1_{\{|{{t^T}\over{\sqrt{n}}}\dot{\ell}(\theta)|
\geq{{\varepsilon}\over{2}}\}}\\
&=& o({{1}\over{n}})+ {{4}\over{\varepsilon^{2}n}} {\rm
E}_{\theta}(t^T\dot{\ell}(\theta))^2 1_{\{|t^T\dot{\ell}(\theta)|
\geq{{\sqrt{n}\varepsilon}\over{2}}\}}\\
&=&
o({{1}\over{n}}).
\end{eqnarray*}
The last equality holds since $t^T\dot{\ell}(\theta)$ is
$P_{\theta}$-a.e. finite. So $(t^T\dot{\ell}(\theta))^2
1_{\{|t^T\dot{\ell}(\theta)|
\geq{{\sqrt{n}\varepsilon}\over{2}}\}} $ vanishes $P_{\theta}$
almost surely and by the dominated convergence theorem  we get
\[
{\rm E}_{\theta}\dot{\ell}^{2}(\theta)
1_{\{|t^T\dot{\ell}(\theta)|
\geq{{\sqrt{n}\varepsilon}\over{2}}\}} \longrightarrow 0.
\]
This completes the proof. 
\hfill$\Box$

\newpage

We now prove the claims in the proof of Theorem \ref{thm:3.1}
\begin{enumerate}
\item
$P_{\theta}\left(\{\max_{1\leq k\leq
n}|T_{nk}|>\varepsilon\}\right) = o(1)$
\item
$\sum_{k=1}^{n}T_{nk}-
\left(t^TS_{n}(\theta)-{{1}\over{4}}t^TI(\theta)t\right)=o_{P_\theta}(1)$
\item
$\sum_{k=1}^{n}T_{nk}^{2}-t^{T}I(\theta)t=o_{P_\theta}(1)$
\item
$\sum_{k=1}^{n}\alpha_{nk}|T_{nk}|^{3}=o_{P_\theta}(1)$
\end{enumerate}
We have
$$
P_{\theta}\left(\{\max_{1\leq k\leq
n}|T_{nk}|>\varepsilon\}\right) \leq
\sum_{k=1}^{n}P_{\theta}\left(\{|T_{nk}|>\varepsilon\}\right) =
nP_{\theta}\left(\{|T_{n1}|>\varepsilon\}\right) = o(1)
$$
by the previous lemma. This proves 1.

\noindent Next we prove 2. We have
\begin{eqnarray*}
\lefteqn{ {\rm E}_{\theta} \left( \sum_{k=1}^{n}T_{nk}-
\left(t^TS_{n}(\theta)-{{1}\over{4}}t^{T}I(\theta)t\right)
\right)^{2}}\\
&=& {\rm E}_{\theta} \left(
\sum_{k=1}^{n}\left(T_{nk}-{{t^T}\over{\sqrt{n}}}\dot{\ell}(X_{k};\theta)
+{{1}\over{4n}}t^{T}I(\theta)t\right)
\right)^{2}\\
&=& \sum_{k=1}^{n}\sum_{j=1}^{n}{\rm E}_{\theta} \left(
T_{nk}-{{t^T}\over{\sqrt{n}}}\dot{\ell}(X_{k};\theta)
+{{1}\over{4n}}t^{T}I(\theta)t
\right)\\
&&\quad\quad\times\left(
T_{nj}-{{t^T}\over{\sqrt{n}}}\dot{\ell}(X_{j};\theta)
+{{1}\over{4n}}t^{T}I(\theta)t
\right)\\
&\stackrel{(a)}{=}& n{\rm E}_{\theta} \left(
T_{n1}-{{t^T}\over{\sqrt{n}}}\dot{\ell}(X_{1};\theta)
+{{1}\over{4n}}t^{T}I(\theta)t
\right)^{2}\\
&& + n(n-1) {\rm E}_{\theta} \left(
T_{n1}-{{t^T}\over{\sqrt{n}}}\dot{\ell}(X_{1};\theta)
+{{1}\over{4n}}t^{T}I(\theta)t
\right)\\
&&\quad\quad\times {\rm E}_{\theta} \left(
T_{n2}-{{t^T}\over{\sqrt{n}}}\dot{\ell}(X_{2};\theta)
+{{1}\over{4n}}t^{T}I(\theta)t
\right)\\
&\stackrel{(b)}{=}& n{\rm E}_{\theta} \left(
T_{n1}-{{t^T}\over{\sqrt{n}}}\dot{\ell}(X_{1};\theta) \right)^{2}
+{{1}\over{2}}t^{T}I(\theta)t{\rm E}_{\theta} \left(
T_{n1}-{{t^T}\over{\sqrt{n}}}\dot{\ell}(X_{1};\theta)
\right)\\
&&\quad+{{1}\over{16n}}(t^{T}I(\theta)t)^2
+o(1)\\
&\stackrel{(c)}{=}&
o(1).
\end{eqnarray*}
(a) because the terms are i.i.d.\\
(b) and (c) because of the previous lemma and lemma \ref{expscore}\\

\noindent By the Chebyshev inequality we have
\[
P_{\theta} \left( \left| \sum_{k=1}^{n}T_{nk} -t^TS_{n}(\theta) +
{{1}\over{4}}t^{T}I(\theta)t \right|\geq\varepsilon \right)
\longrightarrow 0\;\;\;\;\;\forall\varepsilon >0.
\]
Next we prove 3. Note that
\begin{eqnarray*}
\lefteqn{\sum_{k=1}^{n}T_{nk}^{2} - t^{T}I(\theta)t =
\sum_{k=1}^{n}T_{nk}^{2}
-{{1}\over{n}}\sum_{k=1}^{n}t^{T}\dot{\ell}\dot{\ell}^T(X_{k};\theta)t}\\
&+&
{{1}\over{n}}\sum_{k=1}^{n}t^{T}\dot{\ell}\dot{\ell}^T(X_{k};\theta)t-
t^{T}I(\theta)t.
\end{eqnarray*}
By the weak law of large numbers
\[
{{1}\over{n}}\sum_{k=1}^{n}t^{T}\dot{\ell}\dot{\ell}^T(X_{k};\theta)t
\stackrel{P_\theta}{\longrightarrow} t^{T}I(\theta)t
\]
and by the previous lemma and the Cauchy-Schwarz inequality
\begin{eqnarray*}
\lefteqn{ {\rm E}_{\theta} \left| \sum_{k=1}^{n}T_{nk}^{2}-
{{1}\over{n}}\sum_{k=1}^{n}t^{T}\dot{\ell}\dot{\ell}^{T}(X_{k};\theta)t
\right|}\\
&\leq& n{\rm E}_{\theta} \left|
T_{n1}^{2}-\left({{t^T}\over{\sqrt{n}}}\dot{\ell}(X_{1};\theta)\right)^{2}
\right|\\
&\leq& n{\rm E}_{\theta} \left(
T_{n1}-{{t^T}\over{\sqrt{n}}}\dot{\ell}(X_{1};\theta) \right)^{2}
+ 2n{\rm E}_{\theta}
\left|{{t^T}\over{\sqrt{n}}}\dot{\ell}(X_{1};\theta)\right|
\left|T_{n1}-{{t^T}\over{\sqrt{n}}}\dot{\ell}(X_{1};\theta)\right|\\
&\leq& o(1)+ 2n\left( {\rm E}_{\theta}
({t^{T}\over{n}}\dot\ell(X_{1};\theta))^2 \right)^{{1}\over{2}}
\left( {\rm E}_{\theta}
\Big(T_{n1}-{{t^T}\over{\sqrt{n}}}\dot{\ell}(X_{1};\theta)\Big)^{2}
\right)^{{1}\over{2}}\\
&=&
o(1).
\end{eqnarray*}
Here we have used the inequality $|a^2-b^2|\leq (a-b)^2 + 2|b(a-b)|$.

\noindent So for all $\varepsilon >0$
\begin{eqnarray*}
\lefteqn{ P_{\theta} \left( \left| \sum_{k=1}^{n}T_{nk}^{2}
-t^{T}I(\theta)t \right|\geq\varepsilon
\right) }\\
&\leq& P_{\theta} \left( \left| \sum_{k=1}^{n}T_{nk}^{2}-
{{1}\over{n}}\sum_{k=1}^{n}(t^{T}\dot{\ell}(X_{k};\theta))^2
\right|\geq{{\varepsilon}\over{2}}
\right)\\
&&\quad\quad + P_{\theta} \left( \left|
{{1}\over{n}}\sum_{k=1}^{n}(t^{T}\dot{\ell}(X_{k};\theta))^2
-t^{T}I(\theta)t \right|\geq{{\varepsilon}\over{2}}
\right)\\
&\leq& {{2}\over{\varepsilon}}\, {\rm E}_{\theta} \left|
\sum_{k=1}^{n}T_{nk}^{2}-
{{1}\over{n}}\sum_{k=1}^{n}(t^{T}\dot{\ell}(X_{k};\theta))^2
\right|\\
&&\quad\quad + P_{\theta} \left( \left|
{{1}\over{n}}\sum_{k=1}^{n}t^{T}\dot{\ell}\dot{\ell}^T(X_{k};\theta)t
-t^{T}I(\theta)t \right|\geq{{\varepsilon}\over{2}}
\right)\\
&=&
o(1)
\end{eqnarray*}
It remains to prove 4. Since $|\alpha_{nk}|\leq 1$
we have
\begin{eqnarray*}
\lefteqn{ P_{\theta} \left( \left|
\sum_{k=1}^{n}\alpha_{nk}|T_{nk}|^{3} \right|\geq\varepsilon
\right) }\\
&\leq& P_{\theta} \left( \sum_{k=1}^{n}|T_{nk}|^{3}
\geq\varepsilon
\right)\\
&\leq& P_{\theta} \left( \max_{1\leq k\leq
n}|T_{nk}|\sum_{k=1}^{n}|T_{nk}|^{2}\geq\varepsilon
\right)\\
&=& P_{\theta} \left( \max_{1\leq k\leq
n}|T_{nk}|\sum_{k=1}^{n}|T_{nk}|^{2}\geq\varepsilon,
\sum_{k=1}^{n}|T_{nk}|^{2}\leq 1+t^{T}I(\theta)t
\right)\\
&& \quad\quad + P_{\theta} \left( \max_{1\leq k\leq
n}|T_{nk}|\sum_{k=1}^{n}|T_{nk}|^{2}\geq\varepsilon,
\sum_{k=1}^{n}|T_{nk}|^{2}> 1+t^{T}I(\theta)t\right)\\
&\leq& P_{\theta} \left( \max_{1\leq k\leq n}|T_{nk}|
\geq{{\varepsilon}\over{1+t^{T}I(\theta)t}} \right) + P_{\theta}
\left( \sum_{k=1}^{n}|T_{nk}|^{2}> 1+t^{T}I(\theta)t
\right)\\
&=&
o(1),
\end{eqnarray*}
by statements 1. and 3. of this proof. This proves 4. \hfill$\Box$
\end{Note}

\newpage

\section{Notes Chapter \ref{chap:4}}\label{noteschap4}

\begin{Note}\label{note:6.5}

In a regular parametric model, (\ref{3.1}) yields

\begin{eqnarray*}
\lefteqn{\int|s^2(\theta)-s^2(\theta_0)
-2(\theta-\theta_0)^T\dot\ell(\theta_0)s^2(\theta_0)|d\mu}\\
&\leq& \int |s(\theta)-s(\theta_0) -
(\theta-\theta_0)^T\dot\ell(\theta_0)s(\theta_0)|
|s(\theta)+s(\theta_0)|d\mu\\
&&\quad\quad
+\int|(\theta-\theta_0)^T\dot\ell(\theta_0)s(\theta_0)||s(\theta)-s(\theta_0)|d\mu\\
&\leq& \|s(\theta)-s(\theta_0) -
(\theta-\theta_0)^T\dot\ell(\theta_0)s(\theta_0)\|_\mu
\ \{\|s(\theta)\|_\mu+\|s(\theta_0)\|_\mu\}\\
&&\quad\quad +|\theta-\theta_0|\
\||\dot\ell(\theta_0)|s(\theta_0)\|_\mu
\|s(\theta)-s(\theta_0)\|_\mu =o(|\theta-\theta_0|)
\end{eqnarray*}
with $\|\cdot|\|_\mu$ denoting the norm in $\L_2(\mu)$.\\

It follows that for $n$ and $\sigma$ fixed
\begin{eqnarray*}
\lefteqn{\int_{{\mathbb R}^k}\int_{{\cal X}^n}\Big|\Big\{\prod_{i=1}^n p(x_i;\theta+\delta b)\Big\}\Big\{1+b^T\dot\delta\Big\} {\sqrt{n}\over\sigma}\,w_0\Big({{\sqrt{n}(\theta+\delta b-\theta_0)}\over\sigma}\Big) } \\
&& \ -\,\Big\{\prod_{i=1}^n p(x_i;\theta)\Big\}{\sqrt{n}\over\sigma}\,w_0\Big({{\sqrt{n}(\theta-\theta_0)}\over\sigma}\Big) \\
&& \ -\,\Big\{\delta b^T\sum_{i=1}^n \dot\ell(x_i;\theta)+\delta b^T {\sqrt{n}\over\sigma}\,{\dot w_0\over w_0} \,\Big({{\sqrt{n}(\theta-\theta_0)}\over\sigma}\Big) +b^T\dot\delta\Big\} \\
&& \ \times\Big\{\prod_{i=1}^np(x_i;\theta)\Big\} {\sqrt{n}\over\sigma}\, w_0\Big({{\sqrt{n}(\theta-\theta_0)}\over\sigma}\Big)\Big|
\ d\mu(x_1)\cdots d\mu(x_n)d\theta \\
&& \leq \sum_{i=1}^n\int_{{\mathbb R}^k}\int_{{\cal X}^n} \Big|p(x_i;\theta+\delta b)-p(x_i;\theta)-\delta b^T\dot\ell(x_i;\theta)p(x_i;\theta)\Big| \\
&& \ \times \Big\{\prod_{j=1}^{i-1}p(x_j;\theta)\prod_{j=i+1}^np(x_j;\theta+\delta b)\Big\} \Big\{1+b^T\dot\delta\Big\}  
{\sqrt{n}\over\sigma} w_0\Big({{\sqrt{n}(\theta+\delta b-\theta_0)}\over\sigma}\Big) d\mu(x_1)\cdots d\mu(x_n)d\theta \\
&& \ +\int_{{\mathbb R}^k}\int_{{\cal X}^n} \Big\{\prod_{i=1}^np(x_i;\theta)\Big\} \Big|\Big\{1+b^T \dot\delta\Big\} {\sqrt{n}\over\sigma}\,
w_0\Big({{\sqrt{n}(\theta+\delta b-\theta_0)}\over\sigma}\Big) \\
&& \ -{\sqrt{n}\over\sigma}\, w_0\Big({{\sqrt{n}(\theta-\theta_0)}\over\sigma}\Big) -\Big({\sqrt{n}\over\sigma}\Big)^2\delta b^T\dot
w_0\Big(\sqrt{n}(\theta-\theta_0)\Big) \\
&& \ -\,b^T\dot\delta{\sqrt{n}\over\sigma}w_0\Big({{\sqrt{n}(\theta-\theta_0)}\over\sigma}\Big)\Big| \ d\mu(x_1)\cdots d\mu(x_n)d\theta \\
&& \ +\sum_{i=1}^n\int_{{\mathbb R}^k}\int_{{\cal X}^n} \Big|\delta b^T\dot\ell(x_i;\theta)p(x_i;\theta)\Big|
\Big\{\prod_{j=1}^{i-1}p(x_j;\theta)\Big\} \\
&& \ \times\Big|\Big\{\prod_{j=i+1}^np(x_j;\theta+\delta b)\Big\} \Big\{1+b^T\dot \delta\Big\}{\sqrt{n}\over\sigma}\,
w_0\Big({{\sqrt{n}(\theta+\delta b-\theta_0)}\over\sigma}\Big) \\
&& \ -\,\Big\{\prod_{j=i+1}^np(x_j;\theta)\Big\} {\sqrt{n}\over\sigma}\, w_0\Big({{\sqrt{n}(\theta-\theta_0)}\over\sigma}\Big) \Big|
\ d\mu(x_1)\cdots d\mu(x_n)d\theta \ \to \ 0,
\end{eqnarray*}
as $\delta$ and $\dot\delta$ converge to 0. 
Together with formulas
(\ref{new2.23}) and (\ref{new2.24}) this proves the first sentence
of the proof of Theorem \ref{thm:new4.1}. Note that for given
$\sigma$, the sample size $n$ has to be large enough such that
$\sqrt{n}\sigma^{-1}w_0(\sqrt{n}\sigma^{-1}(\theta-\theta_0))$
puts all its mass within $\Theta$. This can be done since $w_0$
has bounded support and $\Theta$ is open.

\end{Note}

\newpage

\section{Notes Chapter \ref{chap:7}}\label{noteschap7}

\begin{Note}\label{note:7}

Let the random vector $U$ have a multivariate normal distribution with mean
$0$ and covariance matrix $\Sigma$. Then its density is equal to
$$
f_U(u)=(2\pi)^{-{1\over 2}k}|\Sigma|^{-{1\over 2}}e^{-{1\over 2}
u^T\Sigma^{-1} u}, \quad u\in {\mathbb R}^k.
$$
Now, if $X$, given $U=u$, is multivariate normally distributed with mean vector
$Iu$ and covariance matrix $I$, then its density is equal to
$$
f_{X|U=u}(x)=(2\pi)^{-{1\over 2}k}|I|^{-{1\over 2}}e^{-{1\over 2}
(x-Iu)^TI^{-1} (x-Iu)}, \quad x,u\in {\mathbb R}^k,
$$
The unconditional density of $X$ can be computed from these two distributions
\begin{eqnarray*}
\lefteqn{f_{X,U}(x,u) =
f_U(u) {f_{X,U}(x,u)\over f_U(u)}= f_U(u)f_{X|U=u}(x)}\\
&=&(2\pi)^{-k}|I|^{-{1\over 2}}|\Sigma|^{-{1\over 2}}
e^{-{1\over 2}[ (x-Iu)^TI^{-1} (x-Iu) + u^T\Sigma^{-1}u]}\\
&=&(2\pi)^{-k}|I|^{-{1\over 2}}|\Sigma|^{-{1\over 2}}
e^{-{1\over 2}[ x^TI^{-1}x -2u^Tx+ u^T(\Sigma^{-1}+I)u]}\\
&=&(2\pi)^{-k}|I|^{-{1\over 2}}|\Sigma|^{-{1\over 2} }
e^{-{1\over 2}{x\choose u}^TA{x\choose u}},
\end{eqnarray*}
where we write
$$
A=
\begin{pmatrix}
   I^{-1} & -{\bf 1} \\
 -{\bf 1} & \Sigma^{-1}+I
\end{pmatrix}
=
\begin{pmatrix}
  I \Sigma I + I & I \Sigma \\
  \Sigma I       & \Sigma
\end{pmatrix}
$$
with ${\bf 1}$ denoting the $k\times k$ identity matrix. This proves (\ref{4.1})

To prove (\ref{4.2}), since the distributions are multivariate normal,
 we only have  to check the expectation and
covariances of $X$ and $CX-U$. Clearly the expectations vanish.
We get
\begin{eqnarray*}
\lefteqn{\ex X(CX-U)^T = \ex X(CX)^T- \ex XU^T } \\
&& = \ex (XX^T)C^T- \ex XU^T = (I\Sigma I+I)C^T - I\Sigma = {\bf 0},
\end{eqnarray*}
and (all matrices are symmetric and so matrix multiplication is commutative)
\begin{eqnarray*}
\lefteqn{\ex (CX-U)(CX-U)^T = \ex (CX)(CX)^T-\ex (CX)U^T-\ex U(CX)^T +\ex UU^T } \\
&& = C\ex XX^T C^T -2C\ex XU^T +\ex UU^T = CI\Sigma -2CI\Sigma +\Sigma \\
&& = C(-I\Sigma + (I\Sigma I+I)I^{-1}) = C.
\end{eqnarray*}
Together with $\ex XX^T = I\Sigma I+I$ this proves (\ref{4.2}).
\end{Note}

\backmatter

\end{document}